\documentclass[amsmath,10pt,a4paper]{article}
\usepackage{amsfonts}
\usepackage{amssymb}
\newtheorem{thm}{Th\'eor\`eme}
\newtheorem{lem}[thm]{Lemme}
\newtheorem{prop}[thm]{Proposition}
\newtheorem{cor}[thm]{Corollaire}
\newtheorem{fait}[thm]{Fait}

\title{Tentative d'\'epuisement de la cohomologie d'une vari\'et\'e de Shimura
par restriction \`a ses sous-vari\'et\'es}

\author{N. Bergeron}

\date{}

\begin{document}

\maketitle

\begin{footnotesize}
\begin{center}
{\bf Abstract}
\end{center}
Let $G$ be a connected semisimple group over ${\Bbb Q}$. Given a maximal compact subgroup $K\subset G({\Bbb R})$ such that $X=G({\Bbb R})/K$ is a Hermitian symmetric domain, and a convenient arithmetic subgroup $\Gamma\subset G({\Bbb Q})$, one constructs a (connected) Shimura variety $S=S(\Gamma)=\Gamma\backslash X$.
If $H\subset G$ is a connected semisimple subgroup such that $H({\Bbb R})\cap K$ is maximal compact, then $Y=H({\Bbb R})/K$ is a Hermitian symmetric subdomain of $X$. For each $g\in G({\Bbb Q})$ one can construct a connected Shimura variety $S(H,g)=(H({\Bbb Q})\cap g^{-1}\Gamma g)\backslash Y$ and a natural holomorphic map $j_g \colon S(H,g)\rightarrow S$ induced by the map $H({\Bbb A})\rightarrow G({\Bbb A}), h\mapsto gh$. Let us assume that $G$ is anisotropic, which implies that $S$ and $S(H,g)$ are compact. Then, for each positive integer $k$, the map $j_g$ induces a restriction map
$$R_g \colon H^{k}(S , {\Bbb C})\rightarrow H^{k}(S(H,g) , {\Bbb C}).$$

In this paper we focus on classical Hermitian domains and give explicit criterions for the injectivity of the product of the maps $R_g$ (for $g$ running through $G({\Bbb Q})$) when restricted to the strongly primitive (in the sense of Vogan and
Zuckerman) part of the cohomology. In the holomorphic case we recover previous results of Clozel and Venkataramana \cite{ClozelVenky}.
We also derive applications of our results to the proofs of new cases of the Hodge conjecture and of new results on the
vanishing of the cohomology of some particular Shimura variety.

The symmetric space $X$ is Hermitian and the quotient $S$ is a compact Kaehler manifold. We relate for each of the
different classical Hermitian domains the Vogan-Zuckerman's decomposition of the cohomology to some natural decomposition which is defined by making use of the Chern classes of some special bundles over $S$, namely the bundles obtained as quotients of some universal bundles over $X$ by the action of $\Gamma$.
The method used here is mainly representation-theoretic, applying Matsushima's formula. We are then reduce to
linear algebra and combinatorics. Some generalizations of the Littlewood-Richardson coefficients are defined and
used.

The proof of the criterion mentioned above is then a consequence of a recent result of Venkataramana \cite{Venky2}.
\end{footnotesize}

\subsection*{Introduction}

\subsubsection*{Les vari\'et\'es de Shimura}

Dans tout le texte on d\'esignera par $G$ un groupe alg\'ebrique r\'eductif, connexe et anisotrope sur ${\Bbb Q}$. Les ad\`eles ${\Bbb A}$ de ${\Bbb Q}$
forment un anneau localement compact, dans lequel ${\Bbb Q}$ se plonge diagonalement comme un sous-anneau.
On peut consid\'erer le groupe $G({\Bbb A})$ des points ad\`eliques de $G$, qui contient $G ({\Bbb Q})$
comme sous-groupe discret.

Nous supposerons, pour simplifier et toujours dans tout le texte, que le groupe r\'eductif $G$ est presque simple sur ${\Bbb Q}$ modulo son
centre. Autrement dit, il n'a pas de sous-groupe distingu\'e, non central et connexe d\'efini sur ${\Bbb Q}$.
Il d\'ecoule de cette hypoth\`ese que tous les facteurs simples de l'alg\`ebre de Lie complexe $\mathfrak{g}$ de
$G$ (modulo son centre) sont isomorphes. Nous supposerons de plus que le groupe $G({\Bbb R})$ des points r\'eels est le
produit (avec intersection finie) d'un groupe compact et d'un groupe r\'eel non compact qui est {\it presque simple}
modulo son centre que l'on suppose compact. Nous noterons ce dernier groupe $G^{{\rm nc}}$ (nc signifie ici non compact).

Un {\it sous-groupe de congruence} de $G({\Bbb Q})$ est un sous-groupe de la forme $\Gamma = G({\Bbb Q}) \cap K_f$
o\`u $K_f$ est un sous-groupe compact ouvert du groupe $G({\Bbb A}_f)$ des points ad\`eliques finis de $G$.
Soit $X_G = G({\Bbb R}) / K_{\infty}$ l'espace sym\'etrique associ\'e au groupe $G$, o\`u $K_{\infty} \subset G({\Bbb R})$
est un sous-groupe compact maximal. Nous supposerons que cet espace est hermitien.

Dans cet article on \'etudie les quotients (compacts, puisque $G$ est anisotrope) $\Gamma \backslash X_G$,
o\`u $\Gamma \subset G({\Bbb Q})$ est un sous-groupe de congruence; ces quotients s'identifient aux composantes
connexes de $G({\Bbb Q}) \backslash G({\Bbb A}) / K_{\infty} . K_f = G({\Bbb Q}) \backslash (X_G \times G({\Bbb A}_f ) )/K_f $,
o\`u $K_f \subset G({\Bbb A}_f )$ est un sous-groupe compact ouvert.
Pr\'ecisemment, d\'esignons par $G_f$ l'adh\'erence de $G({\Bbb Q})$ dans le groupe $G({\Bbb A}_f )$. Un sous-groupe
de congruence $\Gamma \subset G({\Bbb Q})$ s'\'ecrit $\Gamma = G({\Bbb Q}) \cap K_f$ o\`u $K_f$ est l'adh\'erence de
$\Gamma$ dans $G_f$ et,
\begin{eqnarray} \label{quotient}
\Gamma \backslash X_G = G({\Bbb Q}) \backslash (X_G \times G_f ) / K_f .
\end{eqnarray}
Plus exactement, on s'interesse ici \`a la cohomologie
\`a coefficients complexes \footnote{Tous les groupes de cohomologie que nous consid\`ererons seront
\`a coefficients complexes.} $H^* (\Gamma \backslash X_G )$ de ces quotients.

Une fois donn\'e deux sous-groupes de congruence $\Gamma ' \subset \Gamma \subset G({\Bbb Q})$, on obtient un
rev\^etement fini
$$\Gamma ' \backslash X_G \rightarrow \Gamma \backslash X_G$$
qui induit un morphisme injectif
$$H^* (\Gamma \backslash X_G ) \rightarrow H^* (\Gamma ' \backslash X_G )$$
en cohomologie. Les groupes de cohomologies $H^* (\Gamma \backslash X_G )$ (ou
$H^* (G({\Bbb Q}) \backslash (X_G \times G_f ) / K_f )$) forment donc un syst\`eme inductif index\'e par les
sous-groupes de congruence $\Gamma \subset G({\Bbb Q})$ (ou par les sous-groupes compacts ouverts $K_f \subset
G_f$).  En passant \`a la limite (inductive) on d\'efinit
\begin{eqnarray} \label{HSh}
H^* (Sh^0 G) = \lim_{
\begin{array}{c}
\rightarrow \\
\Gamma
\end{array}} H^* (\Gamma \backslash X_G )
= \lim_{
\begin{array}{c}
\rightarrow \\
K_f
\end{array}} H^* (G({\Bbb Q}) \backslash (X_G \times G_f ) / K_f ).
\end{eqnarray}
La notation ci-dessus provient de ce que l'on appelle {\it vari\'et\'e de Shimura}
l'espace topologique
\begin{eqnarray} \label{Sh}
Sh^0 G = \lim_{
\begin{array}{c}
\leftarrow \\
\Gamma
\end{array}} \Gamma \backslash X_G = G({\Bbb Q}) \backslash (X \times G_f ).
\end{eqnarray}
Cet espace est un espace
topologique dont on peut consid\'erer la cohomologie ce C\v{e}ch et il est d\'emontr\'e dans \cite{Rohlfs} que sa
cohomologie co\"{\i}ncide avec (\ref{HSh}).
Pour ce qui nous concerne, il sera suffisant de consid\'erer que la d\'enomination $H^* (Sh^0 G)$ n'est qu'une
notation pour la limite inductive (\ref{HSh}).

Les quotients $\Gamma \backslash X_G$ sont kaehl\'eriens. Dans chacun des groupes
$H^i (\Gamma \backslash X_G )$ pour $0 \leq i \leq \frac{d_G}{2}$, o\`u $d_G$ est la dimension r\'eelle de $X_G$,
on peut consid\'erer la partie primitive de la cohomologie. En passant \`a la limite inductive cela nous d\'efinit la
partie primitive des espaces de cohomologies
$$H_{{\rm prim}}^i (Sh^0 G) , \; \; (0 \leq i \leq \frac{d_G}{2} )$$
de la vari\'et\'e de Shimura associ\'ee.

\paragraph*{Une question d'Arthur}

\`A la toute fin de son c\'el\`ebre article \cite{Arthur}, Arthur pose la question de la possibilit\'e de d\'ecouper les espaces
$H_{{\rm prim}}^i (Sh^0 G)$, $(0 \leq i \leq \frac{d_G}{2} )$, en morceaux identifi\'es \`a des sous-espaces de la cohomologie
primitive en degr\'e m\'edian $H_{{\rm prim}}^{d_H /2} (Sh^0 H)$, attach\'ee \`a des vari\'et\'es de Shimura $Sh^0 H$ de
dimension $d_H$ plus petite que $d_G$.

Cette belle question motive une grande part de cette article. On peut penser \`a cette question comme \`a un
analogue du th\'eor\`eme de Lefschetz pour les vari\'et\'es projectives. Dans son article Arthur propose des candidats
pour les vari\'et\'es de Shimura $Sh^0 H$ de dimension plus petite. Les groupes $H$ devraient \^etre des groupes
endoscopiques. Il est naturel de se demander si, comme dans le cas du th\'eor\`eme de Lefschetz, on ne pourrait
pas \'epuiser une large part de la cohomologie en se restreignant \`a des sous-vari\'et\'es de Shimura.

\subsubsection*{Restriction \`a une sous-vari\'et\'e de Shimura}

Puisque $X_G$ de type hermitien, il existe un \'el\'ement $c$ appartenant au centre de
$K_{\infty}$ tel que $Ad(c)$ induise (sur l'espace tangent $\mathfrak{p}_0$ de $X_G$ au point base $o = eK$)
la multiplication par $i= \sqrt{-1}$. Soit
\begin{eqnarray} \label{dectriang1}
\mathfrak{g} = \mathfrak{k} \oplus \mathfrak{p}^+ \oplus \mathfrak{p}^-
\end{eqnarray}
la d\'ecomposition associ\'ee de $\mathfrak{g} = \mbox{Lie}( G({\Bbb R})) \otimes {\Bbb C}$, avec $\mathfrak{k} =
\mbox{Lie} ( K_{\infty} ) \otimes {\Bbb C}$. Alors
$\mathfrak{p}^+ = \{ X \in \mathfrak{p} \; :  \; Ad (c) X= iX \}$ est l'espace tangent holomorphe \`a $X_G$
au point base $o$.

Soit maintenant $H \subset G$ un sous-groupe r\'eductif connexe d\'efini sur ${\Bbb Q}$. On suppose que
\begin{eqnarray} \label{cond1}
H({\Bbb R})  \cap K_{\infty} \mbox{ {\it est un sous-groupe compact maximal de }} H({\Bbb R}).
\end{eqnarray}
Alors la restriction \`a $H$ de l'involution de Cartan $\theta$ de $G$ est une involution de Cartan de $H({\Bbb R})$.
On a une d\'ecomposition correspondante
\begin{eqnarray}
\mathfrak{h} =  \mathfrak{k}_H \oplus \mathfrak{p}_H
\end{eqnarray}
avec $\mathfrak{p}_H = \mathfrak{p} \cap \mathfrak{h}$. On suppose de plus
\begin{eqnarray} \label{cond2}
\mathfrak{p}_H \mbox{ {\it est stable sous l'action de }} Ad(c).
\end{eqnarray}
Alors $Ad(c)$ d\'efinit une structure complexe $(H({\Bbb R}) \cap K_{\infty})$-invariante sur $\mathfrak{p}_{H,0}$.
L'espace $X_H = H({\Bbb R}) /(H({\Bbb R}) \cap K_{\infty})$ est symm\'etrique hermitien. On a une d\'ecomposition
triangulaire
\begin{eqnarray}
\mathfrak{h} = \mathfrak{k}_H \oplus \mathfrak{p}_H^+ \oplus \mathfrak{p}_H^-
\end{eqnarray}
compatible avec (\ref{dectriang1}). Enfin, le plongement $X_H \hookrightarrow X_G$ est holomorphe.

Consid\'erons maintenant $\Gamma = G({\Bbb Q}) \cap K_f$ un sous-groupe de congruence sans torsion de $G({\Bbb Q})$.
Le quotient $S(\Gamma ) = \Gamma \backslash X_G$ est une vari\'et\'e kaehl\'erienne qui s'identifie \`a
$S(K_f )= G({\Bbb Q}) \backslash (X \times G_f ) / K_f$.

Soit $K_f^H \subset H({\Bbb A}_f )$ un sous-groupe compact ouvert. Si $K_f^H \subset K_f$, il existe une application
naturelle $j : S(K_f^H ) \rightarrow S(K_f )$. Puisque $\Gamma$ est sans torsion, l'application $j$ est finie et non ramifi\'ee.
On rappelle le r\'esultat suivant, d\^u \`a Deligne \cite{Deligne}.

\begin{lem}
\'Etant donn\'e $K_f^H \subset H({\Bbb A}_f )$, il existe un sous-groupe compact ouvert $K_f^1 \subset G({\Bbb A}_f )$
avec $K_f^H \subset K_f^1$, tel que l'application naturelle $j' : S(K_f^H ) \rightarrow S(K_f^1 )$ soit injective.
\end{lem}

En particulier, si l'on prend $K_f^H = K_f \cap H({\Bbb A}_f )$, on obtient une application naturelle finie $j : S(K_f^H )
\rightarrow S(K_f )$. Si l'on remplace $K_f$ par un sous-groupe suffisamment petit $K_f^1$, on obtient un diagramme
\begin{eqnarray}
\begin{array}{ccc}
S(K_f^H ) & \stackrel{j'}{\rightarrow} & S(K_f^1 ) \\
               & \stackrel{j}{\searrow}      & \downarrow \pi \\
               &                                         & S(K_f ) ,
\end{array}
\end{eqnarray}
o\`u $\pi$ est la projection naturelle de rev\^etement et $j'$ est injective.

En passant \`a la limite (inductive) sur les $K_f$, les applications $j$ induisent  l'application de restriction
\begin{eqnarray} \label{res}
\mbox{res}_H^G : H^* (Sh^0 G) \rightarrow H^* (Sh^0 H) .
\end{eqnarray}
Nous dirons que $Sh^0 H$ est une {\it sous-vari\'et\'e de Shimura de $Sh^0 G$}, ce que l'on notera
$Sh^0 H \subset Sh^0 G$.

Nous aurons besoin de consid\'erer \'egalement la stabilisation de cette application de restriction. Expliquons ce
que cela signifie. Soit $g\in G({\Bbb Q})$. Fixons $K_f^H$ un compact ouvert de $H({\Bbb A}_f )$, et consid\'erons
l'application $j_g : X_H \times H_f \rightarrow X_G \times G_f$ donn\'ee par $j_g (x,h) = (gx,gh)$. Il est facile de
v\'erifier que $j_g$ induit une application injective $H({\Bbb Q}) \backslash (X_H \times H_f ) /K_f^H \rightarrow
G({\Bbb Q}) \backslash (X_G \times G_f ) /K_f^H$. En supposant que $K_f^H = K_f \cap H({\Bbb A}) \subset K_f$ (o\`u
$K_f$ est un sous-groupe compact ouvert dans $G({\Bbb A}_f )$), on obtient alors une application naturelle
$j_g : S(K_f^H ) \rightarrow S(K_f )$, en utilisant les notations pr\'ec\'edentes. Cette application est finie et non ramifi\'ee.
On peut \'egalement d\'ecrire $j_g$ comme l'application naturelle $(H({\Bbb R}) \cap g^{-1} \Gamma g) \backslash X_H
\rightarrow  \Gamma \backslash X_G$. On obtient de cette mani\`ere toute une famille, param\`etr\'ee par $g\in G({\Bbb Q})$,
de sous-vari\'et\'es complexes de $\Gamma \backslash X_G$ - les images des applications $j_g$. En cohomologie celles-ci
induisent l'application de {\it restriction stable}
$$H^* (S(\Gamma )) \rightarrow \prod_{g\in G({\Bbb Q})} H^* (S_H (g)),$$
o\`u $S_H (g) = (H({\Bbb R}) \cap g^{-1} \Gamma g) \backslash X_H$, et l'application de restriction est d\'eduite de la famille
d'applications $(j_g )$.
En passant \`a la limite (inductive) sur les $\Gamma$, l'application de restriction stable induit l'application Res$_H^G$
{\it stabilisation} de res$_H^G$~:
\begin{eqnarray} \label{Res}
\mbox{Res}_H^G : H^* (Sh^0 G) \rightarrow \prod_{g\in G({\Bbb Q})} H^* (Sh^0 H).
\end{eqnarray}

\paragraph*{R\'esultats de Venkataramana concernant la restriction et la question d'Arthur}

Dans \cite{Venky2}, Venkataramana d\'emontre le th\'eor\`eme suivant \footnote{Remarquons que le cas de la cohomologie
holomorphe \'etait d\'ej\`a trait\'e dans des travaux ant\'erieurs de Oda \cite{Oda} et de Clozel et Venkataramana
\cite{ClozelVenky}. Enfin, l'\'enonc\'e g\'en\'eral de ce th\'eor\`eme \'etait conjectur\'e par Harris et Li dans \cite{HarrisLi}
o\`u des cas particuliers sont d\'emontr\'es.}.

\begin{thm} \label{venky}
Supposons que $H \subset G$ soient deux groupes alg\'ebriques r\'eductifs connexes et anisotropes sur ${\Bbb Q}$ avec $G^{{\rm nc}} =U(n,1)$ (resp. $G^{{\rm nc}} = O(2,n)$) et
$H^{{\rm nc}} = U(n-1,1)$ (resp. $H^{{\rm nc}} = O(2,n-1)$). Alors,
$$\mbox{Res}_H^G : H^k (Sh^0 G) \rightarrow \prod_{g\in G({\Bbb Q})} H^k (Sh^0 H)$$
est injective pour tout $k\leq n-1$.
\end{thm}

On en d\'eduit imm\'ediatement le corollaire suivant qui r\'epond \`a la question d'Arthur dans le cas des groupes $G$
alg\'ebriques r\'eductifs connexes sur ${\Bbb Q}$ avec $G^{{\rm nc}} = U(n,1)$ ou $G^{{\rm nc}} = O(2,n)$ et provenant, par restriction des
scalaires, d'un groupe unitaire ou orthogonal sur un corps de nombre.

\begin{cor} \label{arthur venky}
Soient $K$ un corps de nombre totalement r\'eel et $G'$ un groupe unitaire (resp. orthogonal) anisotrope sur $K$ compact \`a toutes les places infinies
sauf une o\`u il est isomorphe \`a $U(n,1)$ (resp. $O(2,n)$). Soit $G$ le groupe alg\'ebrique sur ${\Bbb Q}$ obtenu \`a partir de
$G'$ par restriction des scalaires de $K$ \`a ${\Bbb Q}$. Alors, pour tout entier naturel $k\leq n$, il existe une sous-vari\'et\'e
de Shimura $Sh^0 H \subset Sh^0 G$ de dimension complexe $k$ telle que la restriction stable
$$\mbox{Res}_H^G : H^k (Sh^0 G) \rightarrow \prod_{g\in G({\Bbb Q})} H^k (Sh^0 H)$$
soit injective. De plus, l'image de la partie primitive de $H^k (Sh^0 G)$ est envoy\'ee dans la partie primitive de
la cohomologie (de degr\'e m\'edian) de $Sh^0 H$.
\end{cor}

Les groupes $G$ alg\'ebriques r\'eductifs connexes sur ${\Bbb Q}$ avec
$G^{{\rm nc}} \cong U(n,1)$ \`a l'oppos\'ee de ceux consid\'er\'es dans
le corollaire sont les groupes qui ont la propri\'et\'e
(K) de \cite{Clozel} . Il serait int\'eressant de voir si les techniques de \cite{Clozel} jointes au th\'eor\`eme \ref{venky} permettent
de r\'epondre positivement \`a la question d'Arthur en toute g\'en\'eralit\'e pour les groupe unitaires de partie non compacte isomorphe \`a $U(n,1)$.

\bigskip

La d\'emonstration du th\'eor\`eme \ref{venky} repose sur un beau r\'esultat g\'en\'eral \'egalement d\^u \`a Venkataramana \cite{Venky2} donnant
un crit\`ere d'injectivit\'e de l'application de restriction stable en restriction \`a certains sous-espaces. Il s'agit ensuite de comprendre ce
crit\`ere ce qui nous ram\`ene essentiellement \`a de l'alg\`ebre lin\'eaire, facile dans le cas des groupes $U(n,1)$ et $O(2,n)$. Le but de ce texte est de
comprendre ce crit\`ere dans le cas des domaines hermitiens classiques. C'est exactement ce que font Clozel et Venkataramana dans \cite{ClozelVenky} dans le
cas de la cohomologie holomorphe. Nous retrouverons donc leurs r\'esultats, notons au passage que nos arguments sont purement locaux
alors que ceux de Clozel et Venkataramana font intervenir des ph\'enom\`enes globaux.

\medskip

Dans ce qui suit, $M_{m \times n} ({\Bbb C})$ d\'esigne l'espace des matrices avec $m$ lignes et $n$ colonnes et des coefficients complexes. Nous noterons
$E_{i,j} \in M_{m \times n} ({\Bbb C})$ la matrice dont le coefficient appartenant \`a la $i$-\`eme ligne et \`a la $j$-\`eme colonne est \'egal \`a $1$ et dont
tous les autres coefficients sont nuls. Le groupe des matrices carr\'ees $p\times p$ inversibles est not\'e $GL_p$. Le groupe des matrices unitaires $p\times p$ est
not\'es $U(p)$.

Si $E$ est une repr\'esentation d'un groupe, nous noterons $E^*$ sa rep\'esentation contragr\'ediente. Si $e_1 , e_2 , \ldots , e_m$ est une base de $E$, alors sa base duale
de $E^*$ est not\'ee $e_1^* , e_2^* , \ldots , e_m^*$. La $r$-i\`eme puissance ext\'erieure (resp. sym\'etrique) de $E$ est d\'esign\'ee $\bigwedge^r E$ (resp. sym$^r (E)$).

Nous allons nous int\'eresser tour \`a tour aux diff\'erents domaines hermitiens classiques en excluant le cas $G^{{\rm nc}} = O(2,n)$ d\'ej\`a trait\'e dans \cite{Venky2}.
Le cas o\`u $G^{{\rm nc}} = U(p,q)$ est le plus int\'eressant. Il nous occupera une bonne partie de l'article.
Remarquons imm\'ediatement que m\^eme lorsque $G$ est un ${\Bbb Q}$-groupe provenant d'un groupe unitaire sur un corps
de nombre on ne peut esp\'erer r\'epondre positivement \`a la question d'Arthur seulement \`a coup de restrictions (stables)
\`a des sous-vari\'et\'es de Shimura. Le premier \`a l'avoir remarqu\'e est Venkataramana dans \cite{Venky3}, on
rencontrera dans le cours du texte des contres-exemples similaires. Une grande partie de cet article est motiv\'ee par
la volont\'e de comprendre quelle partie de la cohomologie des vari\'et\'es de Shimura unitaires peut \^etre \'epuiser \`a
coup de restrictions (stables) \`a des sous-vari\'et\'es de Shimura.

\subsubsection*{\'Enonc\'es des principaux r\'esultats}

Concernant la question d'Arthur, nos r\'esultats nous permettrons de d\'emontrer le th\'eor\`eme suivant qui g\'en\'eralise le
corollaire \ref{arthur venky} dans le cas unitaire.

\begin{thm} \label{tentative u(p,q)}
Soient $K$ un corps de nombre totalement r\'eel et $G'$ un groupe unitaire  anisotrope sur $K$ compact \`a toutes les places infinies
sauf une o\`u il est isomorphe \`a $U(p,q)$ avec $1 \leq p \leq q$. Soit $G$ le groupe alg\'ebrique sur ${\Bbb Q}$ obtenu \`a partir de
$G'$ par restriction des scalaires de $K$ \`a ${\Bbb Q}$. Alors,
pour tout entier naturel $k< 3p-2$ si $p=q$ et $< p+q-1$ si $p<q$, l'espace
$H^k (Sh^0 G)$ se d\'ecompose en une somme directe finie de sous-espaces telle que chacun de ces sous-espaces
s'injecte naturellement dans la cohomologie d'une vari\'et\'e de Shimura de dimension strictement plus petite que $pq$.
\end{thm}

Nous verrons de plus que les degr\'es $3p-2$ et $p+q-1$ sont optimaux si l'on suppose $Sh^0 H \subset Sh^0 G$.

Dans le cas des groupes $GSp_p$ et $O^*(2p)$ nos m\'ethodes n'apportent pas, concernant la question d'Arthur, beaucoup plus que les r\'esultats de Clozel et Venkataramana que  nous retrouverons. Nous verrons au cours du texte, qu'ils permettent n\'eanmoins de pr\'eciser l'action des classes de Chern sur la cohomologie.

\medskip

Revenons maintenant au cas des vari\'et\'es de Shimura unitaires, {\it i.e} $G^{{\rm nc}} = U(p,q)$.
Le th\'eor\`eme de Lefschetz sur les classes de Hodge de bidegr\'e $(1,1)$ est non vide lorsque $p=1$ et implique que
toute classe de Hodge dans $H^2 (Sh^0 G)$ et donc dans $H^{2q-2} (Sh^0 G)$ est alg\'ebrique. Lorsque $q\geq p>1$, il n'y a pas de
classes de Hodge de degr\'e $2$, en fait et si $q>p+1$ nous verrons que toute classe de Hodge non triviale est de degr\'e $\geq 2p$.
En d\'eveloppant une id\'ee de Venkataramana, nous montrerons le nouveau cas suivant de la conjecture de Hodge.

\begin{thm} \label{hodge u(p,q)}
Soient $K$ un corps de nombre totalement r\'eel et $G'$ un groupe unitaire  anisotrope sur $K$ compact \`a toutes les places infinies
sauf une o\`u il est isomorphe \`a $U(p,q)$ avec $3 \leq 2p+1 \leq q$. Soit $G$ le groupe alg\'ebrique sur ${\Bbb Q}$ obtenu \`a partir de
$G'$ par restriction des scalaires de $K$ \`a ${\Bbb Q}$. Alors, toute classe de Hodge dans $H^{2p(q-1)} (Sh^0 G)$ est alg\'ebrique.
\end{thm}

Enfin, nous d\'emontrerons le r\'esultat d'annulation suivant.

\begin{thm} \label{vanish}
Supposons que $Sh^0 G$ contienne une sous-vari\'et\'e de Shimura $Sh^0 H$, avec $H$ obtenu par restriction des scalaires \`a partir d'un groupe $U(D)$ o\`u $D$ est une
alg\`ebre \`a division de degr\'e premier impair sur une extension quadratique imaginaire d'un corps de nombre totalement r\'eel et tel que $H^{{\rm nc}} = U(p,b)$
avec $b<q$.
Alors, $H^{k}_{{\rm prim}} (Sh^0 G) =0$ pour tout entier $1 \leq k < \min (bp , q )$.
\end{thm}

\subsubsection*{Plan de l'article}

La premi\`ere partie de l'article est consacr\'ee au groupe $U(p,q)$. Dans une premi\`ere section, nous \'etudions la cohomologie des vari\'et\'es de Shimura unitaires
\`a l'aide des travaux de Vogan et Zuckerman mais en adoptant une param\'etrisation par des couples de partitions. Cette param\'etrisation permet de comprendre
relativement facilement  l'action des classes de Chern sur la cohomologie. Le r\'esultat principal de cette section est la proposition \ref{clef}; il d\'ecrit compl\`etement
cette action. Nous appliquons \`a la toute fin de la premi\`ere section ces id\'ees \`a la d\'emonstration d'un th\'eor\`eme g\'en\'eralisant le th\'eor\`eme de Lefschetz fort.
Dans le cas holomorphe ce th\'eor\`eme n'est pas nouveau il est l'objet de l'article \cite{Parthasarathy} de Parthasarathy.

Dans un deuxi\`eme section, nous classifions les types de sous-espaces hermitiens sym\'etriques pouvant apparaitre comme sous-vari\'et\'e totalement g\'eod\'esique de
$X_G$. Puis nous calculons la classe de cohomologie correspondante \`a chacun de ces sous-espaces dans le dual compact de $X_G$.

Les r\'esultats des deux premi\`eres sections nous permettent \`a l'aide du th\'eor\`eme principal de l'article \cite{Venky2} de Venkataramana de d\'emontrer des crit\`eres
d'injectivit\'e de l'application de restriction stable. On d\'eduit de ces crit\`eres le th\'eor\`eme \ref{tentative u(p,q)}.

Les deux derni\`eres parties (plus courtes) traitent respectivement des cas du groupe $GSp_p$ et du groupe $O^* (2p)$. Nous y d\'emontrons les r\'esultats correspondant.

Concluons en remarquant que nos m\'ethodes ne s'\'etendent pas au cas non isotropes ({\it i.e.}  $S(\Gamma )$ non compacte). Mais il est naturel de penser que les r\'esultats
devraient eux s'\'etendre (en consid\'erant cette fois la cohomologie $L^2$). Nous le montrons partiellement dans un travail en pr\'eparation. Enfin, nous laissons au lecteur
le soin de traduire nos r\'esultats au cas des syst\`emes de coefficients non triviaux, les d\'emonstrations se transposent imm\'ediatement.

\medskip

\paragraph*{Remerciements}

Je ne saurai trop souligner ce que les r\'esultats ci-dessus doivent aux travaux de Venkataramana. La lecture de ceux de ses articles figurant dans la bibliographie m'a
incit\'e \`a vouloir comprendre la combinatoire de l'application de restriction stable entre vari\'et\'es de Shimura. Enfin, j'ai b\'en\'efici\'e de nombreuses
conversations sur ce sujet avec Laurent Clozel, je l'en remercie.

\newpage

\section{Cas des vari\'et\'es de Shimura unitaires}

Dans toute cette partie $G$ est un groupe alg\'ebrique r\'eductif connexe et anisotrope sur ${\Bbb Q}$ avec
$G^{{\rm nc}} \cong U(p,q)$, o\`u $p$ et $q$ sont des entiers strictement positifs avec $p\leq q$.
Le rang r\'eel de $G$ est donc $p$. On a
\begin{eqnarray} \label{Gnc}
G^{{\rm nc}} = \left\{ g= \left(
\begin{array}{cc}
A & B \\
C & D
\end{array} \right) \; : \; ^t \! \overline{g} \left(
\begin{array}{cc}
1_p & 0 \\
0 & -1_q
\end{array} \right) g = \left(
\begin{array}{cc}
1_p & 0 \\
0 & -1_q
\end{array} \right) \right\} ,
\end{eqnarray}
o\`u $A\in M_{p\times p} ({\Bbb C})$, $B\in M_{p\times q} ({\Bbb C})$, $C\in M_{q\times p} ({\Bbb C})$ et
$D\in M_{q\times q} ({\Bbb C})$. Soit
$$K = \left\{ g =  \left(
\begin{array}{cc}
A & 0 \\
0 & D
\end{array} \right) \in G^{{\rm nc}} \; : \;  A\in U(p) , \ D\in U(q) \right\} .$$
Le complexifi\'e $K_{{\Bbb C}}$ de $K$ est le groupe
$$K_{{\Bbb C}} = \left\{ g =  \left(
\begin{array}{cc}
A & 0 \\
0 & D
\end{array} \right) \in G^{{\rm nc}}_{{\Bbb C}} \; : \;  A\in GL_p , \ D\in GL_q \right\} .$$
L'involution de Cartan $\theta$ est donn\'ee par $x \mapsto - ^t \! \overline{x}$. Soit
$$T = \left\{ g\in K_{{\Bbb C}} \; : \; g = \left(
\begin{array}{cc}
A & 0 \\
0 & D
\end{array} \right) \mbox{ avec } A, \ D \mbox{ matrices diagonales} \right\} .$$

Puisque les facteurs compact de $G({\Bbb R})$ ne nous int\'eresseront pas, nous noterons $\mathfrak{g}_0$ l'alg\`ebre de Lie de
$G^{{\rm nc}}$ et $\mathfrak{k}_0$ l'alg\`ebre de Lie de $K$. Soit $\mathfrak{g}_0 = \mathfrak{k}_0 \oplus \mathfrak{p}_0$
la d\'ecomposition de Cartan associ\'ee. Si $\mathfrak{l}_0$ est une alg\`ebre de Lie, nous noterons $\mathfrak{l} = \mathfrak{l}_0
\otimes {\Bbb C}$ sa complexification. Rappelons que la multiplication par $i =\sqrt{-1}$ induit une d\'ecomposition
$$\mathfrak{p} = \mathfrak{p}^+ \oplus \mathfrak{p}^- .$$

Nous noterons $(x_1 , \ldots , x_p ; y_1 , \ldots , y_q )$ les \'el\'ements de $T$ ou de son alg\`ebre de Lie.
L'alg\`ebre de Lie $\mathfrak{g}$ est bien \'evidemment $M_{(p+q) \times (p+q)} ({\Bbb C})$, et l'on voit ses
\'el\'ements sous forme de blocs comme dans (\ref{Gnc}). On a alors,
$$\mathfrak{p}^+ = \left\{ \left(
\begin{array}{cc}
0 & B \\
0 & 0
\end{array} \right) \mbox{ avec } B \in M_{p\times q} ({\Bbb C} ) \right\}$$
et
$$\mathfrak{p}^- = \left\{ \left(
\begin{array}{cc}
0 & 0 \\
C & 0
\end{array} \right) \mbox{ avec } C \in M_{q \times p} ({\Bbb C} )\right\} .$$
Soit $E= {\Bbb C}^p$ (resp. $F= {\Bbb C}^q$) la repr\'esentation standard de $U(p)$ (resp. $U(q)$). Alors, comme
repr\'esentation de $K_{{\Bbb C}}$, $\mathfrak{p}^+ = E \otimes F^*$.

Soient $(e_1 , \ldots , e_p )$ et $(f_1 , \ldots , f_q )$ les bases canoniques respectives de $E$ et $F$. Choisissons
comme sous-alg\`ebre de Borel $\mathfrak{b}_K$ dans $\mathfrak{k}$ l'alg\`ebre des matrices dans $\mathfrak{k}$,
qui sont triangulaires sup\'erieures sur $E$ et triangulaires inf\'erieures sur $F$ par rapport \`a ces bases.
Alors l'ensemble des racines simples compactes positives
\begin{eqnarray} \label{racines compactes}
\Phi (\mathfrak{b}_K , \mathfrak{t} ) = \{ x_i - x_j \; : \; 1\leq i < j \leq p \} \cup \{ y_j -y_i \; : \;  1 \leq i < j \leq q \} .
\end{eqnarray}
Les racines de $T$ apparaissant dans $\mathfrak{p}^+$ sont les formes lin\'eaires $x_i - y_j$ avec $1 \leq i \leq p$
et $1 \leq j \leq q$.

\subsection{D\'ecomposition ``\`a la Lefschetz'' de la cohomologie}

Puisque $G$ est anisotrope sur ${\Bbb Q}$, un th\'eor\`eme de Borel et Harish-Chandra \cite{BorelHarishChandra}
affirme que si $\Gamma$ est un sous-groupe de congruence de $G$, la vari\'et\'e $S(\Gamma )$ est compacte.

Soit $\Gamma$ un sous-groupe de congruence de $G$.
Soit ${\cal E}^k (S(\Gamma ))$ l'espace des formes diff\'erentielles de degr\'e $k$ sur $S(\Gamma )$. Puisque le
fibr\'e cotangent $T^* S(\Gamma )$ est isomorphe au fibr\'e $\Gamma \backslash G^{{\rm nc}} \times_K \mathfrak{p}^*
\rightarrow \Gamma \backslash G^{{\rm nc}} /K = S(\Gamma )$, qui est associ\'e au $K$-fibr\'e principal $K \rightarrow \Gamma \backslash G^{{\rm nc}}
\rightarrow \Gamma \backslash G^{{\rm nc}} / K$ et la repr\'esentation (r\'eelle) de $K$ dans $\mathfrak{p}^+$
donn\'ee par
$$\left( \left(
\begin{array}{cc}
A & 0 \\
0 & D
\end{array} \right) , \left(
\begin{array}{cc}
0 & B \\
0 & 0
\end{array} \right) \right) \mapsto \left(
\begin{array}{cc}
0 & ABD^{-1} \\
0 & 0
\end{array} \right) .$$
On a donc~:
\begin{eqnarray} \label{Ek}
{\cal E}^k (S(\Gamma )) \simeq (C^{\infty} (\Gamma \backslash G^{{\rm nc}} )\otimes \bigwedge^k \mathfrak{p}^* )
\simeq \mbox{Hom}_K (\bigwedge^k \mathfrak{p} , C^{\infty} (\Gamma \backslash G^{{\rm nc}}) ) \; \; (k \in {\Bbb N} ).
\end{eqnarray}

Notons $\Delta$ le laplacien de Hodge-de Rham sur la vari\'et\'e riemannienne (localement sym\'etrique) $S(\Gamma )$
(o\`u la m\'etrique est d\'eduite de la forme de Killing sur $\mathfrak{g}_0$). L'espace des formes harmoniques de
degr\'e $k$ est donn\'e par
$${\cal H}^k (S(\Gamma )) := \{ \omega \in {\cal E}^k (S(\Gamma )) \; : \; \Delta \omega = 0 \} .$$
La th\'eorie de Hodge fournit un isomorphisme naturel
$$H^* (S(\Gamma )) \simeq {\cal H}^* (S(\Gamma )) .$$
Soit $(\pi , V_{\pi} )$ un $(\mathfrak{g} , K)$-module irr\'eductible. \`A l'aide de (\ref{Ek}) on d\'efinit une application
lin\'eaire
\begin{eqnarray} \label{Tpi}
T_{\pi} : \left\{
\begin{array}{rcl}
\mbox{Hom}_K (\bigwedge^* \mathfrak{p}, \pi ) \otimes \mbox{Hom}_{\mathfrak{g} , K} (\pi , C^{\infty} (\Gamma \backslash G^{{\rm nc}} ))  & \rightarrow & {\cal E}^* (S(\Gamma )) , \\
\psi \otimes \varphi & \mapsto & \varphi \circ \psi .
\end{array} \right.
\end{eqnarray}

Soit $\hat{G}^{{\rm nc}}$ l'ensemble des classes d'\'equivalence des $(\mathfrak{g} , K)$-modules irr\'eductibles
qui sont unitarisables. Rappelons qu'Harish-Chandra a d\'emontr\'e que $\hat{G}^{{\rm nc}}$ s'identifie naturellement
au dual unitaire de $G^{{\rm nc}}$. Soient $U(\mathfrak{g} )$ l'alg\`ebre enveloppante de l'alg\`ebre de Lie complexe
$\mathfrak{g}$, $Z(\mathfrak{g} )$ son centre et $\Omega \in Z(\mathfrak{g} )$ le casimir d\'efinit par la forme de Killing sur $\mathfrak{g}_0$.
On d\'efinit le sous-ensemble $\hat{G}^{{\rm nc}}_0$ de $\hat{G}^{{\rm nc}}$ par
$$\hat{G}^{{\rm nc}}_0 := \{ \pi \in \hat{G}^{{\rm nc}} \; : \;  \pi (\Omega ) =0 \} ,$$
o\`u l'on a conserv\'e la m\^eme notation $\pi$ pour la repr\'esentation de $U(\mathfrak{g} )$.

L'acion de $G^{{\rm nc}}$ sur $X_G$ induit la repr\'esentation de $U(\mathfrak{g} )$ sur l'espace des formes diff\'erentielles
sur $X_G$. En particulier, le casimir $\Omega ( \in Z (\mathfrak{g} ) \subset U(\mathfrak{g} ) )$ agit sur ${\cal E}^* (X_G )$
comme le laplacien de Hodge-de Rham, puisque la m\'etrique riemannienne sur $X_G$ est induite par la forme de Killing sur $\mathfrak{g}_0$.
Il d\'ecoule de tout ceci que
\begin{eqnarray} \label{im Tpi}
\mbox{Image} (T_{\pi}) \subset {\cal H}^* (S(\Gamma )) \simeq H^* (S(\Gamma ))
\end{eqnarray}
si et seulement si $\pi \in \hat{G}^{{\rm nc}}_0$. On dit dans ce cas que le sous-espace de $H^* (S(\Gamma ))$ correspondant
\`a l'image de $T_{\pi}$ est {\it la $\pi$-composante}, et on \'ecrit $H^* ( \pi : \Gamma )$. Autrement dit,
\begin{eqnarray} \label{pi composante}
H^k (\pi : \Gamma ) := \mbox{Image} (T_{\pi}) \cap H^k (S(\Gamma )) \; \;  (k \in {\Bbb N} ) ,
\end{eqnarray}
via l'isomorphisme (\ref{im Tpi}).

Un r\'esultat d\^u \`a Gel'fand et Piatetski-Shapiro \cite{GGPS} affirme que la repr\'esentation r\'eguli\`ere droite
dans $L^2 (\Gamma \backslash G^{{\rm nc}} )$ admet une d\'ecomposition en somme directe de Hilbert discr\`ete
$$L^2 (\Gamma \backslash G^{{\rm nc}} ) \simeq \sum^{\oplus} \mbox{Hom}_G (\pi , L^2 (\Gamma \backslash G^{{\rm nc}} )) \otimes \pi
= \sum^{\oplus} n_{\Gamma} (\pi ) \pi ,$$
o\`u $\pi$ parcourt cette fois le dual unitaire de $G^{{\rm nc}}$ et la multiplicit\'e
$$n_{\Gamma} (\pi ) := \mbox{dim}_{{\Bbb C}} \mbox{Hom}_G (\pi , L^2 (\Gamma \backslash G^{{\rm nc}} )) < \infty .$$
Alors la formule de Matsushima est r\'esum\'ee dans le lemme suivant.

\begin{lem}[\cite{BorelWallach}, \cite{Matsushima}]
Sous les notations pr\'ec\'edentes. On a
\begin{eqnarray} \label{lem1}
H^* (\pi : \Gamma ) \simeq n_{\Gamma} H^* (\mathfrak{g} , K ; \pi ) ,
\end{eqnarray}
\begin{eqnarray} \label{lem2}
H^* (S(\Gamma )) = \bigoplus_{\pi \in \hat{G}^{{\rm nc}}_0}  H^* (\pi : \Gamma ).
\end{eqnarray}
\end{lem}

En passant \`a la limite (inductive) sur les sous-groupes de congruence $\Gamma \subset G$, nous parlerons
de $\pi$-composante de la cohomologie $H^* (Sh^0 G)$ de la vari\'et\'e de Shimura $Sh^0 G$, ce que nous noterons
$H^* (\pi : Sh^0 G)$. La d\'ecomposition (\ref{lem2}) se traduit alors en
\begin{eqnarray} \label{dec cohom}
H^* (Sh^0 G) = \bigoplus_{\pi \in \hat{G}^{{\rm nc}}_0}  H^* (\pi : Sh^0 G ).
\end{eqnarray}

\medskip

D'apr\`es Parthasarathy \cite{Parthasarathy2}, Kumaresan \cite{Kumaresan} et Vogan-Zuckerman \cite{VoganZuckerman},
les $(\mathfrak{g} , K)$-modules unitarisables ayant des groupes de $(\mathfrak{g}, K)$-cohomologie non triviaux peuvent
\^etre d\'ecrit comme suit. Notons toujours $\mathfrak{t}_0 = $Lie$(T)$ une sous-alg\`ebre de Cartan de $\mathfrak{k}_0$.
On consid\`ere les sous-alg\`ebres paraboliques $\theta$-stable $\mathfrak{q} \subset \mathfrak{g}$~: $\mathfrak{q}  =
\mathfrak{l}  \oplus \mathfrak{u}$ \cite{VoganZuckerman}, o\`u $\mathfrak{l}$ est le centralisateur d'un \'el\'ement
$X\in i \mathfrak{t}_0$ et $u$ est le sous-espace engendr\'e par les racines positives de $X$ dans $\mathfrak{g}$.
Alors $\mathfrak{q} $ est stable sous $\theta$; on en d\'eduit une d\'ecomposition $\mathfrak{u}  = (\mathfrak{u}  \cap
\mathfrak{k} ) \oplus (\mathfrak{u}  \cap \mathfrak{p} )$. Soit $R = \mbox{dim} (\mathfrak{u} \cap \mathfrak{p} )$.

Associ\'e \`a $\mathfrak{q}$, se trouve un $(\mathfrak{g} , K)$-module irr\'eductible bien d\'efini $A_{\mathfrak{q}}$ caract\'eris\'e
par les propri\'et\'es suivantes. Supposons effectu\'e un choix de racines positives pour $(\mathfrak{k} , \mathfrak{t} )$ de fa\c{c}on
compatible avec $\mathfrak{u}$. Soit $e(\mathfrak{q})$ un g\'en\'erateur de la droite $\bigwedge^R (\mathfrak{u} \cap \mathfrak{p} )$.
Alors $e(\mathfrak{q})$ est le vecteur de plus haut poids d'une repr\'esentation irr\'eductible $V(\mathfrak{q})$ de $K$ contenue
dans $\bigwedge^R \mathfrak{p}$; et dont le plus haut poids est donc n\'ecessairement $2\rho (\mathfrak{u} \cap \mathfrak{p} )$.
La classe d'\'equivalence du $(\mathfrak{g} , K)$-module $A_{\mathfrak{q}}$ est alors uniquement caract\'eris\'ee
par les deux propri\'et\'es suivantes.
\begin{eqnarray} \label{Aq1}
\begin{array}{l}
A_{\mathfrak{q}} \mbox{ {\it est unitarisable avec le m\^eme caract\`ere infinit\'esimal que la}} \\
\mbox{{\it  repr\'esentation triviale}}
\end{array}
\end{eqnarray}
\begin{eqnarray} \label{Aq2}
\mbox{Hom}_K (V(\mathfrak{q} ), A_{\mathfrak{q}} ) \neq 0.
\end{eqnarray}
Remarquons que la classe du module $A_{\mathfrak{q}}$ ne d\'epend alors en fait que de l'intersection $\mathfrak{u} \cap \mathfrak{p}$,
autrement dit deux sous-alg\`ebres paraboliques $\mathfrak{q} = \mathfrak{l} \oplus \mathfrak{u}$ et $\mathfrak{q}' =\mathfrak{l}' \oplus
\mathfrak{u}'$ v\'erifiant $\mathfrak{u} \cap \mathfrak{p} = \mathfrak{u}' \cap \mathfrak{p}$ donnent lieu \`a une m\^eme classe de
module cohomologique.

De plus, $V(\mathfrak{q})$ intervient avec multiplici\'e $1$ dans $A_{\mathfrak{q}}$ et $\bigwedge^R (\mathfrak{p} )$, et
\begin{eqnarray} \label{gKcohom}
H^i (\mathfrak{g} , K , A_{\mathfrak{q}} ) \cong \mbox{Hom}_{L\cap K} ( \bigwedge^{i-R} (\mathfrak{l} \cap \mathfrak{p} ), {\Bbb C}).
\end{eqnarray}
Ici $L$ est un sous-groupe de $K$ d'alg\`ebre de Lie $\mathfrak{l}$.

Si $\Gamma$ est un sous-groupe de congruence de $G$, la $A_{\mathfrak{q}}$-composante $H^R (A_{\mathfrak{q}} : \Gamma )$ de $H^R (
S(\Gamma ))$ sera dite {\it fortement primitive}. D'apr\`es ce que nous avons rappel\'e ci-dessus la $A_{\mathfrak{q}}$-composante
fortement primitive est donc la somme sur une base
$\{ \varphi \}$ de Hom$_{\mathfrak{g} , K} (A_{\mathfrak{q}} , C^{\infty} (\Gamma \backslash
G^{nc} ))$ des formes diff\'erentielles $\omega_{\varphi}$ d\'efinies par
$$\omega_{\varphi} (g. \lambda ) = \varphi (\omega (\lambda )) (g) \; \; (\lambda \in \bigwedge^R \mathfrak{p} , \; g\in G^{{\rm nc}} ),$$
o\`u $\omega : \bigwedge^R \mathfrak{p} \rightarrow A_{\mathfrak{q}}$ est une $K$-application non nulle (uniquement d\'efinie \`a un
scalaire pr\`es) qui factorise n\'ecessairement via la composante isotypique $V(\mathfrak{q} )$.
De m\^eme nous parlerons de $H^R (A_{\mathfrak{q}} : Sh^0 G)$ comme de la {\it $A_{\mathfrak{q}}$-composante fortement primitive}.

\subsubsection*{Modules cohomologiques et diagrammes de Young}

Nous avons vu comment associer une sous-alg\`ebre parabolique $\theta$-stable $\mathfrak{q}$ \`a un \'el\'ement $X
= (x_1 , \ldots , x_p ; y_1 , \ldots , y_q ) \in i \mathfrak{t}_0$ (les $x_i$, $y_j$ sont donc tous r\'eels).
Rappelons le choix fix\'e (\ref{racines compactes}) de racines simples compactes positives.
Apr\`es conjugaison par un \'el\'ement de $K$, on peut supposer, et nous le supposerons effectivement par la suite, que $X$ est dominant
par rapport \`a $\Phi (\mathfrak{b}_K , \mathfrak{t} )$, {\it i.e.} que $\alpha (X) \geq 0$ pour tout $\alpha \in \Phi (\mathfrak{b}_K ,
\mathfrak{t} )$; il satisfait alors aux in\'egalit\'es
$$x_1 \geq \ldots \geq x_p \: \mbox{ et } \; y_q \geq \ldots \geq y_1 .$$

Nous allons maintenant associer \`a $X$ un couple de diagrammes de Young (ou, suivant la litt\'erature, diagrammes de Ferrers),
qui coderont compl\`etement le module cohomologique associ\'e.
Rappelons qu'une {\it partition} est une suite d\'ecroissante $\lambda$ d'entiers naturels $\lambda_1 \geq \ldots \geq
\lambda_l \geq 0$. Les entiers $\lambda_1 , \ldots , \lambda_l$ sont des {\it parts}. La {\it longueur} $l(\lambda )$ d\'esigne le
nombre de parts non nulles, et le {\it poids} $|\lambda |$, la somme des parts. On se soucie peu, d'ordinaire, des parts nulles~: on se
permet en particulier, le cas \'ech\'eant, d'en rajouter ou d'en \^oter

Le {\it diagramme de Young} de $\lambda$, que l'on notera \'egalement $\lambda$, s'obtient en superposant, de haut en bas, des lignes
dont l'extr\'emit\'e gauche est sur une m\^eme colonne, et de longueurs donn\'ees par les parts de $\lambda$. Par sym\'etrie
diagonale, on obtient le diagramme de Young de la {\it partition conjugu\'ee}, que l'on notera $\lambda^*$.

Le diagramme de Young de la partition $\lambda = (5,3,3,2)$ et de sa conjugu\'e sont donc~:

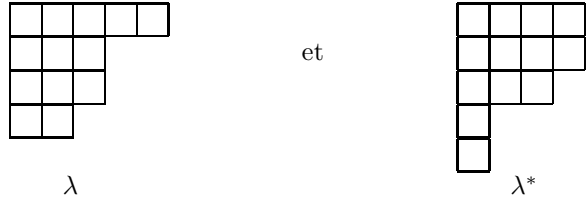
\bigskip

\begin{center}
\begin{tabular}{lcl}

\begin{tabular}{|c|c|c|c|c|} \hline
 & & & &\\ \hline
\end{tabular}
& \hspace{3cm}  &
\begin{tabular}{|c|c|c|c|} \hline
 & & & \\ \hline
\end{tabular}
\\

\begin{tabular}{|c|c|c|} \hline
 & & \\ \hline
\end{tabular}
& et  &
\begin{tabular}{|c|c|c|c|} \hline
 & & & \\ \hline
\end{tabular}
\\

\begin{tabular}{|c|c|c|} \hline
 & & \\ \hline
\end{tabular}
& \hspace{3cm} &
\begin{tabular}{|c|c|c|} \hline
 & & \\ \hline
\end{tabular}
\\

\begin{tabular}{|c|c|} \hline
 & \\ \hline
\end{tabular}
& \hspace{3cm} &
\begin{tabular}{|c|} \hline
 \\ \hline
\end{tabular}
\\

& \hspace{3cm} &
\begin{tabular}{|c|} \hline
 \\ \hline
\end{tabular}
\\

$\quad \quad \lambda$ & \hspace{3cm} & $\quad \quad \lambda^*$

\end{tabular}

\end{center}

\bigskip

Soient $\lambda$ et $\mu$ deux partitions telles que $\mu$ contienne $\lambda$, ce que nous noterons $\lambda \subset \mu$. Notons
$\mu / \lambda$ le compl\'ementaire du diagramme de $\lambda$ dans celui de $\mu$~: c'est une {\it partition gauche} son diagramme
est un {\it diagramme gauche}. Dans la pratique les partitions $\lambda$ que nous rencontrerons seront incluses dans la {\it partition
rectangulaire} $p\times q = (\underbrace{q, \ldots , q}_{p \; {\rm fois}})$, le diagramme gauche $p\times q /\lambda$ est alors
le diagramme de Young d'une partition auquel on a appliqu\'e une rotation d'angle $\pi$; nous noterons $\hat{\lambda}$ cette partition,
la {\it partition compl\'ementaire} de $\lambda$ dans $p \times q$.
Par exemple, la partition $\lambda = (5,3,3,2)$ est incluse dans le rectangle $5\times 5$, et
dans ce rectangle, $\hat{\lambda} = (5 , 3 , 2 , 2)$.

\bigskip

Nous associons maintenant \`a notre \'el\'ement $X \in i \mathfrak{t}_0$ un couple $(\lambda , \mu )$ de partitions comme suit.
\begin{itemize}
\item La partition $\lambda \subset p\times q$ est associ\'ee au sous-diagramme de Young de $p\times q$ constitu\'e des cases de
coordonn\'ees $(i,j)$ telles que $x_i > y_j$.
\item La partition $\mu \subset p\times q$ est associ\'ee au sous-diagramme de Young de $p\times q$ constitu\'e des cases de
coordonn\'ees $(i,j)$ telles que $x_i \geq y_j$.
\end{itemize}

Le lemme suivant est absolument imm\'ediat.

\begin{lem} \label{triv}
Le couple de partitions $(\lambda , \mu)$ associ\'e \`a un \'el\'ement $X \in i \mathfrak{t}_0$ v\'erifie~:
\begin{enumerate}
\item la suite d'inclusion $\lambda \subset \mu \subset p\times q$, et
\item que le diagramme gauche $\mu / \lambda$ est une r\'eunion de diagrammes rectangulaires $p_i \times q_i$, $i=1, \ldots ,m$
ne s'intersectant qu'en des sommets.
\end{enumerate}
R\'eciproquement, \'etant donn\'e un couple de partitions $(\lambda , \mu )$ v\'erifiant 1 et 2, on peut toujours trouver
un \'el\'ement $X \in i \mathfrak{t}_0$ tel que $(\lambda , \mu )$ soit le couple de partitions associ\'e \`a $X$.
\end{lem}

Nous dirons d'un couple de partitions $(\lambda , \mu)$ qu'il est {\it compatible} (ou {\it compatible dans $p\times q$} en cas
d'ambiguit\'e) s'il v\'erifie les points 1 et 2 du lemme \ref{triv}.

Remarquons maintenant que si $X$ et $X'$ sont deux \'el\'ements de $i \mathfrak{t}_0$ de m\^eme couple de partitions associ\'e $(\lambda
,\mu)$ et de sous-alg\`ebres paraboliques associ\'ees respectives $\mathfrak{q}$ et $\mathfrak{q}'$ alors
$\mathfrak{q} \cap \mathfrak{u} = \mathfrak{q}' \cap \mathfrak{u}'$. On d\'eduit donc de la remarque suivant la d\'efinition des modules
$A_{\mathfrak{q}}$ et du lemme \ref{triv} que chaque couple compatible de partitions $(\lambda , \mu )$ d\'efinit sans ambiguit\'e une
classe d'\'equivalence de $(\mathfrak{g} , K)$-modules que nous noterons $A(\lambda , \mu )$. Nous nous autoriserons \`a parler de
``la'' sous-alg\`ebre parabolique $\mathfrak{q}(\lambda , \mu ) = \mathfrak{l}(\lambda , \mu )\oplus \mathfrak{u}(\lambda , \mu )$ de
$(\mathfrak{g}, K)$-module associ\'e $A(\lambda , \mu )$, l'important
pour nous est qu'une telle sous-alg\`ebre existe (d'apr\`es le lemme \ref{triv}). Nous supposerons de plus, ce que l'on peut
toujours faire, que le groupe $L(\lambda , \mu )$ associ\'e \`a la sous-alg\`ebre de Levi $\mathfrak{l} (\lambda , \mu )$
n'a pas de facteurs compacts non ab\'elien. Il est alors facile de voir que
\begin{eqnarray} \label{L}
L(\lambda , \mu )/(L(\lambda ,\mu )\cap K) = \prod_{i=1}^m U(p_i , q_i) /U(p_i ) \times U(q_i ).
\end{eqnarray}

Les r\'esultats de Parthasarathy, Kumaresan et Vogan-Zuckerman mentionn\'es plus haut affirment alors que
$$\hat{G}^{{\rm nc}}_{{\rm VZ}} := \{ A(\lambda , \mu ) \; : \; (\lambda , \mu ) \mbox{ est un couple compatible de partitions} \} (
\subset \hat{G}^{{\rm nc}}_0 \subset \hat{G}^{{\rm nc}})$$
est l'ensemble des $(\mathfrak{g} , K)$-modules ayant des groupes de $(\mathfrak{g} , K)$-cohomologie non nuls.

\bigskip

Comme repr\'esentation de $K_{{\Bbb C}}$, $\mathfrak{p}^+ = E \otimes F^*$ et $\mathfrak{p} = (E \otimes F^*) \oplus (E \otimes F^* )^*$.
Il est bien connu (cf. \cite{Fulton}) qu'\`a chaque partition $\lambda$, il correspond une repr\'esentation irr\'eductible
$E^{\lambda}$ de $GL(E)$.

Consid\'erons la repr\'esentation de $K_{{\Bbb C}}$
\begin{eqnarray} \label{Vlambda}
V(\lambda ) := E^{\lambda} \otimes (F^{\lambda^*})^* .
\end{eqnarray}
C'est une sous-repr\'esentation irr\'eductible de $\bigwedge^{|\lambda |} (E \otimes F^* )$; son vecteur de plus haut poids est
\begin{eqnarray} \label{poids+}
v (\lambda ):= \bigwedge_{i=1}^p \bigwedge_{j=1}^{\lambda_i} e_i \otimes f_j^*
\end{eqnarray}
et son vecteur de plus bas poids est
\begin{eqnarray} \label{poids-}
w (\lambda ):= \bigwedge_{i=1}^p \bigwedge_{j=1}^{\lambda_i} e_{p-i+1} \otimes f_{q-j+1}^* .
\end{eqnarray}
On peut montrer, cf. \cite{Fulton}, que la repr\'esentation
\begin{eqnarray} \label{dec}
\bigwedge \mathfrak{p}^+ = \bigwedge (E \otimes F^* ) = \bigoplus_{\lambda  \subset p\times q} V(\lambda ),
\end{eqnarray}
o\`u chaque sous-espace irr\'eductible $V(\lambda )$ apparait avec multiplicit\'e un.

Soit maintenant $(\lambda , \mu )$ un couple compatible de partitions. Le vecteur
\begin{eqnarray} \label{v(lambda,mu)}
v(\lambda ) \otimes w( \hat{\mu} )^* \in \bigwedge^{|\lambda|} (E\otimes F^* ) \otimes \bigwedge^{|\hat{\mu}|} (E\otimes F^* )^* = \bigwedge^{
|\lambda|,|\hat{\mu}|} \mathfrak{p} \subset \bigwedge^{|\lambda |+|\hat{\mu} |} \mathfrak{p}
\end{eqnarray}
est un vecteur de plus haut poids $2\rho (\mathfrak{u}(\lambda , \mu ) \cap \mathfrak{p})$ et engendre donc sous l'action de
$K_{{\Bbb C}}$ un sous-module irr\'eductible que l'on note $V(\lambda , \mu )$. Ce module est isomorphe \`a $V(\mathfrak{q}(\lambda , \mu
))$.

\subsubsection*{Classes de Chern et diagrammes de Young}

Soit ${\Bbb G}_{p,q}$ la grassmannienne des sous-espaces complexes de dimension $p$ dans ${\Bbb C}^{p+q}$.
Soit $x_0 \in {\Bbb G}_{p,q}$ le point
base correspondant au sous-espace complexe de dimension $p$ constitu\'e des vecteurs de ${\Bbb C}^{p+q}$ dont les $q$ derni\`eres
coordonn\'ees sont toutes nulles. Le groupe $G^{{\rm nc}} = U(p,q)$ agit naturellement sur ${\Bbb G}_{p,q}$ et l'orbite $U(p,q) .x_0$
s'identifie \`a l'espace sym\'etrique hermitien $X_G$.

Cette construction se comprend plus g\'en\'eralement de la fa\c{c}on suivante. Soit $G^{{\rm nc}}_{{\Bbb C}}$ le complexifi\'e du
groupe $G^{{\rm nc}}$. Soit $P_{{\Bbb C}}^-$ le sous-groupe parabolique de $G^{{\rm nc}}_{{\Bbb C}}$ d'alg\`ebre de Lie
$\mathfrak{p}^- \oplus \mathfrak{k}$. Alors, le {\it dual compact} $\hat{X}_G = G^{{\rm nc}}_{{\Bbb C}} / P_{{\Bbb C}}^-$ de
$X_G$ est \'egalement un espace sym\'etrique hermitien et l'inclusion $G^{{\rm nc}} /K =
G^{{\rm nc}} / G^{{\rm nc}} \cap P_{{\Bbb C}}^- \subset G^{{\rm nc}}_{{\Bbb C}} / P^-_{{\Bbb C}}$ r\'ealise $X_G$ comme un
domaine born\'e dans $\hat{X}_G$. Dans notre cas $\hat{X}_G$ s'identifie \`a la grassmannienne ${\Bbb G}_{p,q}$ et on retrouve
le plongement d\'ecrit au pr\'ec\'edent paragraphe. Remarquons que si l'on introduit la sous-alg\`ebre de Lie r\'eelle
$\mathfrak{g}_u = \mathfrak{k}_0 \oplus i \mathfrak{p}_0$ de $\mathfrak{g}$ et $G^{{\rm nc}}_u$ le sous-groupe (compact malgr\'e
le nc) de $G^{{\rm nc}}_{{\Bbb C}}$ d'alg\`ebre de Lie $\mathfrak{g}_u$, alors $\hat{X}_G = G^{{\rm nc}}_u / K$.

On dispose sur la grassmannienne ${\Bbb G}_{p,q}$ d'un {\it fibr\'e tautologique} $\hat{T}$ de rang $p$, dont la fibre au-dessus d'un
sous-espace $W$ de ${\Bbb C}^{p+q}$ est $W$ lui m\^eme. De fa\c{c}on analogue, le {\it fibr\'e quotient} $\hat{Q}$, de rang $q$, a
pour fibre au-dessus de $W$ le quotient ${\Bbb C}^{p+q} /W$.

Soit $\Gamma$ un sous-groupe de congruence de $G$. Le groupe $GL(p+q, {\Bbb C})$ agit sur $\hat{T}$ et $\hat{Q}$. Par restriction,
le groupe $U(p,q)$ agit sur $\hat{T}_{|X_G}$ et sur $\hat{Q}_{|X_G}$. En quotientant par l'action de $\Gamma$ sur $\hat{T}_{|X_G}$
et sur $\hat{Q}_{|X_G}$, on obtient deux fibr\'es sur $S(\Gamma )$ que nous noterons respectivement $T$ et $Q$.

D'apr\`es un th\'eor\`eme classique de Cartan, l'espace $H^* ({\Bbb G}_{p,q})$ peut \^etre identifi\'e avec l'espace des formes
diff\'erentielles $U(p+q)$-invariantes. Soit $\omega$ une forme diff\'erentielle $U(p+q)$-invariante sur ${\Bbb G}_{p,q}$ et soit
$\bar{\omega}$ une forme diff\'erentielle $U(p,q)$-invariante sur $X_G$ \'egale \`a $\omega$ au point $x_0$. Puisque
$\omega$ est, en particulier, $\Gamma$-invariante elle induit une forme (n\'ecessairement ferm\'ee) sur $S(\Gamma )$ qui d\'efinit
donc une classe de cohomologie dans $H^* (S(\Gamma ))$. On a ainsi construit une application
\begin{eqnarray} \label{eta}
\eta : H^* ({\Bbb G}_{p,q} ) \rightarrow H^* (S(\Gamma )).
\end{eqnarray}
Il est bien connu que $\eta$ est injective. Le lemme suivant est lui aussi classique, on peut en trouver une d\'emonstration dans
\cite{Parthasarathy}.

\begin{lem} \label{chern}
Si $\hat{C}_1 , \ldots , \hat{C}_p$ (resp. $\hat{C}_1 ', \ldots , \hat{C}_q '$) sont les classes de Chern du fibr\'e $\hat{T}$ (resp.
$\hat{Q}$), alors $C_i := (-1)^i \eta (\hat{C}_i )$ (resp. $C_i ':= (-1)^i \eta (\hat{C}_i ')$) est la $i$-\`eme classe de Chern
du fibr\'e $T$ (resp. $Q$).
\end{lem}

Nous allons relier les classes de Chern $\hat{C}_i$ et $\hat{C}_i '$ \`a des sous-espaces $K$-invariants de $\bigwedge \mathfrak{p}$.

Remarquons d'abord qu'en utilisant (\ref{dec}) et son dualis\'e~:
$$\bigwedge \mathfrak{p}^- = \bigoplus_{\lambda \subset p\times q} V(\lambda )^* ,$$
on peut d\'ecrire une base $\{ C_{\nu } \; : \; \nu \subset p\times q \}$ de l'espace $\left( \bigwedge \mathfrak{p} \right)^K$
des vecteurs $K$-invariants de $\bigwedge \mathfrak{p}$ param\'etr\'ee par l'ensemble des partitions $\nu \subset p\times q$.
On prend $C_{\nu }:= \sum_l z_l \otimes z_l^* $ o\`u $\{ z_l \}$ est une base de $V(\nu ) \subset \bigwedge \mathfrak{p}^+$ et
$\{ z_l^* \}$ la base duale de $V(\nu )^* \subset \bigwedge \mathfrak{p}^-$.

Soit $\nu \subset p\times q$ une partition. L'\'el\'ement $C_{\nu} \in \bigwedge \mathfrak{p}$ est invariant sous l'action de $K$.
Or le th\'eor\`eme de Cartan mentionn\'e plus haut identifie $\left(\bigwedge \mathfrak{p} \right)^K$ et $H^* ({\Bbb G}_{p,q} )$.
On peut donc voir $C_{\nu}$ comme une classe de cohomologie dans $H^* ({\Bbb G}_{p,q} )$.
D'apr\`es un th\'eor\`eme de Kostant \cite[Theorem 6.15]{Kostant}, $C_{\nu}$ est un multiple non nul de la classe de cohomologie
associ\'ee \`a la sous-vari\'et\'e de Schubert $X_{\nu}$ associ\'ee \`a la partition $\nu \subset p\times q$.

Rappelons qu'une fois un drapeau complet fix\'e
$$0=V_0 \subsetneq \ldots \subsetneq V_{p+q} = {\Bbb C}^{p+q},$$
on associe \`a toute partition $\nu \subset p\times q$, la {\it sous-vari\'et\'e de Schubert}
$$X_{\nu} = \{ W \in {\Bbb G}_{p,q} \; : \; \mbox{dim} (W \cap V_{n+i-\nu_i} ) \geq i , \; 1\leq i \leq m \}.$$
La classe de Schubert associ\'ee $[X_{\nu }] \in H^{2|\nu|} ({\Bbb G}_{p,q} )$ ne d\'epend pas du choix du drapeau. D'apr\`es le
th\'eor\`eme de Kostant cit\'e au pr\'ec\'edent paragraphe $C_{\nu}$ est un multiple non nul de $[X_{\nu} ]$.

D'un autre c\^ot\'e, il est bien connu (cf. \cite{Fulton} ou \cite{Manivel}) que la $k$-i\`eme classe de Chern $\hat{C}_k$
(resp. $\hat{C}_k '$) du fibr\'e $\hat{T}$ (resp. $\hat{Q}$) sur la grassmannienne, est \'egale \`a la classe $[X_{(1^k)}]$
(resp. $[X_{(k)}]$) d'une sous-vari\'et\'e de Schubert associ\'ee \`a la partition $(1^k )=(\underbrace{1 , \ldots , 1}_{k \; {\rm fois}})$
(resp. la partition dont une seule part est non nulle \'egale \`a $k$).

\subsubsection*{Action des classes de Chern sur la cohomologie}

Soit $(\lambda ,\mu )$ un couple compatible de partitions.
Rappelons que le diagramme gauche $\mu / \lambda$ est alors constitu\'e de diagrammes rectangulaires $p_i \times q_i$, $i= 1 , \ldots ,
m$ ne s'intersectant qu'en des sommets. Par exemple si $\mu = (8,8,8,4,4,2)$ et $\lambda = (4,4,4,2,2)$,

\bigskip

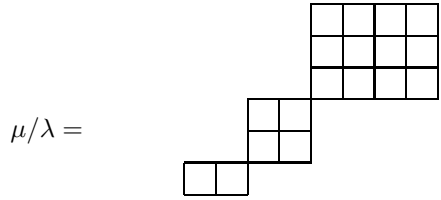
\begin{center}
\begin{tabular}{ccl}
  & \hspace{0.5cm} &
\begin{tabular}{cccc|c|c|c|c|} \cline{5-8}
& & & & & & &  \\ \cline{5-8}
& & & & & & &  \\ \cline{5-8}
& & & & & & &  \\ \cline{5-8}
\end{tabular} \\
$\mu / \lambda = $  & \hspace{0.5cm} &
\begin{tabular}{cc|c|c|} \cline{3-4}
& & & \\ \cline{3-4}
& & & \\ \cline{3-4}
\end{tabular} \\
 & \hspace{0.5cm} &
\begin{tabular}{|c|c|} \hline
& \\ \hline
\end{tabular}
\end{tabular}
\end{center}

\bigskip

\'Etant donn\'e un diagramme gauche $\mu /\lambda$, num\'erotons les cases de droite \`a gauche et de haut en bas; nous appellerons
cet \'etiquetage le {\it num\'erotage inverse} du diagramme gauche. Par exemple, le numerotage inverse de
$(5,4,3,2)/(3,3,1)$ est

\bigskip

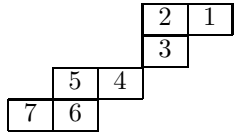
\begin{center}
\begin{tabular}{l}
\hspace{0.405cm} \begin{tabular}{ccc|c|c|} \cline{4-5}
$ $ & $ $  & $ $  & $2$ & $1$ \\ \cline{4-5}
\end{tabular} \\
\hspace{0.405cm} \begin{tabular}{ccc|c|} \cline{4-4}
$ $ & $ $ & $ $ & $3$ \\ \cline{4-4}
\end{tabular} \\
\hspace{0.055cm} \begin{tabular}{c|c|c|} \cline{2-3}
$ $ & $5$ & $4$ \\ \cline{2-3}
\end{tabular} \\
\begin{tabular}{|c|c|} \cline{1-2}
$7$ & $6$ \\ \cline{1-2}
\end{tabular}
\end{tabular}
\end{center}

\bigskip

Nous appellerons {\it sous-diagramme gauche} d'un diagramme gauche $\mu / \lambda$, tout diagramme $\mu ' / \lambda$ o\`u
$\mu '$ est le diagramme d'une partition v\'erifiant $\lambda \subset \mu ' \subset \mu$.
On peut alors donner les d\'efinitions suivantes qui seront fondamentales dans la suite du texte.

\bigskip

\noindent
{\bf D\'efinitions} \begin{itemize}
\item Nous dirons d'une partition $\nu \subset p\times q$ qu'elle est une {\it image} d'un diagramme gauche $\mu / \lambda$ s'il existe
une bijection entre les cases des diagrammes $\nu$ et $\mu / \lambda$
telle que si une case $A$ est au-dessus (au sens large) et \`a gauche (au sens large)
d'une case $B$  dans l'un des diagrammes, les cases correspondantes $A'$ et $B'$ de l'autre diagramme sont dans l'ordre du num\'erotage
inverse.
\item Nous dirons ensuite d'une partition $\nu \subset p\times q$ qu'elle peut {\it s'inscrire} dans un diagramme gauche $\mu / \lambda$
si elle est une image d'un sous-diagramme gauche de $\mu / \lambda$.
\end{itemize}

\bigskip

Nous allons maintenant tenter de donner plus de sens \`a ces d\'efinitions. Commen\c{c}ons par dire qu'elles sont reli\'ees aux
nombres de Littlewood-Richardson $c_{\lambda \nu}^{\mu}$. Nous renvoyons au livre de Fulton \cite{Fulton} pour un grand
nombres de d\'efinitions de ces nombres. En ce qui nous concerne, remarquons imm\'ediatement que la notion d'image
d\'efinie ci-dessus est due \`a Zelevinsky \cite{Zelevinsky} qui montre notamment que $c_{\lambda \nu}^{\mu}$ est \'egal au
nombre d'images entre $\nu$ et $\mu / \lambda$ \footnote{Ici le couple de partitions $(\lambda , \mu )$ n'a pas besoin d'\^etre
n\'ecessairement compatible.}.

Nous aurons besoin d'utiliser plusieurs caract\'erisations diff\'erentes des nombres de Littlewood-Richardson.
Rappelons qu'un {\it tableau de Young} $T$ est le remplissage des cases d'un diagramme de Young  $\mu$ par des entiers $\geq 1$ de
mani\`ere
\begin{enumerate}
\item croissante (au sens large) le long des lignes (de gauche \`a droite), et
\item strictement croissante de haut en bas suivant chaque colonne.
\end{enumerate}
\'Etant donn\'e une partition $\lambda$, nous noterons $U(\lambda )$ le tableau de Young obtenu \`a partir du diagramme $\lambda$
en remplissant toutes les cases de la $i$-\`eme ligne par des $i$.
Nous parlerons \'egalement de {\it tableaux gauches}, {\it i.e.} de tableaux sur des diagrammes gauches.
Nous admettrons dans cette sous-section que le lecteur est familier avec le produit $\cdot$ de deux tableaux de Young ou la rectification
Rect d'un tableau gauche. Rappelons alors que le nombre de Littlewood-Richardson $c_{\lambda \nu}^{\mu}$ est \'egal
\begin{enumerate}
\item au nombre de tableaux $T$ sur $\lambda$ tels que $T \cdot U(\nu ) = U(\mu )$;
\item au nombre de tableaux gauches $S$ sur $\mu / \lambda$ tels que Rect$(S)=U(\nu )$.
\end{enumerate}

\'Etant donn\'e un couple compatible $(\lambda , \mu )$ de partitions, le diagramme gauche $\mu / \lambda$ est r\'eunion de
diagrammes rectangulaires $p_i \times q_i$, $i=1, \ldots , m$. Si pour chaque entier $i$ entre $1$ et $m$ on choisit un
diagramme $\alpha_i \subset p_i \times q_i$, on peut former un sous-diagramme gauche
$\alpha_1 * \ldots * \alpha_m$ de $\mu/\lambda$ en pla\c{c}ant les
diagrammes $\alpha_i$ dans le coin sup\'erieur gauche de chacun des diagrammes $p_i \times q_i$. On notera (de mani\`ere coh\'erente
avec la notation des nombre de Littlewood-Richardson) $c_{\alpha_1 \ldots \alpha_m}^{\nu}$ le nombres d'images entre $\nu$ et
le diagramme gauche $\alpha_1 * \ldots * \alpha_m$.

Nous aurons besoin du lemme suivant.

\begin{lem} \label{LR}
Soit $(\lambda , \mu )$ un couple compatible de partitions tel que $\mu / \lambda = \bigcup_{i=1}^m p_i \times q_i$.
Notons $P_i = p_{i+1} + \ldots + p_m$ et $Q_i = q_{i+1} + \ldots + q_m$. Soit $\nu$ une partition.
\begin{enumerate}
\item On a
$$c_{\lambda  \nu}^{\mu} = \sum_{
\begin{footnotesize}
\begin{array}{l}
\beta_i \subset P_i \times Q_i \\
i=1 ,\ldots , m-2
\end{array}
\end{footnotesize}} \prod_{i=1}^{m-1} c_{p_i \times q_i \beta_i}^{\beta_{i-1}},$$
o\`u on a not\'e $\beta_0 = \nu$ et $\beta_{m-1} = p_m \times q_m$.
\item Le diagramme $\nu$ s'inscrit dans $\mu / \lambda$ si et seulement s'il existe des diagrammes
$\alpha_i \subset p_i \times q_i$, pour $i=1, \ldots , m$ tel que $c_{\alpha_1 \ldots \alpha_m}^{\nu}\neq 0$. Et
\begin{eqnarray} \label{LRmult}
c_{\alpha_1 \ldots \alpha_m}^{\nu} = \sum_{
\begin{footnotesize}
\begin{array}{l}
\beta_i \subset P_i \times Q_i \\
i=1, \ldots , m-2
\end{array}
\end{footnotesize}} \prod_{i=1}^{m-1} c_{\alpha_i \beta_i}^{\beta_{i-1}} ,
\end{eqnarray}
o\`u on a not\'e $\beta_0 = \nu$ et $\beta_{m-1} = \alpha_m$.
\end{enumerate}
\end{lem}
{\it D\'emonstration.} Fixons des diagrammes $\alpha_i \subset p_i \times q_i$, pour $i=1, \ldots ,m$ et
$\beta_i \subset P_i \times Q_i$, pour $i=1, \ldots , m-2$. Nous noterons $\beta_0 = \nu$ et $\beta_{m-1} = \alpha_m$.
Pour tout entier $i=1, \ldots , m-1$, notons alors ${\cal E}(\beta_i )$ l'ensemble
$$\left\{ T_i \mbox{ tableau de Young sur } \alpha_i \; : \; T_i \cdot U(\beta_i ) = U(\beta_{i-1} ) \right\}.$$
Le cardinal de ${\cal E}_i$ est donc \'egal \`a $c_{\alpha_i \beta_i}^{\beta_{i-1}}$.

Soit ${\cal T}$ l'ensemble
$$\left\{ U \mbox{ tableau sur le sous-diagramme gauche }
\alpha_1 * \ldots * \alpha_m \subset \mu / \lambda \; : \; \mbox{Rect}(U) = U(\nu ) \right\}.$$
Alors d'apr\`es Zelevinsky, le cardinal de ${\cal T}$ est \'egal \`a $c_{\alpha_1 \ldots \alpha_m}^{\nu}$.

Consid\'erons maintenant l'application
$$\Phi : \left\{
\begin{array}{lcl}
\bigcup_{
\begin{footnotesize}
\begin{array}{l}
\beta_i \subset P_i \times Q_i \\
i=1, \ldots , m-2
\end{array}
\end{footnotesize}} \prod_{i=1}^{m-1} {\cal E}(\beta_i) & \longrightarrow & {\cal T} \\
(T_1 , \ldots , T_{m-1} ) & \longmapsto & T_1 * \ldots *T_{m-1} *U(\alpha_m ) .
\end{array}
\right.$$
L'application $\Phi$ est clairement injective. Montrons qu'elle est surjective. Il est clair qu'un tableau $U$ dans ${\cal T}$ s'\'ecrit $U= T_1 * \ldots * T_{m-1} *U(\alpha _m )$
o\`u chaque $T_i$, pour $i=1, \ldots , m-1$, est un tableau de Young sur le diagramme $\alpha_i$. Nous allons montrer par
r\'ecurrence que si Rect$(U)=U(\nu )$, il existe une suite de diagrammes $\beta_1 , \ldots , \beta_{m-2}$ telle que  chaque $T_i$
appartienne \`a ${\cal E}(\beta_i )$.

Par d\'efinition du produit $.$ sur les tableaux de Young, il d\'ecoule de l'identit\'e
$$\mbox{Rect}(T_1 * \ldots * T_{m-1} * U(\alpha_m ))=U(\nu )$$
que
\begin{eqnarray} \label{truc}
T_1 \cdot \mbox{Rect}(T_2 * \ldots * T_{m-1} * U(\alpha_m )) =U (\nu ).
\end{eqnarray}
Rect$(T_2 * \ldots * T_{m-1} * U(\alpha_m ))$ est un tableau de Young sur un certain diagramme $\beta_1 \subset P_1 \times Q_1$.
De plus, il d\'ecoule de (\ref{truc}) que Rect$(T_2 * \ldots * T_{m-1} * U(\alpha_m ))$ est n\'ecessairement \'egal \`a $U(\beta_1 )$.
On conclut alors la construction des $\beta_i$ par r\'ecurrence.

La formule (\ref{LRmult}) d\'ecoule imm\'ediatement de la bijectivit\'e de $\Phi$, le point 2. du lemme est donc d\'emontr\'e.
Le cas 1. du lemme s'en suit en posant $\alpha_i = p_i \times q_i$. Ce qui conclut la d\'emonstration du lemme \ref{LR}.

\bigskip

Dans la suite, notons $\mathfrak{p}_L^+$ l'intersection $\mathfrak{p}^+ \cap \mathfrak{l}(\lambda , \mu )$. On introduit alors
$E(G^{{\rm nc}} ,  L(\lambda , \mu ))$ (ou juste $E(G, L)$ lorsqu'il n'y aura pas d'ambiguit\'e)
le sous-espace de $\bigwedge \mathfrak{p}^+$ engendr\'e par les tranlat\'es par $K$ du
sous-espace $\bigwedge \mathfrak{p}_L^+$.

Nous pouvons maintenant \'enoncer et d\'emontrer la proposition clef de cette section.

\begin{prop} \label{clef}
Soient $\lambda$, $\mu$ et $\nu$ trois partitions incluses dans $p\times q$ telles que $(\lambda , \mu )$ forme un couple
compatible. Notons $L = L(\lambda , \mu )$. Alors, les \'enonc\'es suivants sont \'equivalents~:
\begin{enumerate}
\item $C_{\nu } . V(\lambda , \mu ) \neq 0$ dans $\bigwedge \mathfrak{p}$;
\item la partition $\nu$ s'inscrit dans $\mu / \lambda$;
\item $V( \nu ) \subset E(G,L)$.
\end{enumerate}
De plus, les \'el\'ements $\{ C_{\nu} . v(\lambda ) \otimes w(\mu )^* \}$, o\`u $\nu$ d\'ecrit l'ensemble des partitions $\subset p\times q$ qui
s'inscrivent dans $\mu / \lambda$, sont lin\'eairement ind\'ependants.
\end{prop}
{\it D\'emonstration.} Pour simplier les notations, posons $\mathfrak{q} = \mathfrak{q} (\lambda , \mu)$, $\mathfrak{u} =
\mathfrak{u} (\lambda , \mu )$ et $\mathfrak{l} = \mathfrak{l} (\lambda , \mu )$. Rappelons que $v(\lambda ) \otimes w (\hat{\mu })^*$
est un g\'en\'erateur de la droite $\bigwedge^{|\lambda |+|\mu |} (\mathfrak{u} \cap \mathfrak{p} )$. Le lemme suivant, simple
exercice d'alg\`ebre lin\'eaire, est d\^u \`a Venkataramana \cite[Lemma 1.3]{Venky}.

L'inclusion $\mathfrak{l} \cap \mathfrak{p} \rightarrow \mathfrak{p}$ induit par dualit\'e (pour la forme de Killing)
une application $\mathfrak{p} \rightarrow \mathfrak{l} \cap \mathfrak{p}$ que nous dirons ``de restriction''.

\begin{lem} \label{lemvenky1}
Consid\'erons l'application de restriction $B : \left( \bigwedge \mathfrak{p} \right)^T \rightarrow \left( \bigwedge (\mathfrak{l}
\cap \mathfrak{p} )\right)^T$ et le cup-produit $A : \left( \bigwedge \mathfrak{p} \right)^T \rightarrow \bigwedge \mathfrak{p}$
donn\'e par $y \mapsto y \wedge v(\lambda ) \otimes w (\hat{\mu })^*$. Alors les noyaux de $A$ et $B$ sont les
m\^emes.
\end{lem}

Soit $\hat{X}_L = L_u / (L_u \cap K)$ le dual compact de $X_L = L/ L \cap K$.
L'espace sym\'etrique hermitien compact $\hat{X}_L$ se plonge naturellement dans la grassmannienne $\hat{X}_G = {\Bbb G}_{p,q}$.
En cohomologie on peut donc parler de l'application de restriction~:
$$\mbox{res} : H^* (\hat{X}_G ) \rightarrow H^* (\hat{X}_L ) .$$
Rappelons que $C_{\nu}$ s'identifie \`a une classe de cohomologie dans $H^* ({\Bbb G}_{p,q} ) \simeq H^* (\hat{X}_G )$.
Le lemme suivant est \'egalement d\^u \`a Venkataramana \cite[Lemma 1.4]{Venky}, on en esquisse la d\'emonstration pour simplifier
la lecture du texte.

\begin{lem} \label{lemvenky2}
$$C_{\nu} . V(\lambda , \mu ) =0 \Leftrightarrow C_{\nu} \in \mbox{Ker} \left( H^* (\hat{X}_G )
\stackrel{\mbox{res}}{\longrightarrow} H^* (\hat{X}_L ) \right) .$$
\end{lem}
{\it D\'emonstration du lemme \ref{lemvenky2}.} Le $K$-module $V(\lambda , \mu )$ est engendr\'e par $v(\lambda ) \otimes w (\hat{\mu} )^*$,
et $C_{\nu}$ est $K$-invariant, on a donc~:
$$C_{\nu} . V(\lambda , \mu ) =0 \Leftrightarrow C_{\nu }. v(\lambda ) \otimes w (\hat{\mu} ) =0.$$
C'est \'equivalent au fait que $C_{\nu}$ appartient au noyau de l'application $A$ du lemme \ref{lemvenky1}.
D'apr\`es ce dernier c'est donc \'equivalent au fait que $C_{\nu}$ appartient au noyau de l'application $B$.
Mais $B(C_{\nu} )=0$ si et seulement si $C_{\nu}$ est dans le noyau de l'application de restriction
$$\left( \bigwedge \mathfrak{p} \right)^K \rightarrow \left( \bigwedge \mathfrak{p}_L \right)^{K\cap L} .$$
Ce qui conclut la d\'emonstration du lemme \ref{lemvenky2}.

\medskip

Remarquons imm\'ediatement que l'ensemble $\{ C_{\nu} . v(\lambda ) \otimes w(\hat{\mu})^* \; : \; \nu \subset p\times q , \; C_{\nu} .V(\lambda , \mu ) \neq 0 \}$
est identique \`a l'ensemble $\{ B(C_{\nu }) . v(\lambda ) \otimes w(\hat{\mu})^*   \; : \; \nu \subset p\times q , \; C_{\nu} .V(\lambda , \mu ) \neq 0 \}$. Mais, via le
plongement $\bigwedge \mathfrak{p}_L \otimes \bigwedge (\mathfrak{p} \cap \mathfrak{u} ) \rightarrow \bigwedge \mathfrak{p}$, les \'el\'ements
$B( C_{\nu }) . v(\lambda ) \otimes w(\hat{\mu} )^*$ sont simplement les tenseurs d\'ecomposables $B(C_{\nu} ) \otimes (v(\lambda ) \otimes w(\hat{\mu})^* )$ et
sont donc lin\'eairement ind\'ependants.

\medskip

Continuons la d\'emonstration de la proposition \ref{clef}. D'apr\`es le lemme \ref{lemvenky1} et (\ref{L}), et en conservant les
m\^emes notations, on doit comprendre l'application naturelle de restriction~:
\begin{eqnarray} \label{res g}
H^* ({\Bbb G}_{p,q} ) \stackrel{\mbox{res}}{\longrightarrow} \bigotimes_{i=1}^m H^* ({\Bbb G}_{p_i , q_i} ),
\end{eqnarray}
o\`u $p_1 + \ldots + p_m \leq p$, $q_1 + \ldots + q_m \leq q$, $\prod_{i=1}^m {\Bbb G}_{p_i , q_i}$ est naturellement plong\'e
dans ${\Bbb G}_{p,q}$, {\it i.e.} via l'inclusion de $L$ dans $G^{{\rm nc}}$ et o\`u l'on a compos\'e l'application induite
en cohomologie par la restriction et l'isomorphisme de Kunneth.

\begin{lem} \label{res de C}
Soit $\nu \subset p\times q$ une partition. Alors,  l'image de $C_{\nu}$ par l'application de restriction
(\ref{res g}) est~:
$$\mbox{res} (C_{\nu} ) = \sum_{
\begin{footnotesize}
\begin{array}{l}
\alpha_i \subset p_i \times q_i \\
i = 1, \ldots ,m
\end{array}
\end{footnotesize}} c_{\alpha_1 \ldots \alpha_m}^{\nu}  C_{\alpha_1 } \otimes \ldots \otimes C_{\alpha_m }.$$
\end{lem}
{\it D\'emonstration du lemme \ref{res de C}.} Supposons tout d'abord que $m=2$, $p_1 + p_2 =p$ et que $q_1 + q_2 =q$.
Le fibr\'e tautologique $\hat{T}$ au-dessus
de ${\Bbb G}_{p,q}$ induit via le plongement ${\Bbb G}_{p_1,q_1} \times {\Bbb G}_{p_2 , q_2 } \subset {\Bbb G}_{p,q},$
un fibr\'e sur le produit ${\Bbb G}_{p_1,q_1} \times {\Bbb G}_{p_2 , q_2 }$ qui est isomorphe au produit du fibr\'e tautologique
$\hat{T}_1$ sur ${\Bbb G}_{p_1 , q_1}$ par le fibr\'e tautologique $\hat{T}_2$
sur ${\Bbb G}_{p_2 , q_2}$. La formule de Whitney sur les classes de Chern implique alors que, via l'isomorphisme de Kunneth,
la classe de Chern totale de ce fibr\'e est \'egale au produit tensoriel des classes de Chern $c(\hat{T}_1 )$ et $c(\hat{T}_2 )$.
Par fonctorialit\'e des classes de Chern, on en d\'eduit que la restriction de $c(\hat{T} )$ au produit
${\Bbb G}_{p_1,q_1} \times {\Bbb G}_{p_2 , q_2 }$ est donn\'ee par~:
\begin{eqnarray} \label{chern}
c(\hat{T} )_{|({\Bbb G}_{p_1,q_1} \times {\Bbb G}_{p_2 , q_2 })} = c(\hat{T}_1 ) \otimes c(\hat{T}_2 ).
\end{eqnarray}

Notons dor\'enavant, $C_{\alpha}^1$ (resp. $C_{\beta}^2$) la classe de Schubert associ\'ee \`a une partition $\alpha \subset p_1 \times
q_1$ (resp. $\beta \subset p_2 \times q_2$) dans la cohomologie $H^* ({\Bbb G}_{p_1 , q_1} )$ (resp. $H^* ({\Bbb G}_{p_2 , q_2})$).
L'\'equation (\ref{chern}) implique alors que pour tout entier $k \leq p$,
\begin{eqnarray} \label{relations}
(C_{(1^k)})_{|({\Bbb G}_{p_1,q_1} \times {\Bbb G}_{p_2 , q_2 })} = \sum_{a+b=k} C_{(1^a )}^1 \otimes C_{(1^b )}^2 .
\end{eqnarray}

Notons $\Lambda_m$ l'anneau des polyn\^omes sym\'etriques \`a coefficients entiers de $m$ variables. Rappelons qu'une base de
$\Lambda_m$ est fournit par les {\it fonctions de Schur} dont on renvoie \`a \cite{Fulton} ou \cite{Manivel} pour une
d\'efinition. Plus pr\'ecisemment, lorsque $\lambda$ d\'ecrit l'ensemble des partitions de longueur $m$ au plus, les
fonctions de Schur $s_{\lambda}$ forment une base de $\Lambda_m$.

Il est alors classique (cf. \cite{Fulton}, \cite{Manivel}) que l'application
$$\varphi_{p,q} : \Lambda_p \rightarrow H^* ({\Bbb G}_{p,q} ),$$
qui \`a la fonction de Schur $s_{\lambda}$ associe la classe de Schubert $C_{\lambda}$ si $\lambda \subset p\times q$, et
z\'ero sinon, est un morphisme d'anneaux surjectif.
Nous consid\`ererons de m\^eme les morphismes d'anneaux surjectifs~:
$$\varphi_{p_1,q_1} : \Lambda_{p_1} \rightarrow H^* ({\Bbb G}_{p_1 ,q_1 } ),$$
et
$$\varphi_{p_2 ,q_2 } : \Lambda_{p_2} \rightarrow H^* ({\Bbb G}_{p_2 ,q_2 } ).$$

Num\'erotons $(x_1 , \ldots , x_{p_1} , y_1 , \ldots , y_{p_2} )$ les $p$ variables d'une fonction de $\Lambda_p$, alors
$\Lambda_p$ s'identifie au produit tensoriel $\Lambda_{p_1} \otimes \Lambda_{p_2}$ et d'apr\`es \cite[\S 5.2, Exercice 4]{Fulton},
\begin{eqnarray} \label{schur}
s_{\nu } (x_1 , \ldots , x_{p_1} , y_1 , \ldots , y_{p_2} ) = \sum_{\alpha , \beta} c_{\alpha \beta}^{\nu}
s_{\alpha} (x_1 , \ldots , x_{p_1}) s_{\beta} (y_1 , \ldots , y_{p_2} ) .
\end{eqnarray}
Puisque la longueur de la partition sous-jacente \`a un produit de tableaux de Young est toujours inf\'erieure \`a la somme
des longueurs des partitions sous-jacentes \`a ces deux tableaux,
il d\'ecoule de l'identit\'e (\ref{schur}) et des morphismes $\varphi_{p,q}$, $\varphi_{p_1 , q_1 }$ et $\varphi_{p_2 , q_2 }$ que
l'application
$$R : \left\{
\begin{array}{lcl}
H^* ({\Bbb G}_{p,q} ) & \longrightarrow & H^* ({\Bbb G}_{p_1 , q_1 }) \otimes H^* ({\Bbb G}_{p_2 , q_2 }) \\
C_{\nu } & \longmapsto & \sum_{\alpha , \beta } c_{\alpha \beta }^{\nu} C_{\alpha}^1 \otimes C_{\beta}^2
\end{array} \right. $$
est un morphisme d'anneaux. Remarquons que
$$R( C_{(1^k )}) = \sum_{a+b=k} C_{(1^a )}^1 \otimes C_{(1^b )}^2 .$$
Puisque les classes $C_{(1^k )}$ pour $k=0, \ldots , p$ engendrent l'anneau $H^* ({\Bbb G}_{p,q} )$ et compte tenu des
relations (\ref{relations}), le morphisme $R$ est n\'ecessairement \'egal au morphisme de restriction
$$\mbox{res} : H^* ({\Bbb G}_{p,q} )  \rightarrow  H^* ({\Bbb G}_{p_1 , q_1 }) \otimes H^* ({\Bbb G}_{p_2 , q_2 }).$$

Finalement, lorsque $m=2$, $p_1 + p_2 =p$ et $q_1 + q_2 =q$, nous avons d\'emontr\'e que pour toute
partition $\nu \subset p \times q$,
\begin{eqnarray} \label{m=2}
\mbox{res} (C_{\nu }) = \sum_{
\begin{footnotesize}
\begin{array}{l}
\alpha \subset p_1 \times q_1 \\
\beta \subset p_2 \times q_2
\end{array}
\end{footnotesize}} c_{\alpha \beta }^{\nu} C_{\alpha}^1 \otimes C_{\beta}^2 .
\end{eqnarray}

Par r\'ecurrence sur $m$ \footnote{Quitte \`a rajouter des facteurs ${\Bbb G}_{p_i , q_i }$ que l'on ``oubliera'' de fa\c{c}on \`a
ce que $p=p_1 + \ldots + p_m$ et $q=q_1 + \ldots + q_m$.}, on en d\'eduit que si $\nu \subset p\times q$ est une partition, alors
l'image de $C_{\nu }$ par l'application de restriction (\ref{res g}) est~:
$$\mbox{res} (C_{\nu}) = \sum_{
\begin{footnotesize}
\begin{array}{l}
\alpha_i \subset p_i \times q_i \\
i = 1 , \ldots , m \\
\beta_j \subset P_j \times Q_j \\
j=1 , \ldots , m-2
\end{array}
\end{footnotesize}} \left( \prod_{i=1}^{m-1} c_{\alpha_i \beta_i}^{\beta_{i-1}} \right) C_{\alpha_1} \otimes \ldots \otimes C_{\alpha_m} ,$$
o\`u on a pos\'e $\beta_0 = \nu$ et $\beta_{m-1} = \alpha_m$.
Le lemme \ref{res de C} d\'ecoule alors du lemme \ref{LR}.

\bigskip

Concluons maintenant la d\'emonstration de la proposition \ref{clef}. L'\'equivalence entre les points 1. et 2. d\'ecoule des lemmes
\ref{lemvenky2} et \ref{res de C}. L'\'equivalence de ces deux points avec le point 3. d\'ecoule du lemme suivant d\^u \`a Venkataramana
\cite[Lemma 1.5]{Venky}.

\begin{lem} \label{lemvenky3}
En conservant les notations pr\'ec\'edentes, le noyau de l'application de restriction
$$\mbox{res} : H^* (\hat{X}_G ) \rightarrow H^* (\hat{X}_L )$$
contient la classe de Schubert $C_{\nu}$ si et seulement si $V(\nu ) \subset E(G,L)^{\perp}$, o\`u l'orthogonal est pris par rapport
au produit scalaire induit par la forme de Killing.
\end{lem}

Puisque le $K$-sous-espace $V(\nu )$ de $\bigwedge \mathfrak{p}^+$ est irr\'eductible et de multiplicit\'e $1$, il est
n\'ecessairement inclus dans $E(G,L)$ ou dans son orthogonal. La proposition \ref{clef} d\'ecoule des lemmes \ref{lemvenky2}
et \ref{lemvenky3}.

\bigskip

\subsubsection*{D\'ecomposition ``\`a la Lefschetz''}

Fixons maintenant un sous-groupe de congruence $\Gamma$ dans $G$.
Soit $(\lambda , \mu)$ un couple compatible de partition. Nous noterons $H^{\lambda , \mu} (S(\Gamma )) = H^{|\lambda| +|\hat{\mu}|}
(A(\lambda , \mu ) : \Gamma )$ la $A(\lambda , \mu )$-composante fortement primitive de la cohomologie de $S(\Gamma )$. La formule
de Matsushima et la classification de Vogan-Zuckerman impliquent le th\'eor\`eme suivant.

\begin{thm} \label{VZ}
Soit $\Gamma$ un sous-groupe de congruence dans $G$.
Pour chaque couple d'entiers $(i,j)$ avec $i+j \leq pq$, on a~:
$$H^{i,j} (S(\Gamma )) = \bigoplus_{
\begin{footnotesize}
\begin{array}{l}
(\lambda , \mu) \mbox{ couple compatible} \\
\mbox{de partitions avec} \\
|\lambda |\leq i \ |\hat{\mu}| \leq j
\end{array}
\end{footnotesize}} \bigoplus_{
\begin{footnotesize}
\begin{array}{l}
\nu_i \subset p_i \times q_i \\
i= 1 , \ldots , m
\end{array}
\end{footnotesize}} E^{\lambda , \mu}_{\nu_1 , \ldots , \nu_m } (S(\Gamma)),$$
o\`u le diagramme gauche $\mu/\lambda$ est r\'eunion de diagrammes rectangulaires $p_i \times q_i$ ne s'intersectant
qu'en des sommets et chaque $E^{\lambda , \mu}_{\nu_1 , \ldots , \nu_m } (S(\Gamma))$ est isomorphe \`a $H^{\lambda , \mu} (S(\Gamma))$.
\end{thm}

Cette d\'ecomposition est en g\'en\'eral plus fine que celle induite par l'action des classes de Chern. Ainsi par exemple pour $p=q=3$, si l'on consid\'ere $\lambda_1 =(2,2)$, $\mu_1 =(3,2,1)$, $\lambda_2 =(2,1,1)$ et $\mu_2 = (3,2,1)$. Il existe bien \'evidemment une image entre les tableaux gauches $\mu_1 / \lambda_1$ et $\mu_2 /\lambda_2$ et
nous verrons au cours de la d\'emonstration du th\'eor\`eme suivant que l'action naturelle des  classes de Chern de $S(\Gamma )$ sur la cohomologie ne distingue alors pas le
sous-espace $H^{\lambda_1 ,\mu_1} (S(\Gamma ))$ du sous-espace  $H^{\lambda_2 ,\mu_2} (S(\Gamma ))$. On peut n\'eanmoins montrer le th\'eor\`eme suivant.

\begin{thm} \label{dec de Lefschetz}
Soit $\Gamma$ un sous-groupe de congruence dans $G$ et
soit $\eta : H^* ({\Bbb G}_{p,q} ) \rightarrow H^* (S(\Gamma ))$ l'application d\'efinie en (\ref{eta}). Fixons
$\lambda$, $\mu$ et $\nu$ trois partitions incluses dans $p \times q$ telles que le couple $(\lambda , \mu )$ soit
compatible. Alors,
\begin{enumerate}
\item pour toute classe fortement primitive $s \in H^{\lambda , \mu} (S(\Gamma ))$,
$\eta (C_{\nu} ).s =0$ si et seulement si la partition $\nu$ ne s'inscrit pas dans $\mu / \lambda$, et
\item lorsque le diagramme gauche $\mu / \lambda$ est rectangulaire \'egal \`a $a \times b$ avec $a=p$ ou $b=q$,
$$H^{\lambda , \mu } (S(\Gamma )) = \left\{ s \in H^{|\lambda|, |\hat{\mu}|} (S(\Gamma)) \; : \; \eta (C_{\nu}).s=0, \  \forall
\nu \not\subset a \times b \right\}.$$
\item si $s \in H^{\lambda , \mu } (S( \Gamma ))$ est une classe non nulle, les \'el\'ements
$$\{ C_{\nu} .s \; : \; \nu \subset p\times q , \; \nu \mbox{ s'inscrit dans } \mu / \lambda \}$$
sont lin\'eairement ind\'ependants.
\end{enumerate}
\end{thm}
{\it D\'emonstration.}  Soit $s$ une classe fortement primitive dans $H^{\lambda , \mu} (S(\Gamma ))$. Alors d'apr\`es la formule de Matsushima,
$$s \in \mbox{Hom}_K (V(\lambda , \mu ) , C^{\infty} (\Gamma \backslash G^{{\rm nc}} )).$$
On en d\'eduit que
$$\eta (C_{\nu} ) .s =0 \Leftrightarrow C_{\nu} . V( \lambda , \nu ) =0 \mbox{ dans } \bigwedge \mathfrak{p} .$$
Le point 1. du th\'eor\`eme \ref{dec de Lefschetz} d\'ecoule donc de la proposition \ref{clef}. Le point 3. s'en d\'eduit pareillement.

Supposons maintenant le diagramme gauche $\mu /\lambda$ rectangulaire \'egal \`a $a \times b$ et  consid\'erons une classe $s \in H^{|\lambda| , |\hat{\mu}|} (S(\Gamma ))$ telle que
$\eta (C_{\nu }).s=0$ pour tout diagramme $\nu \not\subset a \times b$. D'apr\`es la formule de Matsushima, on peut supposer que
$s \in \mbox{Hom}_K (V(\alpha , \beta ) , C^{\infty} (\Gamma \backslash G^{{\rm nc}} ))$, pour un certain couple compatible de partitions $(\alpha , \beta )$ v\'erifiant
$|\alpha | = |\lambda|$ et $|\beta|=|\mu|$. Alors le diagramme gauche $\beta / \alpha$ comporte $ab$ cases. S'il est constitu\'e de plusieurs rectangles, l'une des partitions $\nu = (b+1)$ ou $\nu = (1^a)$ s'y inscrit n\'ecessairement et $C_{\nu} . V( \alpha , \beta ) \neq 0$, ce qui est absurde puisque $\nu \not\subset a \times b$ et donc $\eta (C_{\nu} ) .s =0$.
Le diagramme gauche $\beta / \alpha$ est donc rectangulaire et comporte $ab$ cases. De la m\^eme mani\`ere on montre que le diagramme rectangulaire
$\beta / \alpha$ est n\'ecessairement \'egal \`a $a \times b$. Il reste \`a comprendre sa position dans le diagramme $p\times q$. Mais sous l'hypoth\`ese du point 2. du th\'eor\`eme \ref{dec de Lefschetz} celle-ci est d\'etermin\'ee par la connaissance de $|\lambda|$ et de $|\hat{\mu}|$.
Finalement, on a n\'ecessairement $\alpha = \lambda$ et $\beta = \mu$. Ce qui conclut la d\'emonstration du th\'eor\`eme \ref{dec de Lefschetz}.

\bigskip

Le th\'eor\`eme \ref{dec de Lefschetz} contient la d\'ecomposition de Lefschetz usuelle \footnote{Et montre que celle-ci est optimale pour $p=1$, o\`u l'optimalit\'e signifie que
la d\'ecomposition en $K$-types induite par l'action des classes de Chern est une d\'ecomposition en irr\'eductible. Signalons au passage l'article pr\'ecurseur de Chern \cite{Chern}
sur ces questions. Les deux th\'eor\`emes ci-dessus peuvent \^etre interpr\'et\'e comme la compl\`etion du programme de Chern dans le cas des vari\'et\'es localement model\'ees sur
$X_G$.} et implique le corollaire suivant d\^u \`a Parthasarathy \cite[Theorem 2 et 3 et Corollary 2.24]{Parthasarathy}.

\begin{cor} \label{partha}
\begin{enumerate}
\item Le groupe $H^{l,0} (S(\Gamma )) $ est trivial si $l$ n'est pas de la forme $pq-ab$ pour certains entiers $a$ et $b$ tels que $0\leq a \leq p$ et $0\leq b \leq q$.
\item Pour $l\neq pq$, soit $I_l = \left\{ (a,b) \in {\Bbb N}^2 \; : \; 1 \leq a \leq p , \, 1 \leq b \leq q \mbox{ et } pq-ab=l \right\}$. Pour $l=pq$, soit $I_l = \left\{ (0,0) \right\}$.
Soit $H(a,b) = \left\{ s \in H^{pq-ab,0} (S(\Gamma )) \; : \right.$ $\left. C_{a+1} .s = \ldots = C_p .s = 0 \mbox{ et } C_{b+1} ' .s = \ldots = C_{q} '.s =0 \right\}$. Alors,
$H^{l,0} (S(\Gamma )) = \bigoplus_{(a,b) \in I_l} H(a,b)$.
\item Pour toute classe non nulle $s \in H(a,b)$ et $y \in \eta (H^* ({\Bbb G}_{p,q} ))$, $y.s =0$ si et seulement si $y$ appartient \`a l'id\'eal engendr\'e par
$C_{a+1} , \ldots , C_{p}$ et $C_{b+1}' = \ldots = C_q '$.
\item Soit $s$ une classe non nulle dans $H(a,b)$. Alors, les \'el\'ements $\{ C_1^{k_1} C_2^{k_2} \ldots C_a^{k_a} .s \}$ (resp. $\{ (C_1 ')^{k_1} \ldots (C_b ')^{k_b} .s$)
o\`u $C_1^{k_1} C_2^{k_2} \ldots C_a^{k_a}$ (resp. $(C_1 ')^{k_1} \ldots (C_b ')^{k_b}$) parcourt l'ensemble des mon\^omes de degr\'e total $\leq b$ (resp. $\leq a$) en
$C_1, C_2, \ldots , C_a$ (resp. $C_1 ' , \ldots , C_b '$) sont lin\'eairements ind\'ependants.
\end{enumerate}
\end{cor}

\subsection{Sous-vari\'et\'es de Shimura de $Sh^0 G$}

Notre principal but est de comprendre l'application de restriction stable de la cohomologie de $Sh^0 G$ vers une sous-vari\'et\'e de Shimura.
Dans cette section nous classifions les vari\'et\'es de Shimura pouvant apparaitre comme sous-vari\'et\'e de $Sh^0 G$.

\medskip

Commen\c{c}ons par d\'ecrire quelques sous-groupes r\'eels du groupe $U(p,q)$.

Le groupe $U(p,q)$ pr\'eserve la forme hermitienne standard $h(x,y) = \sum |x_{\mu }|^2 - \sum |y_{\nu}|^2$ sur la somme directe ${\Bbb C}^p \oplus {\Bbb C}^q$.
Soient $a,b,i , j$ des entiers strictement positifs. Supposons $ia+jb \leq p$ et $ib+ja \leq q$. L'espace
${\Bbb C}^p \oplus {\Bbb C}^q$ se d\'ecompose en une somme directe
$$\underbrace{({\Bbb C}^a \oplus {\Bbb C}^b ) \oplus \ldots \oplus ({\Bbb C}^a \oplus {\Bbb C}^b )}_{i \ {\rm fois}} \oplus
\underbrace{({\Bbb C}^b \oplus {\Bbb C}^a ) \oplus \ldots \oplus ({\Bbb C}^b \oplus {\Bbb C}^a )}_{j \ {\rm fois}} \oplus
({\Bbb C}^{p-ia-jb} \oplus {\Bbb C}^{q-ib-ja} ) .$$

On r\'ealise alors le groupe $U(a,b)$ comme le sous-groupe du groupe $U(p,q)$ qui pr\'eserve chacun des $i$ sous-espaces ${\Bbb C}^a \oplus {\Bbb C}^b$ en agissant de
mani\`ere standard, pr\'eserve  les $j$ sous-espaces ${\Bbb C}^b \oplus {\Bbb C}^a$ en agissant comme son conjugu\'e complexe apr\`es permutation des
facteurs ${\Bbb C}^a$ et ${\Bbb C}^b$ et agit
trivialement sur ${\Bbb C}^{p-ia-jb} \oplus {\Bbb C}^{q-ib-ja}$. Pour simplifier, nous noterons ce plongement de $U(a,b)$ dans $U(p,q)$~:
\begin{eqnarray} \label{UdansU}
g \in U(a,b) \mapsto (\underbrace{g, \ldots ,g}_{i \ {\rm fois}} , \underbrace{\tilde{g} , \ldots , \tilde{g}}_{j \ {\rm fois}} , id ) .
\end{eqnarray}

Rappelons que le groupe $GSp_a$ est le sous-groupe de $U(a,a)$ d\'efini par~:
$$\left\{ g = \left(
\begin{array}{cc}
A & B \\
C & D
\end{array} \right) \in U(a, a ) \; : \; {}^t g \left(
\begin{array}{cc}
0 & 1_{a} \\
-1_{a} & 0
\end{array} \right) g = \left(
\begin{array}{cc}
0 & 1_{a} \\
-1_{a} & 0
\end{array} \right) \right\} .$$
L'application compos\'ee
\begin{eqnarray} \label{SdansU}
g \in GSp_a \mapsto g \in U(a,a) \mapsto (\underbrace{g , \ldots ,g}_{i \ {\rm fois}} ,\underbrace{\tilde{g} , \ldots , \tilde{g}}_{j \ {\rm fois}} , id )
\end{eqnarray}
r\'ealise donc un plongement de $GSp_a$ dans $U(p,q)$.

Rappelons enfin que le groupe $O^*(2a)$ est le sous-groupe de $U(a,a)$ d\'efini par~:
$$\left\{ g = \left(
\begin{array}{cc}
A & B \\
C & D
\end{array} \right) \in U(a, a ) \; : \; {}^t g \left(
\begin{array}{cc}
0 & 1_{a} \\
1_{a} & 0
\end{array} \right) g = \left(
\begin{array}{cc}
0 & 1_{a} \\
1_{a} & 0
\end{array} \right) \right\} .$$
L'application compos\'ee
\begin{eqnarray} \label{O*dansU}
g \in O^*(2a) \mapsto g \in U(a,a) \mapsto (\underbrace{g , \ldots ,g}_{i \ {\rm fois}} ,\underbrace{\tilde{g} , \ldots , \tilde{g}}_{j \ {\rm fois}} , id )
\end{eqnarray}
r\'ealise donc un plongement de $O^*(2a)$ dans $U(p,q)$.

Soient $p_1 , \ldots , p_m , q_1 , \ldots , q_m$ des entiers strictement positifs tels que
$p_1 + \ldots + p_m \leq p$ et $q_1 + \ldots + q_m \leq q$. On r\'ealise alors le groupe $U(p_1 , q_1 ) \times \ldots \times U(p_m , q_m )$
comme le sous-groupe du groupe $U(p,q)$ qui pr\'eserve les sous-espaces $E_i \oplus F_i$ pour chaque $i$
et agit trivialement sur l'orthogonal de la somme directe de tous ces sous-espaces. \`A conjugaison par un \'el\'ement de $U(p,q)$ pr\`es, le plongement
ainsi d\'efini ne d\'epend pas du choix des sous-espaces $E_i$, $F_i$. Pour simplifier, nous le noterons~:
\begin{eqnarray} \label{UidansU}
(g_1 , \ldots , g_m ) \in \prod_j U(p_j ,q_j ) \mapsto (g_1 , \ldots , g_m , id ) \in U(p,q) .
\end{eqnarray}

Remarquons que les plongements (\ref{UdansU}), (\ref{SdansU}) et (\ref{O*dansU}) de $GSp_a$ v\'erifient bien la condition (\ref{cond2}).

\begin{prop} \label{ssvar}
Soit $Sh^0 H$ une sous-vari\'et\'e de Shimura de $Sh^0 G$. Rappelons que, par d\'efinition, $H$ est alors
un sous-groupe de $G$ qui est alg\'ebrique r\'eductif connexe d\'efini sur ${\Bbb Q}$ et qui
v\'erifie (\ref{cond1}) et (\ref{cond2}). Pour simplifier, nous supposerons $H$ maximal parmi tous les sous-groupes de $G$ donnant la
m\^eme sous-vari\'et\'e de Shimura. Notons $H^{{\rm nc}}$ le produit des facteurs non compact de $H({\Bbb R})$.
Alors, le groupe $H^{{\rm nc}}$ est un produit de groupes $U(a_j , b_j )$,  $GSp_{a_j}$ et $O^* (2a_j )$ et, \`a conjugaison par un \'el\'ement de
$G^{{\rm nc}}$ pr\`es, son plongement dans $G^{{\rm nc}} = U(p,q)$ s'obtient en composant au produit direct des morphismes (\ref{UdansU}), (\ref{SdansU}) et
(\ref{O*dansU}) le morphisme (\ref{UidansU}).
\end{prop}
{\it D\'emonstration.} On voit dor\'enavant le ${\Bbb R}$-groupe $G^{{\rm nc}} = U(p,q)$ comme $U_{p+q} ({\Bbb R}, h)$, o\`u $h$ est la forme hermitienne standard non
d\'eg\'en\'er\'ee et de signature $(p,q)$ sur ${\Bbb C}^{p+q}$.

Soit $A$ l'alg\`ebre sur ${\Bbb R}$~: $M_{p+q} ({\Bbb C})$ munie de l'involution $\tau (X) = {}^t \overline{X}$. Le centre de $A$ est le corps ${\Bbb C}$ et puisque
la restriction de $\tau$ \`a celui-ci est non triviale, l'involution $\tau$ est dite de deuxi\`eme esp\`ece.
L'alg\`ebre complexifi\'ee  $A \otimes_{{\Bbb R}} {\Bbb C} = M_{p+q} ({\Bbb C}) \oplus M_{p+q} ({\Bbb C})$ munie de l'involution $\tilde{\tau} (X,Y) = ({}^t Y , {}^t X)$.
Le groupe $GL_{p+q}$ s'identifie \`a $\{ Z \in A \otimes_{{\Bbb R}} {\Bbb C} \; : \; Z \tilde{\tau} (Z) = 1_{p+q} \}$ et (cf. \cite{PlatonovRapinchuk}) Aut$(GL_{p+q} )$ s'identifie
naturellement au groupe des automorphismes de l'alg\`ebre $A \otimes_{{\Bbb R}} {\Bbb C}$ qui commutent \`a l'involution $\tilde{\tau}$.
Le groupe $G^{{\rm nc}}$ est une forme r\'eelle du groupe $GL_{p+q}$ obtenue en le tordant par un cocycle $a \in H^1 ({\Bbb R} , {\rm Aut}_{{\Bbb C}} (GL_{p+q} ))$.
Consid\'erons $a$ comme un cocycle dans $H^1 ({\Bbb R} , {\rm Aut}_{{\Bbb C}} (A \otimes_{{\Bbb R}} {\Bbb C} ))$, nous pouvons alors construire
l'alg\`ebre tordue ${}_a (A \otimes_{{\Bbb R}} {\Bbb C})$. Puisque l'image par le cocyle $a$ de la conjugaison complexe (g\'en\'erateur de Gal$({\Bbb C} / {\Bbb R} )$)
est un automorphisme qui commute \`a l'involution $\tilde{\tau}$, l'involution de ${}_a (A \otimes_{{\Bbb R}} {\Bbb C})$ d\'eduite de $\tilde{\tau}$ commute
\`a la conjugaison complexe. Notons alors $B$ la sous-alg\`ebre r\'eelle de ${}_a (A \otimes_{{\Bbb R}} {\Bbb C})$ fix\'ee par la conjugaison complexe.
L'alg\`ebre $B$ est n\'ecessairement simple sur ${\Bbb C}$ avec une involution de seconde esp\`ece. Plus pr\'ecisemment,
$B= M_{p+q} ({\Bbb C} )$ munie de l'involution tordue $\sigma (X) = H^{-1} {}^t \overline{X} H$, o\`u $H$ est la matrice de $h$ dans la base canonique de ${\Bbb C}^{p+q}$.
Et,
\begin{eqnarray} \label{G}
G^{{\rm nc}} = \left\{ X  \in B \; : \; \sigma (X )X = 1_{p+q}  \right\} .
\end{eqnarray}

Le sous-groupe r\'eel $H^{{\rm nc}}$ de $G^{{\rm nc}}$ s'obtient donc en tordant un sous-groupe alg\'ebrique $R$ de $GL_{p+q}$ invariant par la conjugaison complexe.
Il correspond \`a un tel sous-groupe la sous-alg\`ebre de $A \otimes_{{\Bbb R}} {\Bbb C}$ pr\'eserv\'ee par Aut$R$ et \'equipp\'ee de l'involution induite
par $\tilde{\tau}$. Cette alg\`ebre \`a involution est n\'ecessairement une somme directe de facteurs $(C, \nu)$ isomorphent (comme alg\`ebre \`a un involution)
\`a l'un des trois types suivant d'alg\`ebres (cf. \cite{PlatonovRapinchuk})~:
\begin{enumerate}
\item $C=M_{k} ({\Bbb C})$, $\nu (X) = {}^t X$;
\item $C=M_{k} ({\Bbb C})$ avec $k$ pair et $\nu (X) = J_k {}^t X J_k^{-1}$, o\`u $J_k = \left(
\begin{array}{cc}
0 & 1_{k/2} \\
-1_{k/2} & 0
\end{array} \right)$;
\item $C=M_k ({\Bbb C}) \oplus M_k ({\Bbb C})$, $\nu (X,Y) = ({}^t Y , {}^t X)$.
\end{enumerate}
Dans la suite, nous supposerons pour simplifier qu'il n'y a qu'un seul facteur $(C, \nu )$. Le cas g\'en\'eral se traiterait pareillement et explique l'existence des plongements
(\ref{UidansU}).

La torsion de l'alg\`ebre $A \otimes_{{\Bbb R}} {\Bbb C}$ par le cocycle $a$ induit une torsion ${}_a C$ de l'alg\`ebre $C$. Notons $D$ la sous-alg\`ebre r\'eelle de
${}_a C$ fix\'ee par la conjugaison complexe. L'alg\`ebre $D$ est r\'eelle, simple sur son centre (\'egal \`a ${\Bbb R}$ ou ${\Bbb C}$) et \'equipp\'ee d'une involution $\theta$.
Donc (cf. \cite{PlatonovRapinchuk}) l'alg\`ebre \`a involution $(D, \theta )$ est l'une des formes suivantes~:
\begin{enumerate}
\item $D = M_{k} ({\Bbb R} )$, $\theta (X) = P {}^t X P^{-1}$ et $P \in GL_k ({\Bbb R})$ sym\'etrique;
\item $D = M_k ({\Bbb R} )$ avec $k$ pair, $\theta (X) = P {}^t X P^{-1}$ et $P \in GL_k ({\Bbb R})$ antisym\'etrique;
\item $D = M_{k/2} ({\Bbb H} )$ avec $k$ pair, $\theta (X) = P {}^t \overline{X} P^{-1}$ et $P \in GL_{k/2} ({\Bbb H})$ hermitienne;
\item $D = M_{k/2} ({\Bbb H} )$ avec $k$ pair, $\theta (X) = P {}^t \overline{X} P^{-1}$ et $P \in GL_{k/2} ({\Bbb H})$ antihermitienne;
\item $D = M_k ({\Bbb C} )$, $\theta (X) = P {}^t \overline{X} P^{-1}$ et $P \in GL_k ({\Bbb C})$ hermitienne.
\end{enumerate}
(Ici ${\Bbb H}$ d\'esigne l'alg\`ebre des quaternions sur ${\Bbb R}$.)

Les deux premiers cas sont exclus puisque l'espace sym\'etrique associ\'e \`a $H$ doit \^etre hermitien et v\'erifier (\ref{cond2}).

Maintenant, les seules repr\'esentations irr\'eductibles de l'alg\`ebre $M_{k} ({\Bbb C})$ sont la repr\'esentation standard ${\Bbb C}^k$ et sa
conjugu\'ee. \`A conjugaison pr\`es les seuls plongements de l'alg\`ebre
$M_k ({\Bbb C})$ dans $M_{p+q} ({\Bbb C})$ sont diagonaux par blocs et de la forme suivante~:
\begin{eqnarray} \label{plgt}
X \in M_k ({\Bbb C}) \mapsto \left(
\begin{array}{ccccccc}
0 &  &  &  &  &  & \\
  & X &  &  &  &  & \\
  &    & \ddots & & & & \\
 &  &  & X & & &  \\
 &  &  &  & \overline{X} & &  \\
 &  &  &  &  & \ddots & \\
 & & & & & & \overline{X}
\end{array} \right) \in M_{p+q} ({\Bbb C}) ,
\end{eqnarray}
o\`u la matrice $X$ (resp. $\overline{X}$) apparait $i$ (resp. $j$) fois sur la diagonale pour certains entier $i$ et $j$ de somme $i+j \in [1, (p+q)/k ]$.

(Le seul plongement d'alg\`ebre de $M_k ({\Bbb C}) \oplus M_l ({\Bbb C})$ dans $M_{k+l} ({\Bbb C})$ est le plongement diagonal par blocs. Lorsque $C$ a plusieurs facteurs, les diff\'erents plongements de $C$ dans $M_{p+q} ({\Bbb C})$ sont donc obtenus en plongeant chaque facteur de $C$ diagonalement par bloc
et \`a l'aide d'un plongement du type (\ref{plgt}). Rappelons que pour simplifier, nous supposons que $C$ n'a qu'un seul facteur.)

Dans la suite, nous travaillerons toujours \`a conjugaison pr\`es sans le pr\'eciser.

Lorsque l'alg\`ebre $C =M_k ({\Bbb C})$ est du premier type ci-dessus, elle se plonge n\'ecessairement dans $A \otimes_{{\Bbb R}} {\Bbb C} = M_{p+q} ({\Bbb C}) \oplus M_{p+q} ({\Bbb C})$ comme
$$X \mapsto \left(
\begin{array}{ccccccc}
0 &  &  &  &  &  & \\
  & X &  &  &  &  & \\
  &    & \ddots & & & & \\
 &  &  & X & & &  \\
 &  &  &  & \overline{X} & &  \\
 &  &  &  &  & \ddots & \\
 & & & & & & \overline{X}
\end{array} \right)
\oplus
\left(
\begin{array}{ccccccc}
0 &  &  &  &  &  & \\
  & X &  &  &  &  & \\
  &    & \ddots & & & & \\
 &  &  & X & & &  \\
 &  &  &  & \overline{X} & &  \\
 &  &  &  &  & \ddots & \\
 & & & & & & \overline{X}
\end{array} \right) ,$$
o\`u $X$ (resp. $\overline{X}$) apparait $i$ (resp. $j$) fois sur chaque diagonale.
Dans ce cas l'alg\`ebre r\'eelle $D$ d\'eduite de $C$ est du type 3 ci-dessus (on a d\'ej\`a vu qu'elle ne peut pas \^etre du type 1). L'entier $k$ est donc
n\'ecessairement pair \'egal \`a $2a$ et le groupe r\'eel $H^{{\rm nc}}$ est alors $O^* (2a)$.
Le plongement correspondant de $O^* (2a)$ dans $U(p,q)$ est le plongement (\ref{O*dansU}).

Lorsque l'alg\`ebre $C=M_k  ({\Bbb C})$ est du deuxi\`eme type (alors $k$ est pair), elle se plonge n\'ecessairement dans $A \otimes_{{\Bbb R}} {\Bbb C} =M_{p+q} ({\Bbb C}) \oplus M_{p+q} ({\Bbb C})$ comme
$$X \mapsto \left(
\begin{array}{ccccccc}
0 &  &  &  &  & &  \\
 & X &  &  &  & &  \\
 &  & \ddots &  &  & &  \\
 &  &  & X & & & \\
 &  &  &  & \overline{X} & & \\
 & & & & & \ddots & \\
 & & & & & & \overline{X}
\end{array} \right) $$
$$ \oplus \left(
\begin{array}{ccccccc}
0 &  &  &  &  & &  \\
 & J_k X J_k^{-1} &  &  &  & &  \\
 &  & \ddots  &  &  & &  \\
 &  &  & J_k X J_k^{-1} & & &  \\
 &  &  &  & J_k \overline{X} J_k^{-1} & & \\
 & & & & & \ddots & \\
 & & & & & & J_k \overline{X} J_k^{-1}
\end{array} \right) ,$$
o\`u $X$ et $J_k X J_k^{-1}$ (resp. $\overline{X}$ et $J_k \overline{X} J_k^{-1}$) apparaissent $i$ (resp. $j$) fois chacun sur les diagonales respectives.
Dans ce cas l'alg\`ebre $D$ d\'eduite de $C$ est du type 4 ci-dessus (on a d\'ej\`a vu qu'elle ne peut pas \^etre du type 2). L'entier $k$ est n\'ecessairement \'egal \`a
$2a$ et le groupe r\'eel $H^{{\rm nc}}$ est alors $GSp_a$. Le plongement correspondant dans $U(p,q)$ est alors le plongement
(\ref{SdansU}).

Enfin lorsque l'alg\`ebre $C=M_k ({\Bbb C}) \oplus M_k ({\Bbb C})$ (troisi\`eme type), elle se plonge n\'ecessairement dans $A \otimes_{{\Bbb R}} {\Bbb C} = M_{p+q} ({\Bbb C}) \oplus M_{p+q} ({\Bbb C})$ comme
$$X \oplus Y \mapsto
\left(
\begin{array}{ccccccc}
0 &  &  &  &  &  & \\
  & X &  &  &  &  & \\
  &    & \ddots & & & & \\
 &  &  & X & & &  \\
 &  &  &  & \overline{X} & &  \\
 &  &  &  &  & \ddots & \\
 & & & & & & \overline{X}
\end{array} \right)
\oplus
\left(
\begin{array}{ccccccc}
0 &  &  &  &  &  & \\
  & Y &  &  &  &  & \\
  &    & \ddots & & & & \\
 &  &  & Y & & &  \\
 &  &  &  & \overline{Y} & &  \\
 &  &  &  &  & \ddots & \\
 & & & & & & \overline{Y}
\end{array} \right) ,$$
o\`u $X$ (resp. $Y$) apparait $i$ fois sur la premi\`ere (resp. deuxi\`eme) diagonale et
$\overline{X}$ (resp. $\overline{Y}$) apparait $j$ fois sur la premi\`ere (resp. deuxi\`eme) diagonale.
Dans ce cas l'alg\`ebre $D$ d\'eduite de $C$ est du type 5 ci-dessus.
Le groupe r\'eel $H^{{\rm nc}}$ est alors $U(a,b)$ pour deux entiers $a$, $b$ de somme $k$. Le plongement correspondant dans $U(p,q)$ est alors le plongement
(\ref{UdansU}).

La proposition \ref{ssvar} s'en d\'eduit en mettant bout \`a bout ces diff\'erents cas.

\bigskip

Afin d'appliquer les r\'esultats de la section pr\'ec\'edente, nous aurons besoin de d\'eterminer la classe duale $[\hat{X}_H ]$ dans $H^* (\hat{X}_G )$ associer
\`a chaque sous-vari\'et\'e de Shimura $Sh^0 H$ de $Sh^0 G$. Ce calcul ne faisant, en fait, appel qu'aux groupes r\'eels, d'apr\`es la proposition \ref{ssvar}
il nous suffira de d\'eterminer $[\hat{X}_H]$ dans chacun des cas suivant~:
\begin{enumerate}
\item $H^{{\rm nc}} = U(p_1 ,q_1 ) \times \ldots \times U(p_m ,q_m )$ avec $p_j , q_j \geq 1$, $p_1 + \ldots + p_m \leq p$ et $q_1 + \ldots + q_m \leq q$;
\item $H^{{\rm nc}} = GSp_p$ et $p=q$;
\item $H^{{\rm nc}} = O^* (2p)$ et $p=q$.
\end{enumerate}
C'est l'objet des trois propositions suivantes.

\begin{prop} \label{dualU}
Soient $p_j$, $q_j$ avec $j=1, \ldots, m$ des entiers strictement positifs tels que $p_1 + \ldots + p_m \leq p$ et $q_1 + \ldots + q_m \leq q$.
Supposons que $H^{{\rm nc}} = U(p_1 , q_1 ) \times \ldots \times U(p_m , q_m )$. Alors, \`a un multiple scalaire non nul pr\`es, la classe duale
$$[\hat{X}_H ] = \sum_{\nu \subset p\times q} c^{\nu}_{p_1 \times q_1 \ldots p_m \times q_m} C_{\hat{\nu}} \in H^{pq-\sum_i p_i q_i} ({\Bbb G}_{p,q} ).$$
\end{prop}
{\it D\'emonstration.} Notons $C$ la classe duale $[\hat{X}_H ]$ \`a la sous-vari\'et\'e $\hat{X}_H = {\Bbb G}_{p_1 , q_1} \times \ldots \times {\Bbb G}_{p_m , q_m} $ de
$\hat{X}_G = {\Bbb G}_{p,q}$. Un classe
$C_{\nu} \in H^* ({\Bbb G}_{p,q} )$ v\'erifie~:
\begin{eqnarray} \label{cc}
C_{\nu} . C = \mbox{res}(C_{\nu} ) \wedge C ,
\end{eqnarray}
o\`u res$:H^* ({\Bbb G}_{p,q} ) \rightarrow H^* (\hat{X}_H )$ est l'application naturelle de restriction.
En particulier, $C_{\nu} . C \neq 0$ si et seulement si res$(C_{\nu} ) \neq 0$. Ce qui d'apr\`es le lemme \ref{lemvenky2} et la proposition \ref{clef}
est \'equivalent au fait que $\nu$ s'inscrive dans le diagramme gauche $(p_1\times q_1 ) * \ldots * (p_m \times q_m )$.

Rappelons (cf. par exemple \cite{Fulton}) que la classe duale de $C_{\nu}$ dans $H^* ({\Bbb G}_{p,q} )$ est \'egale \`a $C_{\hat{\nu}}$.
La classe $C$ est donc \'egale \`a une combinaison lin\'eaire \`a coefficients non nuls des classes $C_{\hat{\nu}}$ o\`u $\nu$ d\'ecrit
l'ensembles des partitions images du diagramme gauche $(p_1 \times q_1 ) * \ldots * (p_m \times q_m )$ (on a n\'ecessairement
$|\nu | = p_1 q_1 + \ldots + p_m q_m$).

Finalement le lemme \ref{res de C} implique que si $|\nu | = p_1 q_1 + \ldots + p_m q_m$,
$$\mbox{res}(C_{\nu} ) = c^{\nu}_{p_1 \times q_1 \ldots p_m \times q_m}  C_{p_1 \times q_1} \otimes \ldots \otimes C_{p_m \times q_m} \in H^* (\hat{X}_H ) = \bigotimes_{i=1}^m
H^* ({\Bbb G}_{p_i , q_i } ).$$
\`A l'aide de (\ref{cc}) on constate alors qu'\`a un multiple scalaire non nul pr\`es, on doit donc n\'ecessairement avoir~:
$$C= \sum_{\nu \subset p\times q} c^{\nu}_{p_1 \times q_1 \ldots p_m \times q_m} C_{\hat{\nu}} \in H^{pq-\sum_i p_i q_i} ({\Bbb G}_{p,q} ).$$
Ce qui conclut la d\'emonstration de la proposition \ref{dualU}.

\bigskip

\begin{prop} \label{dualS}
Supposons $p=q$ et $H^{{\rm nc}} = GSp_p$. Alors, \`a un multiple scalaire non nul pr\`es, la classe duale
$$[\hat{X}_H ] = C_{\nu} \in H^{\frac{p(p-1)}{2}} ({\Bbb G}_{p,p}),$$
o\`u $\nu$ est la partition $(p-1 , p-2 , \ldots , 1)$ de diagramme de Young
$$p-1 \ {\rm cases}  \; \left\{
\begin{tabular}{l}
\begin{tabular}{|c|c|c|c|c|} \hline
 & & & &\\ \hline
\end{tabular}
\\
\begin{tabular}{|c|c|c|c|} \hline
 & & & \\ \hline
\end{tabular}
\\
 \begin{tabular}{|c|c|c|} \hline
 & & \\ \hline
\end{tabular}
\\
\begin{tabular}{|c|c|} \hline
 & \\ \hline
\end{tabular}
\\
\begin{tabular}{|c|} \hline
 \\ \hline
\end{tabular}
\end{tabular} \right. $$
\end{prop}
{\it D\'emonstration.} Supposons $p=q$, alors $G^{{\rm nc}} = U(p,p)$ et (d'apr\`es la proposition \ref{ssvar}) on peut supposer
que $H^{{\rm nc}} = GSp_p$ est donn\'e par~:
$$\left\{ g = \left(
\begin{array}{cc}
A & B \\
C & D
\end{array} \right) \in U(p,p) \; : \; {}^t g \left(
\begin{array}{cc}
0 & 1_p \\
-1_p & 0
\end{array} \right) g = \left(
\begin{array}{cc}
0 & 1_p \\
-1_p & 0
\end{array} \right) \right\} .$$
Alors, le groupe $K \cap H^{{\rm nc}} = U(p)$ se plonge dans $K$ par l'application~:
$$g \mapsto \left(
\begin{array}{cc}
g & 0 \\
0 & {}^t g^{-1}
\end{array} \right) ,$$
et l'action de $U(p)$ sur $\mathfrak{p}^+ = {\Bbb C}^p \otimes ({\Bbb C}^p )^*$ est isomorphe \`a la repr\'esentation $\rho \otimes \rho$, o\`u $\rho$ est la repr\'esentation
standard de $U(p)$ sur ${\Bbb C}^p$. Via le plongement de $H^{{\rm nc}}$ dans $G^{{\rm nc}}$, $\mathfrak{p}^+ \cap \mathfrak{h}$ se trouve identifi\'e avec
$$\mbox{sym}^2 (\rho ) = \mbox{sym}^2 ({\Bbb C}^p) \subset {\Bbb C}^p \otimes ({\Bbb C}^p )^*.$$
(On identifie sym$^2 ({\Bbb C}^p )$ avec l'espace des matrices sym\'etriques.)

Fixons maintenant une partition $\lambda$ de poids $|\lambda | = \frac{p(p+1)}{2}$. (Remarquons que $\frac{p(p+1)}{2}$ est la dimension de $\mathfrak{p}^+ \cap \mathfrak{h}$.)
Notons $D^{\otimes k}$ la repr\'esentation de dimension un de $U(p)$ donn\'ee par le d\'eterminant \`a la puissance $k$.

\begin{lem} \label{res qd H=S}
L'image de $V(\lambda ) (\subset \bigwedge \mathfrak{p}^+ )$ par l'application de restriction dans $\bigwedge (\mathfrak{p}^+ \cap \mathfrak{h} )$ est nulle sauf si
$\lambda = (p, p-1 , \ldots , 1)$ auquel cas elle est isomorphe au $U(p)$-module $D^{\otimes (p+1)}$.
\end{lem}
{\it D\'emonstration.} Vue comme repr\'esentation de $U(p) (= K \cap H^{{\rm nc}} \subset K)$ le module
$$V(\lambda ) = ({\Bbb C}^p )^{\lambda } \otimes ({\Bbb C}^p )^{\lambda^*} .$$
Il est classique (cf. par exemple \cite{Fulton}) que sa d\'ecomposition en irr\'eductibles est~:
\begin{eqnarray} \label{deco en irred}
V(\lambda ) = \bigoplus_{\nu} [ ({\Bbb C}^p )^{\nu} ]^{\oplus c_{\lambda \lambda^*}^{\nu} }.
\end{eqnarray}

L'espace $\mathfrak{p}^+ \cap \mathfrak{h} =$ sym$^2 ({\Bbb C}^p )$ s'identifie avec l'espace des matrices sym\'etriques, il est engendr\'e par les matrices
$E_{i,j} + E_{j,i} $ ($1\leq i,j \leq p$). On constate donc que si $\nu \not\subset p\times (p+1)$,
la partition $\nu$ ne peut pas \^etre le poids d'un vecteur du $U(p) (= K \cap H^{{\rm nc}})$-module $\bigwedge (\mathfrak{p}^+ \cap \mathfrak{h} )$.
En particulier, si $\nu \not\subset p\times (p+1)$, $\nu$ ne peut pas \^etre le plus haut poids d'un sous-module irr\'eductible de $\bigwedge (
\mathfrak{p}^+ \cap \mathfrak{h}  )$.

Puisque dans la d\'ecomposition (\ref{deco en irred}) chaque partition $\nu$ donnant lieu \`a un sous-module non trivial doit avoir un poids $|\nu| = |\lambda |+ |\lambda^* | =p(p+1)$.
L'image de chacun de ces sous-modules par l'application de restriction dans $\bigwedge (\mathfrak{p}^+ \cap \mathfrak{h} )$ est nulle sauf \'eventuellement le sous-module
$$[ ({\Bbb C}^p )^{\nu_0 } ]^{\oplus c_{\lambda \lambda^*}^{\nu_0}} ,$$
o\`u $\nu_0 = p\times (p+1)$.

Il nous faut donc comprendre pour quelles partitions $\lambda$, le coefficient de Littlewood-Richardson $c_{\lambda \lambda^*}^{p\times (p+1)}$
est non nul.  Supposons donc que la partition $\lambda$ v\'erifie que $c_{\lambda \lambda^*}^{p\times (p+1)}\neq 0$. Commen\c{c}ons par remarquer que
$\lambda$ et $\lambda^*$ sont alors n\'ecessairement contenues dans $p\times (p+1)$. Ensuite remarquons que le calcul de Schubert implique si $\lambda $ et $\mu$
sont deux partitions contenues dans une partition rectangulaire $a\times b$ et de somme des poids $|\lambda |+ |\mu | =ab$ alors~:
$$c_{\lambda , \mu}^{a\times b} = \delta_{\mu , \hat{\lambda}} .$$
(Ici $\delta$ est le symbole de Kronecker \'egal \`a un si ses deux arguments sont identiques, \`a z\'ero sinon.)

Revenons maintenant au coefficient $c_{\lambda \hat{\lambda}}^{p\times (p+1)}$. D'apr\`es ce que nous venons de voir, le coefficient $c_{\lambda \hat{\lambda}}^{p\times (p+1)}$
est non nul si et seulement si $c_{\lambda \hat{\lambda}}^{p\times (p+1)}=1$ ce qui est \'equivalent \`a ce que $\lambda^* = \hat{\lambda}$ (ici la notation $\hat{\lambda}$
d\'esigne la partition compl\'ementaire de $\lambda$ dans le rectangle $p\times (p+1)$). Ceci se traduit en
$$\lambda_{p+1-i} + \max \{ j \; : \; \lambda_j \geq i \} = p+1 ,$$
pour tout entier $i=1, \ldots , p$. D'o\`u l'on d\'eduit $\lambda_i = p+1-i$ pour $i=1 , \ldots ,p$ et $\lambda = (p , p-1 , \ldots , 1)$.

Nous avons donc d\'emontr\'e que l'image de $V(\lambda )$ par l'application de restriction dans $\bigwedge (\mathfrak{p}^+ \cap \mathfrak{h} )$
est nulle sauf si $\lambda = (p,p-1 , \ldots ,1)$. Montrons maintenant que lorsque $\lambda = (p,p-1 , \ldots ,1)$, l'image de $V(\lambda )$ dans
$\bigwedge (\mathfrak{p}^+ \cap \mathfrak{h})$ est non nulle.
Supposons donc $\lambda = (p,p-1, \ldots ,1)$. Il est clair que le vecteur
$$\bigwedge_{j=1}^p \bigwedge_{i=1}^j e_i \otimes f_j^* \in V(\lambda ) .$$
(Appliquer la permutation $(1 \; 2 \; \ldots \; p)$ sur les colonnes de $M_p ({\Bbb C}) = E \otimes F^*$, via l'action d'un \'el\'ement de $K$.)

\noindent
L'image de ce vecteur dans  $\bigwedge (\mathfrak{p}^+ \cap \mathfrak{h} ) = \bigwedge \mbox{sym}^2 ({\Bbb C}^p )$ est un multiple non nul de
$$\bigwedge_{1\leq i\leq j\leq p} (E_{i,j} + E_{j,i} )$$
qui engendre un $U(p)$-module isomorphe au module $D^{\otimes (p+1)}$. Ce qui conclut la d\'emonstration du lemme \ref{res qd H=S}.

\medskip

Continuons la d\'emonstration de la proposition \ref{dualS}.
Notons $C$ la classe duale $[\hat{X}_H ]$ \`a la sous-vari\'et\'e $\hat{X}_H$ de $\hat{X}_G = {\Bbb G}_{p,q}$.
D'apr\`es le lemme \ref{lemvenky2} (qui reste valable dans ce cas, cf. \cite{Venky}) une classe $C_{\nu} \in H^* ({\Bbb G}_{p,q} )$
v\'erifie $C_{\nu} . C \neq 0$ si et seulement si res$(C_{\nu}) \neq 0$ o\`u res$: H^* ({\Bbb G}_{p,q} ) \rightarrow H^* (\hat{X}_H )$ est l'application naturelle de restriction.
Mais chaque classe $C_{\nu}$ est $K$-invariante et donc sa restriction \`a $\hat{X}_H$ est $(K \cap H^{{\rm nc}} = U(p))$-invariante. Il suffit donc de v\'erifier res$(C_{\nu}) \neq 0$
au point base. Ce qui nous ram\`ene au lemme \ref{res qd H=S}.

On conclut que si $\nu \subset p\times p$ est une partition de poids $|\nu | = \frac{p(p+1)}{2}$,
$$C_{\nu} . C \neq 0 \Leftrightarrow \nu = (p, p-1 , \ldots , 1) .$$
La classe $C$ est donc n\'ecessairement \'egale \`a un multiple scalaire non nul de $C_{\hat{\nu}}$. Ce qui d\'emontre la proposition \ref{dualS}.

\bigskip

\begin{prop} \label{dualO}
Supposons $p=q$ et $H^{{\rm nc}} = O^* (2p)$. Alors, \`a un multiple scalaire non nul pr\`es, la classe duale
$$[\hat{X}_H ] = C_{\nu} \in H^{\frac{p(p+1)}{2}} ({\Bbb G}_{p,p} ), $$
o\`u $\nu$ est la partition $(p, p-1 , \ldots , 1)$ de diagramme de Young
$$p \ {\rm cases}  \; \left\{
\begin{tabular}{l}
\begin{tabular}{|c|c|c|c|c|} \hline
 & & & &\\ \hline
\end{tabular}
\\
\begin{tabular}{|c|c|c|c|} \hline
 & & & \\ \hline
\end{tabular}
\\
 \begin{tabular}{|c|c|c|} \hline
 & & \\ \hline
\end{tabular}
\\
\begin{tabular}{|c|c|} \hline
 & \\ \hline
\end{tabular}
\\
\begin{tabular}{|c|} \hline
 \\ \hline
\end{tabular}
\end{tabular} \right. $$
\end{prop}
{\it D\'emonstration.} Supposons $p=q$, alors $G^{{\rm nc}} = U(p,p)$ et (d'apr\`es la proposition \ref{ssvar}) on peut supposer que $H^{{\rm nc}} = O^* (2p)$ est donn\'e par~:
$$\left\{ g = \left(
\begin{array}{cc}
A & B \\
C & D
\end{array} \right) \in U(p,p) \; : \; {}^t g \left(
\begin{array}{cc}
0 & 1_p \\
1_p & 0
\end{array} \right) g = \left(
\begin{array}{cc}
0 & 1_p \\
1_p & 0
\end{array} \right) \right\} .$$
Alors, le groupe $K \cap H^{{\rm nc}} = U(p)$ se plonge dans $K$ par l'application~:
$$g \mapsto \left(
\begin{array}{cc}
g & 0 \\
0 & {}^t g^{-1}
\end{array} \right) ,$$
et l'action de $U(p)$ sur $\mathfrak{p}^+ = {\Bbb C}^p \otimes ({\Bbb C}^p )^*$ est isomorphe \`a la repr\'esentation $\rho \otimes \rho$, o\`u $\rho$ est la repr\'esentation
standard de $U(p)$ sur ${\Bbb C}^p$. Via le plongement de $H^{{\rm nc}}$ dans $G^{{\rm nc}}$, $\mathfrak{p}^+ \cap \mathfrak{h}$ se trouve identifi\'e avec
$$\bigwedge^2 \rho = \bigwedge^2 {\Bbb C}^p \subset {\Bbb C}^p \otimes ({\Bbb C}^p )^*.$$
(On identifie $\bigwedge^2 {\Bbb C}^p $ avec l'espace des matrices antisym\'etriques.)

Fixons maintenant une partition $\lambda$ de poids $|\lambda | = \frac{p(p-1)}{2}$. (Remarquons que $\frac{p(p-1)}{2}$ est la dimension de $\mathfrak{p}^+ \cap \mathfrak{h}$.)
La d\'emonstration du lemme \ref{res qd H=S} se traduit facilement en une d\'emonstration du lemme suivant.

\begin{lem} \label{res qd H=O}
L'image de $V(\lambda ) (\subset \bigwedge \mathfrak{p}^+ )$ par l'application de restriction dans $\bigwedge (\mathfrak{p}^+ \cap \mathfrak{h} )$ est nulle sauf si
$\lambda = (p-1 , p-2 , \ldots , 1)$ auquel cas elle est isomorphe au $U(p)$-module $D^{\otimes (p)}$.
\end{lem}

On conclut alors la d\'emonstration de la proposition \ref{dualO} en suivant mot pour mot la fin de la d\'emonstration de la proposition \ref{dualS}.

\subsection{Restriction stable \`a une sous-vari\'et\'e de Shimura}

Dans cette section nous d\'emontrons des crit\`eres d'injectivit\'e de l'application de restriction
stable (\ref{Res}) aux diff\'erentes sous-vari\'et\'es de Shimura de $Sh^0 G$. Les crit\`eres seront ind\'ependant des ${\Bbb Q}$-structures des groupes alg\'ebriques en question.
La proposition \ref{ssvar} permet alors de r\'eduire le probl\`eme aux trois types suivants de sous-vari\'et\'es de Shimura $Sh^0 H \subset Sh^0 G$~:
\begin{enumerate}
\item $H^{{\rm nc}} = U(p_1 ,q_1 ) \times \ldots \times U(p_m ,q_m )$ avec $p_j , q_j \geq 1$, $p_1 + \ldots + p_m \leq p$ et $q_1 + \ldots + q_m \leq q$;
\item $H^{{\rm nc}} = GSp_p$ et $p=q$;
\item $H^{{\rm nc}} = O^* (2p)$ et $p=q$.
\end{enumerate}
C'est l'objet des deux th\'eor\`emes qui suivent ainsi que d'une remarque finale concernant le cas $H^{{\rm nc}} = O^* (2p)$.

\begin{thm} \label{res stable U}
Soit $Sh^0 H$ une sous-vari\'et\'e de Shimura de $Sh^0 G$ avec $H^{{\rm nc}} = U(p_1 ,q_1 ) \times \ldots \times U(p_m ,q_m )$, $p_j , q_j \geq 1$, $p_1 + \ldots + p_m \leq p$
et $q_1 + \ldots + q_m \leq q$. Soient $\lambda$ et $\mu$ deux partitions incluses dans $p\times q$ telles que le couple $(\lambda , \mu )$ soit compatible. Alors, l'application
$${\rm Res}_H^G : H^* (Sh^0 G) \rightarrow \prod_{G({\Bbb Q} )} H^* (Sh^0 H )$$
de restriction stable est {\bf injective en restriction \`a $H^{\lambda , \mu} (Sh^0 G)$} s'il existe une partition $\nu$ image du diagramme gauche $(p_1 \times q_1 ) *
\ldots * (p_m \times q_m )$ telle que $\hat{\nu}$ s'inscrive dans le diagramme gauche $\mu / \lambda$.
\end{thm}
{\it D\'emonstration.} La d\'emonstration repose sur un th\'eor\`eme de Venkataramana, sur la proposition \ref{clef} et sur la proposition \ref{dualU}.
Plus pr\'ecisemment, d'apr\`es \cite[Theorem 6]{Venky2}, si $s$ est un classe dans $H^{\lambda , \mu} (Sh^0 G)$ dont la restriction stable Res$_H^G (s)$ \`a $Sh^0 H$ est triviale, alors
$$[\hat{X}_H ] . V(\lambda , \mu ) =0 \mbox{ dans } \bigwedge \mathfrak{p} .$$

D'apr\`es la proposition \ref{dualU},
$$[\hat{X}_H ] = \sum_{\nu \subset p\times q} c^{\nu}_{p_1 \times q_1 \ldots p_m \times q_m} C_{\hat{\nu}} \in H^{pq-\sum_i p_i q_i} ({\Bbb G}_{p,q} ).$$

Puisque d'apr\`es la proposition \ref{clef}, $C_{\hat{\nu}} . V(\lambda , \mu ) \neq 0$ si et seulement si $\hat{\nu}$ s'inscrit dans le diagramme gauche $\mu / \lambda$ et que les
\'el\'ements $\{ C_{\hat{\nu}} . v(\lambda ) \otimes w (\hat{\mu} )^* \}$, o\`u $\hat{\nu}$ parcourt l'ensemble des partitions qui s'inscrivent dans le diagramme gauche
$\mu /  \lambda$, sont lin\'eairement ind\'ependants, on conclut que $[\hat{X}_H ] . V(\lambda , \mu ) \neq 0$ si et seulement s'il existe une partition $\nu \subset p\times q$
image du diagramme gauche $(p_1 \times q_1 ) * \ldots * (p_m \times q_m )$ (autrement dit $c^{\nu}_{p_1 \times q_1 \ldots p_m \times q_m}\neq 0$) telle que
$\hat{\nu}$ s'inscrive dans le diagramme gauche $\mu / \lambda$. Ce qui cl\^ot la d\'emonstration du th\'eor\`eme \ref{res stable U}.

\bigskip

Rappelons que le sous-espace $H^{\lambda , \mu} (Sh^0 G)$ apparait dans la cohomologie holomorphe si et seulement si $\mu = p\times q$. La partition $\lambda$
est alors naturellement param\'etr\'ee par un couple d'entier $(r,s)$ avec $0 \leq r \leq p$ et $0 \leq s \leq q$ tels que
$$\lambda = ( \underbrace{q, \ldots , q}_{r \ {\rm fois}} , \underbrace{s , \ldots , s}_{p-r \ {\rm fois}} )$$
de diagramme de Young~:
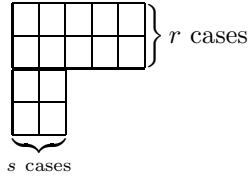
$$
\begin{array}{l}
\left. \hspace{0,035cm}
\begin{array}{|c|c|c|c|c|} \hline
 & & & &\\ \hline
 & & & &\\ \hline
\end{array}
\right\}  r \ {\rm cases} \\
\underbrace{
\begin{array}{|c|c|} \hline
 &  \\ \hline
 &  \\ \hline
\end{array}}_{s \ {\rm cases}}
\end{array}
$$
(Ici $p=4$ et $q=5$.)

Dans ce cas (et pour souligner le fait qu'il apparait dans la cohomologie holomorphe) nous noterons $H^{(r,s),0} (Sh^0 G)$ le sous-espace $H^{\lambda , \mu} (Sh^0 G)$ de
la cohomologie holomorphe de degr\'e $|\lambda | =  rq + s(p-r)$ (remarquons que $|\hat{\mu}|=0$).

\begin{cor}[Clozel-Venkataramana] \label{CV}
Soit $Sh^0 H$ une sous-vari\'et\'e de Shimura de $Sh^0 G$ avec $H^{{\rm nc}} = U(p_1 , q_1 ) \times \ldots \times U(p_m , q_m )$, $p_j , q_j \geq 1$, $p_1 + \ldots + p_m \leq p$
et $q_1 + \ldots + q_m \leq q$. Soit $(r,s)$ un couple d'entiers naturels avec $r\leq p$ et $s\leq q$. Alors, l'application
$${\rm Res}_H^G : H^* (Sh^0 G) \rightarrow \prod_{G({\Bbb Q} )} H^* (Sh^0 H )$$
de restriction stable est {\bf injective} en restriction \`a $H^{(r,s), 0} (Sh^0 G)$ {\bf si et seulement si} soit $p_1 + \ldots +p_m =p$,
$r=0$ et $s \leq q_i$ pour chaque $i$, soit $q_1 + \ldots + q_m = q$, $s=0$ et $r\leq p_i$ pour chaque $i$.
\end{cor}
{\it D\'emonstration.} Commen\c{c}ons par montrer que la condition est suffisante. Ici le diagramme gauche $\mu / \lambda$
est en fait le diagramme rectangulaire $(p-r)(q-s)$. Une partition $\alpha$ s'inscrit donc dans $\mu / \lambda$ si et seulement si $l(\alpha ) \leq p-r$ et $l(\alpha^* ) \leq q-s$.

Quitte \`a r\'eordonner les facteurs $U(p_i , q_i )$ on peut supposer $q_1 \leq \ldots \leq q_m$.
Il est alors clair que la partition
$$\nu_0 = (\underbrace{q_m , \ldots , q_m}_{p_m \ {\rm fois}} , \ldots , \underbrace{q_1 , \ldots , q_1}_{p_1\  {\rm fois}} )$$
est une image du diagramme gauche $(p_1 \times q_1 ) * \ldots * (p_m \times q_m )$. Or,
$$l(\hat{\nu}_0 ) = \left\{
\begin{array}{ll}
p - p_m & \mbox{ si } q_m = q \\
p & \mbox{ sinon}
\end{array} \right.$$
et
$$l(\hat{\nu}_0^* ) = \left\{
\begin{array}{ll}
q & \mbox{ si } p_1 + \ldots + p_m < p \\
q- q_1 & \mbox{ sinon} .
\end{array} \right. $$
Donc si $p_1 + \ldots + p_m = p$, $r=0$ et $s \leq q_i$ pour chaque $i$, la partition $\nu_0$ s'inscrit dans $\mu / \lambda$ et, d'apr\`es le th\'eor\`eme \ref{res stable U},
l'application de restriction stable Res$_H^G$ est donc injective en restriction \`a $H^{(r,s) ,0} (Sh^0 G)= H^{\lambda , \mu } (Sh^0 G)$.

On d\'emontre de m\^eme que Res$_H^G$ est injective en restriction \`a $H^{(r,s), 0} (Sh^0 G)$ lorsque $q_1 + \ldots + q_m = q$, $s=0$ et $r \leq p_i$ pour chaque $i$.

Montrons maintenant que ces conditions suffisantes sont en fait n\'ecessaires. Il nous faut pour cela d'abord comprendre plus pr\'ecisemment
$l(\hat{\nu} )$ et $l(\hat{\nu}^*)$ lorsque $\nu$ est une image du diagramme gauche $(p_1 \times q_1) * \ldots *(p_m \times q_m )$.
C'est l'objet du lemme suivant.

\begin{lem} \label{lnu}
Soient $p_1 , \ldots , p_m , q_1 , \ldots ,q_m$ des entiers $\geq 1$.
Soit $\nu$ une partition image du diagramme gauche $(p_1 \times q_1 ) * \ldots * (p_m \times q_m )$. Alors,
\begin{enumerate}
\item $l(\nu ) \leq p_1 + \ldots +p_m$,
\item $\nu_{p_1 + \ldots + p_m} \leq \min q_i$, et
\item si $\nu_{p_1 + \ldots +p_m} \neq 0$ alors $\nu_1 <q_1 + \ldots + q_m$.
\end{enumerate}
\end{lem}
{\it D\'emonstration du lemme \ref{lnu}.} Pour simplifier nous supposerons
$q_1 \leq \ldots \leq q_m$.
Il suffit alors de montrer que toute partition $\nu$ image du diagramme gauche $(p_1 \times q_1) * \ldots * (p_m \times q_m)$ v\'erifie
\begin{itemize}
\item $l(\nu ) \leq p_1 + \ldots + p_m$,
\item $\nu_{p_1 + \ldots + p_m} \leq q_1$, et
\item si $\nu_{p_1 + \ldots +p_m} \neq 0$ alors $\nu_1 <q_1 + \ldots + q_m$.
\end{itemize}
Puisque la longueur de la partition sous-jacente au produit de deux tableaux de Young est toujours inf\'erieure \`a la somme des longueurs des partitions
sous-jacentes \`a ces deux tableaux, le premier point se d\'emontre imm\'ediatement par r\'ecurrence sur $m$.

Montrons les deux derniers points \'egalement par r\'ecurrence sur $m$. Le cas $m=1$ est trivial. On passe facilement de $m$ \`a $m+1$
en remarquant que l'on peut supposer $p_{m+1} =1$ puisque le diagramme rectangulaire $p_{m+1} \times q_{m+1}$ est bien \'evidemment une
image du diagramme gauche $\underbrace{(1 \times q_{m+1}) * \ldots * (1 \times q_{m+1} )}_{p_{m+1} \  {\rm fois}}$.
Il existe alors une partition $\alpha$ image du diagramme gauche $(p_1 \times q_1 ) * \ldots * (p_m \times q_m )$ telle que
$\nu$ soit une image du diagramme gauche $\alpha * (1 \times q_{m+1} )$. Le diagramme $\nu$ est alors obtenu en ajoutant $q_{m+1}$ cases
au diagramme $\alpha$ sans jamais en mettre deux dans la m\^eme colonne. Or d'apr\`es l'hypoth\`ese de r\'ecurrence,
$\alpha_{p_1 + \ldots + p_m} \leq q_1$ et si $\alpha_{p_1 + \ldots + p_m} \neq 0$ alors $\alpha_1 < q_1 + \ldots +q_m$.
On ne peut donc pas rajouter plus de $q_1$ cases \`a la $(p_1 + \ldots + p_m +1)$-i\`eme ligne (qui est vide) du diagramme $\alpha$.
Donc $\nu_{p_1 + \ldots +p_m +1} \leq q_1$. Enfin, on montre de m\^eme que si $\nu_{p_1 + \ldots + p_m +1} \neq 0$ alors $\nu_1 < q_1 + \ldots + q_m$.
Ce qui conclut la r\'ecurrence et la d\'emonstration du lemme \ref{lnu}.

\medskip

Continuons la d\'emonstration du corollaire, d'apr\`es \cite[Proposition 2.3]{ClozelVenky} le crit\`ere du th\'eor\`eme de Venkataramana
que l'on a utilis\'e dans la d\'emonstration du th\'eor\`eme \ref{res stable U} est, dans le cas de la cohomologie {\bf holomorphe}, \'egalement n\'ecessaire et \'equivalent \`a ce que
$E(G,H) \supset V(\lambda )$. D'apr\`es la proposition \ref{clef} cette derni\`ere assertion est \'equivalente au fait que $\lambda$ s'inscrive
dans le diagramme gauche $(p_1 \times q_1 ) * \ldots * (p_m \times q_m )$. Le lemme \ref{lnu} appliqu\'e tour \`a tour \`a $\lambda$
et \`a son conjugu\'e implique alors que soit $p_1 + \ldots + p_m = p$, $r=0$ et $s \leq \min q_i$
soit $q_1 + \ldots + q_m =q$, $s=0$ et $r \leq \min p_i$. Ce qui conlut la d\'emonstration du corollaire.

\bigskip

Remarquons que notre d\'emonstration du r\'esultat de Clozel et  Venkataramana est local alors que la d\'emonstration de \cite[Proposition 3.A.10]{ClozelVenky}
utilise un argument global.

\medskip

Le lemme suivant est imm\'ediat. Conjugu\'e au th\'eor\`eme \ref{res stable U} il donne un crit\`ere plus simple \`a v\'erifier pour l'injectivit\'e de l'application de restriction stable Res$_H^G$ en restriction aux sous-espaces consid\'er\'es dans le th\'eor\`eme.

\begin{lem} \label{critere simple}
Tout diagramme de Young obtenu en assemblant des blocs rectangulaires $p_i \times q_i$ est une image du diagramme gauche $(p_1 \times q_1 ) * \ldots * (p_m \times q_m )$.
\end{lem}

Remarquons n\'eanmoins que l'on n'obtient pas toutes les images du diagramme gauche $(p_1 \times q_1 ) * \ldots * (p_m \times q_m)$ de cette mani\`ere. Par exemple, la partition $(3,1)$ est une image du diagramme gauche $(1 \times 2) * (1 \times 2)$.

\bigskip

La d\'emonstration du th\'eor\`eme \ref{res stable U} en rempla\c{c}ant la proposition \ref{dualU} par la proposition \ref{dualS} implique le th\'eor\`eme suivant.

\begin{thm} \label{res stable S}
Soit $Sh^0 H$ une sous-vari\'et\'e de Shimura de $Sh^0 G$ avec $H^{{\rm nc}} = GSp_p$ et $p=q$. Soient $\lambda$ et $\mu$ deux partitions incluses dans $p\times p$ telles
que le couple $(\lambda , \mu )$ soit compatible. Alors, l'application
$${\rm Res}_H^G : H^* (Sh^0 G) \rightarrow \prod_{G({\Bbb G})} H^* (Sh^0 H)$$
de restriction stable est {\bf injective en restriction \`a $H^{\lambda , \mu} (Sh^0 G)$} si la partition $\nu =  (p-1 , p-2 , \ldots , 1)$ de diagramme de Young
$$p-1 \ {\rm cases}  \; \left\
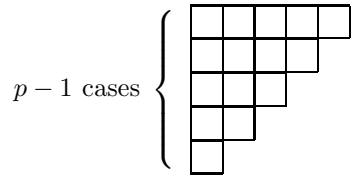
{
\begin{tabular}{l}
\begin{tabular}{|c|c|c|c|c|} \hline
 & & & &\\ \hline
\end{tabular}
\\
\begin{tabular}{|c|c|c|c|} \hline
 & & & \\ \hline
\end{tabular}
\\
 \begin{tabular}{|c|c|c|} \hline
 & & \\ \hline
\end{tabular}
\\
\begin{tabular}{|c|c|} \hline
 & \\ \hline
\end{tabular}
\\
\begin{tabular}{|c|} \hline
 \\ \hline
\end{tabular}
\end{tabular} \right. $$
s'inscrit dans le diagramme gauche $\mu / \lambda$.
\end{thm}

Et comme au-dessus, on en d\'eduit le corollaire suivant.

\begin{cor}[Clozel-Venkataramana]
Soit $Sh^0 H$ une sous-vari\'et\'e de Shimura de $Sh^0 G$ avec $H^{{\rm nc}} = GSp_p$ et $p=q$.
Soit $(r,s)$ un couple d'entiers naturels avec $r\leq p$ et $s\leq q$. Alors, l'application
$${\rm Res}_H^G : H^* (Sh^0 G) \rightarrow \prod_{G({\Bbb G})} H^* (Sh^0 H)$$
de restriction stable est {\bf injective} en restriction \`a $H^{(r,s),0} (Sh^0 G)$ {\bf si et seulement si} $r,s \leq 1$.
\end{cor}

\bigskip

Remarquons enfin que d'apr\`es le lemme \ref{lnu} la partition $(p, p-1 , \ldots ,1)$ ne peut \^etre image d'un diagramme gauche
$(p_1 \times q_1 ) * \ldots * (p_m \times q_m)$ avec $p_1 + \ldots + p_m \leq p$ et $q_1 + \ldots + q_m\leq p$. Dans le cas o\`u $H^{{\rm nc}} =O^* (2p)$, le crit\`ere d'injectivit\'e de l'application de restriction stable que l'on obtiendrait en suivant la m\^eme m\'ethode serait vide, ou presque puisqu'il permet quand m\^eme de retrouver le r\'esultat suivant de
Clozel et Venkataramana~:

\medskip

\noindent
{\it Soit $Sh^0 H$ une sous-vari\'et\'e de Shimura de $Sh^0 G$ avec $H^{{\rm nc}} = O^* (2p)$ et $p=q$.
Alors, l'application
$${\rm Res}_H^G : H^* (Sh^0 G) \rightarrow \prod_{G({\Bbb Q})} H^* (Sh^0 H)$$
de restriction stable est identiquement {\bf nulle} en restriction \`a la cohomologie {\bf holomorphe} de degr\'e strictement positif.}

\bigskip

Concluons cette section par une r\'eciproque partielle au th\'eor\`eme \ref{res stable U}, un crit\`ere d'annulation de l'application de restriction.

\begin{thm} \label{res nulle}
Soient $Sh^0 H$ une sous-vari\'et\'e de Shimura de $Sh^0 G$ avec $H = H_1 \times \ldots \times H_m$, $H_i^{{\rm nc}} = U(p_i , q_i )$, $p_i , q_i \geq 1$, $p_1 + \ldots + p_m \leq p$
et $q_1 + \ldots + q_m \leq q$. Soient $\lambda$ et $\mu$ deux partitions incluses dans $p\times q$ telles que le couple $(\lambda , \mu )$ soit compatible.
Fixons pour chaque entier $i=1, \ldots ,m$, un couple compatible de partitions $(\lambda_i , \mu_i )$ dans $p_i \times q_i$. Supposons
$$c_{\lambda_1 \ldots \lambda_m}^{\lambda} c_{\hat{\mu}_1 \ldots \hat{\mu}_m}^{\hat{\mu}} = 0.$$
(Autrement dit, supposons soit que $\lambda$ n'est pas une image du diagramme gauche $\lambda_1 * \ldots * \lambda_m$ soit que $\hat{\mu}$ n'est pas une image du
diagramme gauche $\hat{\mu}_1 * \ldots * \hat{\mu}_m$.)

Alors, la projection de l'image de l'application de restriction
$${\rm res}_H^G : H^{\lambda , \mu} (Sh^0 G) \rightarrow H^{|\lambda | + |\hat{\mu}|} (Sh^0 H)$$
dans la composante de K\"unneth
$$H^{\lambda_1 , \mu_1} (Sh^0 H_1 ) \otimes \ldots \otimes H^{\lambda_m , \mu_m} (Sh^0 H_m )$$
est {\bf nulle}.
\end{thm}
{\it D\'emonstration.} Remarquons que l'on peut supposer $|\lambda_1 | + \ldots + |\lambda_m | = |\lambda |$ et $|\hat{\mu}_1| + \ldots + |\hat{\mu}_m| =|\hat{\mu}|$.
Supposons par exemple $c_{\lambda_1 \ldots \lambda_m}^{\lambda} =0$ (la d\'emonstration serait similaire dans le cas
$c_{\hat{\mu}_1 \ldots \hat{\mu}_m}^{\hat{\mu}}$). D'apr\`es la formule de Matsushima, si $s$ est une classe de cohomologie,
$$s \in \mbox{Hom}_K (V(\lambda , \mu ), C^{\infty} (G({\Bbb Q}) \backslash G({\Bbb A}))) ,$$
o\`u, rappelons le, $V(\lambda , \mu ) \subset \bigwedge^{|\lambda|} \mathfrak{p}^+ \otimes \bigwedge^{|\hat{\mu}|} \mathfrak{p}^-$.

Notons $\mathfrak{p}_H^{\pm}$ les espaces $\mathfrak{p}^{\pm} \cap \mathfrak{h}$.

La restriction res$_H^G (s)$ de $s$ \`a la sous-vari\'et\'e $Sh^0 H$ est (repr\'esent\'ee par) une forme diff\'erentielle ferm\'ee qui appartient \`a
$${\rm Hom}_K ( \bigwedge^{|\lambda|} \mathfrak{p}_H^+ \otimes \bigwedge^{|\hat{\mu}|} \mathfrak{p}_H^- , C^{\infty} (H({\Bbb Q} ) \backslash H({\Bbb A}))),$$
et est obtenue par la formule~:
$${\rm res}_H^G (s) ( \xi ) = s( {\rm pr}_{V(\lambda , \mu )} (\xi )) \quad (\mbox{pour tout } \xi \in \bigwedge^{|\lambda|} \mathfrak{p}_H^+ \otimes \bigwedge^{|\hat{\mu}|} \mathfrak{p}_H^- )$$
o\`u pr$_{V(\lambda , \mu)} : \bigwedge^{|\lambda|} \mathfrak{p}_H^+ \otimes \bigwedge^{|\hat{\mu}|} \mathfrak{p}_H^-  \rightarrow V(\lambda , \mu )$ est une projection $K$-\'equivariante.

La formule de K\"unneth correspond \`a la d\'ecomposition de $\bigwedge^{|\lambda|} \mathfrak{p}_H^+ \otimes \bigwedge^{|\hat{\mu}|} \mathfrak{p}_H^-$ en~:
$$\bigwedge^{|\lambda|} \mathfrak{p}_H^+ \otimes \bigwedge^{|\hat{\mu}|} \mathfrak{p}_H^- = \bigoplus_{
\begin{array}{l}
k_1 + \ldots + k_m = |\lambda| \\
l_1 + \ldots + l_m = |\hat{\mu}|
\end{array}} (\bigwedge^{k_1} \mathfrak{p}_{H_1}^+ \otimes \bigwedge^{l_1} \mathfrak{p}_{H_1}^- ) \otimes \ldots \otimes (\bigwedge^{k_m} \mathfrak{p}_{H_m}^+
\otimes \bigwedge^{l_m} \mathfrak{p}_{H_m}^- ) .$$

L'espace $V(\lambda_1 , \mu_1 ) \otimes \ldots \otimes V(\lambda_m , \mu_m )$ est inclus dans $\bigwedge^{|\lambda|} \mathfrak{p}_H^+ \otimes \bigwedge^{|\hat{\mu}|} \mathfrak{p}_H^- $. Le th\'eor\`eme \ref{res nulle} d\'ecoulera de la d\'emonstration de
\begin{eqnarray} \label{pr}
{\rm pr}_{V(\lambda , \mu )} (V(\lambda_1 , \mu_1 ) \otimes \ldots \otimes V(\lambda_m , \mu_m )) = 0.
\end{eqnarray}

Mais pour chaque entier $i = 1, \ldots ,m$, le
sous-espace $V(\lambda_i , \mu_i )$ est inclus dans
$$V(\lambda_i ) \otimes V(\hat{\mu}_i )^* \subset \bigwedge^{|\lambda_i |} \mathfrak{p}_{H_i}^+ \otimes \bigwedge^{|\hat{\mu}_i|} \mathfrak{p}_{H_i}^- .$$
Soit $F(\lambda_1 , \ldots , \lambda_m )$ le $K$-module engendr\'e par $V(\lambda_1 ) \otimes \ldots \otimes V(\lambda_m ) \subset \bigwedge^{|\lambda|} \mathfrak{p}^+$.
Pour d\'emontrer (\ref{pr}), il nous suffit alors de d\'emontrer que~:
\begin{eqnarray} \label{KHom1}
{\rm Hom}_K ( V(\lambda ) , F(\lambda_1 , \ldots , \lambda_m ) )=0 .
\end{eqnarray}
Puisque le $K$-module $V(\lambda)$ et le $K_H$-module $V(\lambda_1 ) \otimes \ldots \otimes V(\lambda_m )$ sont tous deux irr\'eductibles, (\ref{KHom1})
est \'equivalent \`a~:
\begin{eqnarray} \label{KHom2}
{\rm Hom}_{K_H} (V(\lambda_1 ) \otimes \ldots \otimes V(\lambda_m ) , V(\lambda )) = 0 .
\end{eqnarray}
Or lorsque $m=2$ et d'apr\`es \cite[(20) p.122]{Fulton}, vu comme $(K_H = K_{H_1} \times K_{H_2})$-module $V(\lambda )$ contient un sous-module irr\'eductible isomorphe au
module $V(\lambda_1 ) \otimes V(\lambda_2)$ si et seulement si $c_{\lambda_1 \lambda_2}^{\lambda} \neq 0$. On conclut alors la d\'emonstration du
th\'eor\`eme \ref{res nulle} par r\'ecurrence sur $m$ et \`a l'aide du lemme \ref{LR}.

\bigskip

\subsection{Exemples et applications}

Concluons cette premi\`ere partie par l'\'etude de quelques familles d'exemples. L'application de restriction stable de la cohomologie de $Sh^0 G$ vers la cohomologie de
ses sous-vari\'et\'es de Shimura est simple \`a comprendre en petit degr\'e en raison du fait \'evident suivant.

\begin{fait} \label{fait}
{\rm Soit $(\lambda , \mu )$ un couple compatible de partitions incluses dans $p\times q$. Supposons le diagramme gauche
$\mu/\lambda = (p_1 \times q_1 )* \ldots *(p_m \times q_m )$. Supposons $|\lambda|+|\hat{\mu}| < 3p-2$ si $q=p$ et $|\lambda|+|\hat{\mu}| <p+q-1$ si $q>p$.
Alors (l'entier $m$ est n\'ecessairement $\leq 2$ et) soit $p_1 + \ldots + p_m =p$, soit $q_1 + \ldots + q_m = q$, soit $p=q$, $m=1$ et $p_1 = q_1 = p-1$.}
\end{fait}

En effet, remarquons que puisque chaque $p_i$ et chaque $q_i$ est $\geq 1$,
si $p_1 + \ldots + p_m \leq p-1$ et $q_1 + \ldots + q_m \leq q-1$ nous devons avoir
$p_1 q_1 + \ldots + p_m q_m \leq (p-m)(q-m)+ m-1$ et donc $|\lambda| +|\hat{\mu}| \geq pq-(p-m)(q-m)-m+1 = m(p+q) -m^2 -m+1$.
Il est alors facile de v\'erifier que $|\lambda|+|\hat{\mu}| (\geq m(p+q)-m^2 -m+1) \geq p+q-1$ avec \'egalit\'e si et seulement si $m=1$, $p_1 = p-1$ et $q_1 = q-1$. Enfin,
si $p=q$ et $m>1$ ou $p_1 <p-1$ ou $q_1 <p-1$, $|\lambda|+|\hat{\mu}| (\geq 2mp-m^2 -m+1) \geq 3p-2$. Ce qui conclut la d\'emonstration du fait.

\bigskip

On d\'eduit facilement du fait \ref{fait} et des th\'eor\`emes \ref{res stable U} et \ref{res stable S} la d\'emonstration du th\'eor\`eme \ref{tentative u(p,q)} que  l'on reformule de la fa\c{c}on suivante.

\medskip

\noindent
{\bf Th\'eor\`eme \ref{tentative u(p,q)}'} {\it Soient $K$ un corps de nombre totalement r\'eel et $G'$ un groupe unitaire  anisotrope sur $K$ compact \`a toutes les places infinies
sauf une o\`u il est isomorphe \`a $U(p,q)$ avec $1 \leq p \leq q$. Soit $G$ le groupe alg\'ebrique sur ${\Bbb Q}$ obtenu \`a partir de
$G'$ par restriction des scalaires de $K$ \`a ${\Bbb Q}$. Alors,
pour tout sous-espace de cohomologie fortement primitive $H^{\lambda , \mu} (Sh^0 G)$, de degr\'e $|\lambda| + |\hat{\mu}| < 3p-2$ si $p=q$ et
$|\lambda|+|\hat{\mu}| < p+q-1$ si $p<q$, de la
vari\'et\'e de Shimura $Sh^0 G$, il existe une sous-vari\'et\'e de Shimura $Sh^0 H \subset Sh^0 G$ telle que l'application de
restriction stable
$${\rm Res}_H^G : H^k (Sh^0 G)  \rightarrow H^k (Sh^0 H)$$
soit injective en restriction \`a $H^{\lambda , \mu } (Sh^0 H)$.}

\medskip

\noindent
{\it D\'emonstration.} Supposons comme d'habitude le diagramme gauche $\mu / \lambda$ \'egal \`a $(p_1 \times q_1) * \ldots * (p_m \times q_m )$.
Sous  les hypoth\`eses du th\'eor\`eme sur le degr\'e $|\lambda| + |\hat{\mu}| $, et d'apr\`es le fait \ref{fait}, soit $p_1 + \ldots + p_m= p$ et dans dans ce cas d'apr\`es le th\'eor\`eme \ref{res stable U}
on peut prendre n'importe quelle sous-vari\'et\'e $Sh^0 H$ avec $H^{{\rm nc}} = U(p,q-1)$, soit $q_1 + \ldots + q_m =q$ et dans dans ce cas d'apr\`es le th\'eor\`eme \ref{res stable U}
on peut prendre n'importe quelle sous-vari\'et\'e $Sh^0 H$ avec $H^{{\rm nc}} = U(p-1 , q)$, soit
$p=q$ et le diagramme gauche $\mu / \lambda$ est rectangulaire \'egal \`a $(p-1) \times (q-1)$
et dans dans ce cas d'apr\`es le th\'eor\`eme \ref{res stable S} on peut prendre n'importe quelle sous-vari\'et\'e $Sh^0 H$ avec $H^{{\rm nc}} = GSp_p$.

\bigskip

Remarquons avec Venkataramana que l'on ne peut esp\`erer r\'epondre positivement \`a la question d'Arthur simplement \`a coup de restriction \`a des sous-vari\'et\'es de
Shimura. On peut plus g\'en\'eralement montrer le r\'esultat suivant qui montre qu'en un certain sens le th\'eor\`eme \ref{tentative u(p,q)} est optimal.

\begin{prop} \label{negative}
Soient $K$ un corps de nombre totalement r\'eel et $G'$ un groupe unitaire  anisotrope sur $K$ compact \`a toutes les places infinies
sauf une o\`u il est isomorphe \`a $U(p,q)$ avec $1 \leq p \leq q$. Soit $G$ le groupe alg\'ebrique sur ${\Bbb Q}$ obtenu \`a partir de
$G'$ par restriction des scalaires de $K$ \`a ${\Bbb Q}$. Alors, il existe une classe de cohomologie (holomorphe) {\bf non triviale} et
de degr\'e $3p-2$ si $p=q$ et $p+q-1$ si $p<q$, dont la restriction stable \`a n'importe quelle sous-vari\'et\'e de Shimura $Sh^0 H \subset Sh^0 G$ soit {\bf nulle}.
\end{prop}
{\it D\'emonstration.} Lorsque $q>p$, il suffit de consid\'erer une classe de cohomolologie holomorphe non nulle dans
$H^{(1,1),0} (Sh^0 G)$, une telle classe existe d'apr\`es un th\'eor\`eme d'Anderson \cite{Anderson} (dans \cite{Li} Li g\'en\'eralise amplement ce th\'eor\`eme \footnote{Pr\'ecisemment, Li montre que $H^{\lambda, \mu} (Sh^0 G)$ est non
trivial lorsque le diagramme gauche $\mu / \lambda$ est rectangulaire.}).
On a classifi\'e toutes les sous-vari\'et\'es de Shimura possible de $Sh^0 G$ (cf. proposition \ref{ssvar}) et d'apr\`es les r\'esultats de Clozel et Venkataramana red\'emontr\'es
dans la section pr\'ec\'edente, la restriction stable de cette classe \`a n'importe quelle sous-vari\'et\'e de Shimura est nulle.

Lorsque $q=p$, il suffit par exemple de consid\'erer une classe de cohomologie holomorphe non nulle dans $H^{(2,1),0} (Sh^0 G)$. L\`a encore,
une telle classe existe d'apr\`es \cite{Anderson}. D'apr\`es les r\'esultats de Clozel et Venkataramana et la classification des sous-vari\'et\'es de Shimura de $Sh^0 G$,
la restriction stable de cette classe \`a n'importe quelle sous-vari\'et\'e de Shimura est nulle.

\bigskip

Remarquons maintenant qu'\`a partir d'une id\'ee de Venkataramana, les th\'eor-\\\`emes \ref{res stable U} et \ref{res nulle} permettent de d\'emontrer le th\'eor\`eme \ref{hodge u(p,q)}
dont nous rappelons l'\'enonc\'e ci-dessous.

\medskip

\noindent
{\bf Th\'eor\`eme \ref{hodge u(p,q)}} {\it
Soient $K$ un corps de nombre totalement r\'eel et $G'$ un groupe unitaire anisotrope sur $K$ compact \`a toutes les places infinies
sauf une o\`u il est isomorphe \`a $U(p,q)$ avec $3 \leq 2p+1 \leq q$. Soit $G$ le groupe alg\'ebrique sur ${\Bbb Q}$ obtenu \`a partir de
$G'$ par restriction des scalaires de $K$ \`a ${\Bbb Q}$.
Alors, toute classe de Hodge dans $H^{2pq-2p} (Sh^0 G)$ est alg\'ebrique.}

\medskip

\noindent
{\it D\'emonstration.} Commen\c{c}ons par remarquer que si $q\geq 2p+1$, la cohomologie de bidegr\'e $(p,p)$ de $Sh^0 G$ est somme directe de
la cohomologie invariante $\oplus_{|\nu| = p} {\Bbb C} \eta (C_{\nu} )$ (o\`u $\nu$ d\'esigne une partition) et de la cohomologie fortement primitive
$H^{(p\times 1) , (p\times (q-1))} (Sh^0 G)$.

Il est classique que les classes $\eta (C_{\nu} )$ sont alg\'ebriques. Nous ne nous occuperons donc que de la partie non invariante $H^*_{{\rm n.i.}} (Sh^0 G)$ (fortement primitive) de la
cohomologie.

Consid\'erons une sous-vari\'et\'e de Shimura $Sh^0 H \subset Sh^0 G$ avec
$H= H_1 \times \ldots \times H_p$ et chaque $H_i^{{\rm nc}} = U(1,2)$ (une telle sous-vari\'et\'e existe puisque $G$ provient d'un groupe unitaire sur un corps de nombre).
Puisque $q\geq 2p+1$ et d'apr\`es le th\'eor\`eme \ref{res stable U} (et le lemme \ref{critere simple}), l'application de restriction stable
Res$_H^G $ est injective en restriction \`a $H^{(p\times 1) , (p\times (q-1)} (Sh^0 G)$.

Mais, d'apr\`es le th\'eor\`eme \ref{res nulle}, l'image Res$_H^G (H^{(p\times 1) , (p\times (q-1)} (Sh^0 G) )$ est contenue dans la composante de K\"unneth
$\prod_{G({\Bbb Q})} H^{1,1} (Sh^0 H_1 ) \otimes \ldots \otimes H^{1,1} (Sh^0 H_p )$ de $\prod_{G({\Bbb Q})} H^{p,p} (Sh^0 H)$. Le th\'eor\`eme de
Lefschetz sur les classes de bidegr\'e $(1,1)$ implique alors que la restriction stable Res$_H^G (\alpha )$ de toute classe $\alpha \in
H^{(p \times 1) , (p\times q-1)} (Sh^0 G)$ appartient au sous-espace engendr\'e par les translat\'es de Hecke de classes $[ S_1 \times \ldots \times S_p ]$,
o\`u chaque $S_i$ est une courbe dans $Sh^0 H_i$ \footnote{Remarquons que $S_i$ n'est pas n\'ecessairement une sous-vari\'et\'e de Shimura.}.

D'apr\`es l'injectivit\'e de l'application de restriction stable, une classe $\alpha \in H^{(p\times 1) , (p\times (q-1))} (Sh^0 G)$ est triviale si et seulement si
$g\alpha \wedge [V] =0$ pour toute sous-vari\'et\'e alg\'ebrique $V$ de $Sh^0 H$ et pour tout $g \in G({\Bbb Q})$.
Or les $G({\Bbb Q})$-modules $H^{2p}_{{\rm n.i.}} (Sh^0 G)$ et $H^{2pq-2p}_{{\rm n.i.}} (Sh^0 G)$ sont en dualit\'e (d\'eduite de la dualit\'e de Poincar\'e).
L'espace engendr\'e par les translat\'es de Hecke dans $Sh^0 G$ des classes $[V]$ de sous-vari\'et\'es alg\'ebriques de $Sh^0 H$ contient donc l'espace
$H^{2pq-2p}_{{\rm n.i.}} (Sh^0 G) \cap H^{pq-p,pq-p} (Sh^0 G)$. Ce qui prouve que toute classe de Hodge dans $H^{2pq-2p} (Sh^0 G)$ est alg\'ebrique.

\bigskip

Concluons cette partie en remarquant que le th\'eor\`eme \ref{res stable U} implique imm\'edia-\\tement, \`a l'aide des th\'eor\`emes de Clozel contenus dans \cite{Clozel}, de
nouveaux r\'esultats d'annulation de la cohomologie de certaines vari\'et\'es de Shimura unitaires. Nous nous contenterons d'\'enoncer un r\'esultat particulier mais suffisamment frappant, renvoyant le lecteur \`a \cite{Clozel} pour d\'eduire du th\'eor\`eme \ref{res stable U} des \'enonc\'es plus g\'en\'eraux.

\begin{thm} \label{cohom 0}
Supposons que $Sh^0 G$ contienne une sous-vari\'et\'e de Shimura $Sh^0 H$, avec $H$ obtenu par restriction des scalaires \`a partir d'un groupe $U(D)$ o\`u $D$ est une
alg\`ebre \`a division de degr\'e premier impair sur une extension quadratique imaginaire d'un corps de nombre totalement r\'eel et tel que $H^{{\rm nc}} = U(p-a,q)$ (resp. $=U(p,q-b)$)
avec $a\geq 1$ (resp. $b\geq 1$). Soit $(\lambda , \mu)$ un couple de partitions compatible de diagramme gauche $\mu / \lambda = (p_1 \times q_1 )* \ldots * (p_m \times q_m )$.
Supposons que
\begin{itemize}
\item $q_1 + \ldots + q_m =q$ (resp. $p_1 + \ldots + p_m =p$),
\item $p_i \geq a$ (resp. $q_i \geq b$) pour chaque $i=1, \ldots , m$, et
\item $|\lambda|+|\hat{\mu}| < pq-aq$ (resp. $<pq-bp$).
\end{itemize}
Alors, $H^{\lambda , \mu} (Sh^0 G) =0$.
\end{thm}

On v\'erifiera facilement que ce dernier r\'esultat implique le th\'eor\`eme \ref{vanish} annonc\'e en introduction.

\bigskip

\section{Cas des vari\'et\'es de Shimura associ\'e au groupe $GSp$}

Dans cette partie $G$ est un groupe alg\'ebrique r\'eductif connexe et anisotrope sur ${\Bbb Q}$ avec $G^{{\rm nc}} = GSp_p$, o\`u $p$ est un entier strictement positif.
Rappelons alors que
\begin{eqnarray} \label{GSp}
G^{{\rm nc}} = \left\{ g = \left(
\begin{array}{cc}
A & B \\
C & D
\end{array} \right) \in U(p,p) \; : \; {}^t g \left(
\begin{array}{cc}
0 & 1_p \\
-1_p & 0
\end{array} \right) g = \left(
\begin{array}{cc}
0 & 1_p \\
-1_p & 0
\end{array} \right) \right\}.
\end{eqnarray}
Soit $K$ le groupe $U(p)$ plong\'e dans $G^{{\rm nc}}$ via l'application
$$k \mapsto \left(
\begin{array}{cc}
k & 0 \\
0 & {}^t k^{-1}
\end{array} \right) .$$
Le complexifi\'e $K_{{\Bbb C}}$ est le groupe $GL_p$. L'involution de Cartan $\theta$ est donn\'ee par $x \mapsto -{}^t \overline{x}$. Soit $T$ le sous-groupe
de $K_{{\Bbb C}}$ constitu\'e des matrices diagonales.

Comme dans la premi\`ere partie, nous noterons $\mathfrak{g}_0$, $\mathfrak{k}_0 ,\ldots$ les alg\`ebres de Lie de $G^{{\rm nc}}$, $K, \ldots$ et $\mathfrak{g}$, $\mathfrak{k}, \ldots$
leur complexifications. Nous noterons $(x_1 , \ldots , x_p ; -x_1 , \ldots \\ , -x_p )$ les \'el\'ements de l'alg\`ebre de Lie de $T$ (vus comme \'el\'ements de $\mathfrak{g}$).

Dans l'alg\`ebre de Lie complexe $\mathfrak{g}$, on a~:
$$\mathfrak{p}^+ = \left\{
\left(
\begin{array}{cc}
0 & B \\
0 & 0
\end{array}
\right) \ {\rm avec } \ B \in M_{p\times p} ({\Bbb C} ) \mbox{ sym\'etrique} \right\}$$
et
$$\mathfrak{p}^- = \left\{
\left(
\begin{array}{cc}
0 & 0 \\
C & 0
\end{array}
\right) \ {\rm avec } \ C \in M_{p\times p} ({\Bbb C} ) \mbox{ sym\'etrique} \right\} .$$
Nous noterons $E={\Bbb C}^p$ la repr\'esentation standard de $U(p)$.
Alors, comme nous l'avons d\'ej\`a remarqu\'e au cours de la premi\`ere partie, comme repr\'esentation de $K_{{\Bbb C}}$, $\mathfrak{p}^+ = {\rm sym}^2 (E )$.
De la m\^eme mani\`ere, $\mathfrak{p}^- = {\rm sym}^2 (E^* )$.

Soit $(e_1 , \ldots , e_p)$ la base canonique de $E$. Choisissons comme sous-alg\`ebre de Borel $\mathfrak{b}_K$ dans $\mathfrak{k}$ l'alg\`ebre des matrices dans
$\mathfrak{k}$ qui sont triangulaires sup\'erieures sur $E$ par rapport \`a cette base. Alors l'ensemble des racines simples compactes positives
\begin{eqnarray} \label{rac cpctes}
\Phi (\mathfrak{b}_K , \mathfrak{t} ) = \left\{ x_i - x_j \; : \; 1\leq i < j \leq p \right\} .
\end{eqnarray}
Les racines simples positives de $T$ apparaissant dans $\mathfrak{p}^+$ sont les formes lin\'eaires $x_i + x_j$ avec $1\leq i \leq j \leq p$.
Dans la suite nous noterons $e_{i,j} = e_i . e_j$ la base canonique de sym$^2 ( E)$. (Vue comme matrice sym\'etrique l'\'el\'ement $e_{i,j}$ est donc \'egal
\`a $E_{i,j} + E_{j,i}$.)

\subsection{D\'ecomposition ``\`a la Lefschetz'' de la cohomologie}

\subsubsection*{Modules cohomologiques et diagrammes de Young}

Nous avons vu dans la premi\`ere partie comment associer une sous-alg\`ebre parabolique $\theta$-stable $\mathfrak{q}$ \`a un \'el\'ement
$X= (x_1 , \ldots , x_p ; -x_1 , \ldots , -x_p )  \in i \mathfrak{t}_0$ (les $x_i$ sont donc tous r\'eels). Rappelons le choix fix\'e (\ref{rac cpctes}) de racines simples
compactes positives. Apr\`es conjugaison par un \'el\'ement de $K$, on peut supposer, et nous le supposerons effectivement par la suite, que $X$ est dominant par
rapport \`a $\Phi (\mathfrak{b}_K , \mathfrak{t} )$, {\it i.e.} que $\alpha (X) \geq 0$ pour tout $\alpha \in \Phi (\mathfrak{b}_K , \mathfrak{t} )$; il satisfait alors aux in\'egalit\'es
$$x_1 \geq \ldots \geq x_p .$$

Nous associons maintenant \`a notre \'el\'ement $X\in i \mathfrak{t}_0$ un couple $(\lambda ,\mu )$ de partitions comme suit.
\begin{itemize}
\item La partition $\lambda \subset p\times p$ est associ\'ee au sous-diagramme de Young de $p\times p$ constitu\'e des cases de coordonn\'ees $(i,j)$ telles que
$x_i + x_j >0$.
\item  La partition $\mu \subset p\times p$ est associ\'ee au sous-diagramme de Young de $p\times p$ constitu\'e des cases de coordonn\'ees $(i,j)$ telles que
$x_i + x_j \geq 0$.
\end{itemize}

Il est imm\'ediat que le couple de partitions $(\lambda , \mu)$ ainsi d\'efini est compatible, que $(\lambda^* , \mu^* ) = (\lambda , \mu )$ et que, r\'eciproquement, tout
couple de partitions compatible dans $p\times p$ et v\'erifiant $(\lambda^* , \mu^* ) = (\lambda , \mu )$ est associ\'e \`a un \'el\'ement $X$ dans $i \mathfrak{t}_0$.

Nous dirons d'une partition $\lambda$ qu'elle est {\it sym\'etrique} si $\lambda^* =\lambda$. Remarquons que si $(\lambda , \mu)$
est un couple compatible de partitions sym\'etriques dans $p\times p$, le diagramme gauche $\mu / \lambda$ est sym\'etrique et s'\'ecrit donc
$(a_1 \times b_1) * \ldots * (a_{m} \times b_{m} )* (p_0 \times p_0 )* (b_{m} \times a_{m} )* \ldots *(b_1 \times a_1)$ pour un certain $m\geq 1$ et des
entiers $a_i$, $b_i$ et $p_0 \geq 1$.

Comme dans la premi\`ere partie, on d\'eduit alors de la remarque suivant la d\'efinition des modules $A_{\mathfrak{q}}$ que chaque couple
compatible de partitions sym\'etriques $(\lambda , \mu )$ d\'efinit sans ambiguit\'e une classe d'\'equivalence de $(\mathfrak{g} , K)$-modules que nous
noterons $A(\lambda , \mu )$. Nous nous autoriserons encore \`a parler de ``la'' sous-alg\`ebre parabolique $\mathfrak{q} (\lambda , \mu ) = \mathfrak{l} (\lambda , \mu )
\oplus  \mathfrak{u} (\lambda , \mu )$ de $(\mathfrak{g} ,K)$-module associ\'e $A(\lambda , \mu )$, l'important pour nous est qu'une telle sous-alg\`ebre existe. Nous supposerons
de plus, ce que l'on peut toujours faire, que le groupe $L(\lambda , \mu )$ associ\'e \`a la sous-alg\`ebre de Levi $\mathfrak{l} (\lambda , \mu )$ n'a pas de facteurs compacts non
ab\'elien. Il est alors facile de voir que
\begin{eqnarray} \label{LL}
L(\lambda , \mu )/ (L(\lambda , \mu ) \cap K)  = GSp_{p_0} /U(p_0 ) \times \prod_{i=1}^m U(a_i, b_i ) / U(a_i ) \times U(b_i ) ,
\end{eqnarray}
o\`u le plongement du groupe $GS_{p_0} \times \prod_{i=1}^m U(a_i , b_i )$ dans $GSp_p$ est (\`a conjugaison dans $GSp_p$
pr\`es) induit par le plongement
$$\begin{array}{lcl}
U(p_0 , p_0 ) \times U(a_1 , b_1) \times \ldots \times U(a_m , b_m ) & \longrightarrow & U(p,p) \\
(g_0 , g_1 , \ldots , g_m ) & \longmapsto & (g_0 , g_1 , \tilde{g}_1 , \ldots , g_m , \tilde{g}_m , id ).
\end{array} $$

Les r\'esultats de Parthasarathy, Kumaresan et Vogan-Zuckerman d\'ecrits dans la premi\`ere partie affirment alors que
$$\hat{G}^{{\rm nc}}_{{\rm VZ}} := \{ A(\lambda , \mu ) \; : \; (\lambda , \mu) \mbox{ compatible et } \lambda , \mu \mbox{ sym\'etriques} \}$$
est l'ensemble des $(\mathfrak{g} ,K)$-modules ayant des groupes de $(\mathfrak{g},K)$-cohomologie non nuls.

\medskip

\'Etant donn\'e une partition sym\'etrique $\lambda \subset p\times p$, nous noterons $\lambda^+$ la partition $( \lambda_1^+ , \ldots , \lambda_p^+)$
o\`u pour $i=1, \ldots ,p$
\begin{eqnarray*}
\lambda_i^+ & = & | \{ j \geq i \; : \; (i,j) \in \lambda \} | , \\
                       & = & \max (0 , \lambda_i -i+1 ).
\end{eqnarray*}
Nous noterons $\overline{\lambda}$ le diagramme de Young $\subset p\times (p+1)$ obtenu en rajoutant une case \`a chaque ligne intersectant la diagonale. Remarquons
alors que $|\overline{\lambda}|$ est n\'ecessairement pair \'egal \`a $2|\lambda^+ |$.
Si $(\lambda , \mu)$ est un couple compatible de partitions sym\'etriques $\subset p\times p$ dont le
diagramme gauche associ\'e $\mu / \lambda = (a_1 \times b_1) * \ldots * (a_{m} \times b_{m} )* (p_0 \times p_0 )* (b_{m} \times a_{m} )* \ldots *(b_1 \times a_1)$, nous
noterons enfin $(\mu/ \lambda )^+$ le sous-diagramme gauche \'egal \`a $(p_0 \times p_0 ) * (b_{m} \times a_{m} )* \ldots *(b_1 \times a_1)$.

Consid\'erons maintenant la repr\'esentation de $K_{{\Bbb C}}$
\begin{eqnarray} \label{U}
V( \lambda ) := E^{\overline{\lambda}} .
\end{eqnarray}
C'est une sous-repr\'esentation irr\'eductible de $\bigwedge^{|\lambda^+|} {\rm sym}^2 (E)$; son vecteur de plus haut poids est
\begin{eqnarray} \label{ulambda}
v(\lambda ):= \bigwedge_{i=1}^p \bigwedge_{j=i}^{\lambda_i} e_{i,j}
\end{eqnarray}
et son vecteur de plus bas poids est
\begin{eqnarray} \label{ylambda}
w(\lambda ) := \bigwedge_{i=1}^p \bigwedge_{j=i}^{\lambda_i} e_{p-i+1 , p-j+1} .
\end{eqnarray}
Puisque d'apr\`es \cite[Lemma 3.5]{HottaWallach} il existe une correspondance bijective entre les sous-espaces irr\'eductibles de $\bigwedge \mathfrak{p}^+$ et les syst\`emes
positifs de racines contenant l'ensemble (\ref{rac cpctes}), la repr\'esentation
\begin{eqnarray} \label{dec en irred}
\bigwedge \mathfrak{p}^+ = \bigwedge {\rm sym}^2 (E) = \bigoplus_{\begin{array}{l}
\lambda \subset p\times p \\
\lambda = \lambda^*
\end{array} } V(\lambda ) ,
\end{eqnarray}
o\`u chaque sous-espace irr\'eductible $V(\lambda )$ apparait avec multiplicit\'e un.

Soit maintenant $(\lambda , \mu )$ un couple compatible de partitions sym\'etriques. Le vecteur
\begin{eqnarray} \label{uy}
v(\lambda ) \otimes w( \hat{\mu} )^* \in \bigwedge^{|\lambda^+|} ({\rm sym}^2 (E) ) \otimes \bigwedge^{|\hat{\mu}^+ |} ({\rm sym}^2 (E))^* = \bigwedge^{|\lambda^+| , |\hat{\mu}^+|}
\mathfrak{p} \subset \bigwedge^{|\lambda^+| + |\hat{\mu}^+|} \mathfrak{p}
\end{eqnarray}
est un vecteur de plus haut poids $2\rho (\mathfrak{u} (\lambda , \mu ) \cap \mathfrak{p} )$ et engendre donc sous l'action de $K_{{\Bbb C}}$ un sous-module irr\'eductible que l'on
note $V(\lambda , \mu )$. Ce module est isomorphe \`a $V(\mathfrak{q}(\lambda , \mu ))$.

\subsubsection*{Classes de Chern et diagrammes de Young}

Consid\'erons $G^{{\rm nc}}_{{\Bbb C}}$ le complexifi\'e du groupe $G^{{\rm nc}} = GSp$. Il s'identifie au sous-groupe de $GL(2p, {\Bbb C})$
constitu\'e des transformations de ${\Bbb C}^{2p}$ qui pr\'eservent la forme symplectique
$\omega = z_1 \wedge z_{p+1} + z_2 \wedge z_{p+2} + \ldots + z_p \wedge z_{2p}$.
Nous noterons $GSp(p)$ le groupe compact obtenu en intersectant $G^{{\rm nc}}_{{\Bbb C}}$ avec $U(2p)$ (c'est un sous-groupe compact maximal).

Soit toujours ${\Bbb G}_{p,p}$ la grassmanienne  des sous-espaces complexes de dimension $p$ dans ${\Bbb C}^{2p}$. Le dual compact ${\Bbb GS}_p := \hat{X}_G$ de
$X_G$ s'identifie au sous-ensemble de ${\Bbb G}_{p,p}$ constitu\'e de tous les espaces totalement isotropes (par rapport \`a la forme symplectique $\omega$)
de ${\Bbb C}^{2p}$. Le groupe $U(2p)$ agit transitivement sur ${\Bbb G}_{p,p}$, le sous-ensemble ${\Bbb GS}_p$ est une sous-vari\'et\'e projective lisse de ${\Bbb G}_{p,p}$
invariante sous l'action du groupe $GSp(p)$ et cette action est transitive. Soit $x_0 \in {\Bbb GS}_p$ le point correspondant au sous-espace de dimension $p$ de ${\Bbb C}^{2p}$
constitu\'e des points dont les $p$ derni\`eres coordonn\'ees sont nulles.

Soit $\nu \subset p\times p$ une partition. On lui a associ\'e dans la premi\`ere partie une classe $C_{\nu } \in H^* ({\Bbb G}_{p,p} )$. Nous noterons dor\'enavant
$C_{\nu}^{{\rm Gr}}$ ces classes de cohomologie. Rappelons que pour tout entier $k\geq1$, la classe de cohomologie $C_{(1^k )}^{{\rm Gr}}$ (resp. $C_{(k)}^{{\rm Gr}}$) est
un multiple non nul de la $k$-i\`eme classe de Chern $\hat{C}_k$ (resp. $\hat{C}_k '$) du fibr\'e $\hat{T}$ (resp. $\hat{Q}$) sur la grassmannienne ${\Bbb G}_{p,p}$.

En utilisant (\ref{dec en irred}) et son dualis\'e~:
$$\bigwedge \mathfrak{p} = \bigoplus_{\begin{array}{l}
\lambda \subset p\times p \\
\lambda = \lambda^*
\end{array} } V(\lambda )^* ,$$
on peut d\'ecrire une base $\{ C_{\nu} \; : \; \nu \subset p\times p , \; \nu^* = \nu \}$ de l'espace $(\bigwedge \mathfrak{p} )^K$ des vecteurs $K$-invariants de $\bigwedge \mathfrak{p}$
param\'etr\'ee par l'ensemble des partitions sym\'etriques $\nu \subset p\times p$. On prend $C_{\nu} := \sum_l z_l \otimes z_l^* $ o\`u $\{ z_l \}$ est une base de
$V(\nu ) \subset \bigwedge \mathfrak{p}^+$ et $\{ z_l^* \}$ la base duale de $V(\nu )^* \subset \bigwedge \mathfrak{p}^-$.

Soit $\nu \subset p\times p$ une partition sym\'etrique.
Le th\'eor\`eme d\'ej\`a mentionn\'e dans la premi\`ere partie identifie $(\bigwedge \mathfrak{p} )^K$ et $H^* ({\Bbb GS}_p )$. On peut donc voir $C_{\nu}$ comme
une classe de cohomologie dans $H^* ({\Bbb GS}_p )$.

Le plongement ${\Bbb GS}_p \rightarrow {\Bbb G}_{p,p}$ induit l'application (de restriction) $H^*({\Bbb G}_{p,p} ) \rightarrow H^* ({\Bbb GS}_p )$ en cohomologie.

\begin{prop} \label{res chern}
Soit $\lambda \subset p\times p$ une partition. L'image de $C_{\lambda}^{{\rm Gr}}$ dans $H^* ({\Bbb GS}_p )$ par l'application de restriction
$H^*({\Bbb G}_{p,p} ) \rightarrow H^* ({\Bbb GS}_p )$ est non nulle si et seulement si $\lambda $ ou $\lambda^* = \nu^+$ pour une certaine partition sym\'etrique
$\nu \subset p \times p$. Et alors, l'image de $C_{\lambda }^{{\rm Gr}}$ dans $H^* ({\Bbb GS}_p )$ est \'egale \`a $C_{\nu }$.
\end{prop}
{\it D\'emonstration.} Dans la suite nous aurons \`a consid\'erer \`a la fois le cas du groupe $U(p,p)$ et du groupe $GSp_p$. Pour distinguer ces deux cas nous utiliserons
les notations \'evidentes $K_{{\rm Gr}}$, $\mathfrak{p}_{{\rm Gr}}$, $\mathfrak{p}_{{\rm Gr}}^+$, $\mathfrak{p}_{{\rm Gr}}^-$ et $V_{{\rm Gr}} (\lambda ) \subset \bigwedge \mathfrak{p}_{{\rm Gr}}^+$ lorsque nous parlerons du cas du groupe $U(p,p)$.

L'inclusion de l'espace des matrices sym\'etriques dans l'espace de toutes les matrices induit une inclusion $\mathfrak{p}^- \rightarrow \mathfrak{p}_{{\rm Gr}}^-$ qui
commute \`a l'action de $K$. Par dualit\'e (pour la forme de Killing) cette inclusion induit \`a son tour l'application $\mathfrak{p}_{{\rm Gr}}^+ \rightarrow \mathfrak{p}^+$ que nous
dirons ``de restriction''. Cette derni\`ere application induit
\begin{eqnarray} \label{gr vers s}
\bigwedge \mathfrak{p}_{{\rm Gr}}^+ \rightarrow \bigwedge \mathfrak{p}^+ .
\end{eqnarray}

Puisque vue comme forme diff\'erentielle sur ${\Bbb G}_{p,p}$ (resp. ${\Bbb GS}_p$), la classe $C_{\lambda}^{{\rm Gr}}$ (resp. $C_{\nu}$) est  invariante sous l'action du groupe
$U(2p)$ (resp. $GSp(p)$), il nous suffit de v\'erifier la proposition \ref{res chern} au point base $x_0 \in {\Bbb GS}_p$. C'est l'objet du lemme suivant.

\begin{lem} \label{rr}
L'image de $V_{{\rm Gr}} (\lambda ) $ dans $\bigwedge \mathfrak{p}^+$ sous l'application (\ref{gr vers s}) est non triviale si et seulement
si $\lambda$ ou $\lambda^* = \nu^+$ pour une certaine partition sym\'etrique $\nu \subset p\times p$. Et alors, son image est $V(\nu )$.
\end{lem}
{\it D\'emonstration.} D'apr\`es (\ref{Vlambda}), $V_{{\rm Gr}} (\lambda ) = E^{\lambda } \otimes (E^{\lambda^*} )^*$.
Vu comme $K$-module, $V_{{\rm Gr}} (\lambda )$ est donc \'egal \`a
\begin{eqnarray} \label{reduc}
E^{\lambda } \otimes E^{\lambda^*} .
\end{eqnarray}

Le $K$-module (\ref{reduc}) se d\'ecompose en irr\'eductible comme suit~:
$$E^{\lambda } \otimes E^{\lambda^* } = \bigoplus_{\mu} (E^{\mu} )^{\oplus c_{\lambda \lambda^*}^{\mu}} .$$

Commen\c{c}ons par remarquer que si $\nu \subset p\times p$, l'image de $V_{{\rm Gr}} (\nu^+ )$ (resp. $V_{{\rm Gr}} ((\nu^+ )^*)$) dans $\bigwedge \mathfrak{p}^+$
sous l'application (\ref{gr vers s}) est non triviale \'egale \`a $V(\nu )$. En effet, il est clair que le vecteur
$$\bigwedge_{i=1}^p \bigwedge_{j=i}^{\nu_i} e_i \otimes e_j^* \ \left(\mbox{resp. }  \bigwedge_{i=1}^p \bigwedge_{j=i}^{\nu_i} e_j \otimes e_i^* \right)$$
appartient \`a $V_{{\rm Gr}} (\nu^+)$ (resp. $V_{{\rm Gr}} ((\nu^+ )^*)$). Mais l'image de ce vecteur dans $\bigwedge \mathfrak{p}^+$ est un multiple non nul de $v(\nu )$.

Il nous reste alors \`a montrer que si $\mu \subset p\times p$ est une partition telle que
\begin{itemize}
\item $c_{\lambda \lambda^*}^{\mu} \neq 0$,
\item $\mu = \overline{\nu}$ pour une certaine partition sym\'etrique $\nu \subset p\times p$, et
\item l'image de $E^{\mu} \subset E^{\lambda} \otimes E^{\lambda^*} $ sous l'application (\ref{gr vers s}) est non triviale,
\end{itemize}
alors $\lambda$ ou $\lambda^* = \nu^+$.

Nous le montrons en raisonnant par l'absurde. Supposons donc la partition $\mu = \overline{\nu}$ avec $\lambda$ et $\lambda^* \neq \nu$. Alors, le $K$-module
$$(E^{\lambda} \otimes E^{\lambda^*} ) \oplus (E^{\nu^+} \otimes E^{(\nu^+ )^*} ) $$
apparait comme $K$-sous-module de $\bigwedge ({\rm sym}^2 (E) \oplus \bigwedge^2 E)$.
Or si un \'el\'ement $x \in E^{\lambda} \otimes E^{\lambda^*}$ et un \'el\'ement $y\in \bigwedge (E \otimes E)$ s'envoient sur un m\^eme vecteur dans
$\bigwedge \mathfrak{p}^+$ (sous l'application (\ref{gr vers s})) on a n\'ecessairement $y\in (E^{\lambda^*} \otimes E^{\lambda} )\oplus ( E^{\lambda} \otimes E^{\lambda^*} )$.
Le module $V(\mu )$ apparait donc dans $\bigwedge {\rm sym}^2 (E)$ avec multiplicit\'e $\geq 2$
(puisque $c_{\nu^+ (\nu^+ )^*}^{\mu} \neq 0$ et $c_{\lambda \lambda^*}^{\mu} \neq 0$) ce qui contredit la d\'ecomposition (\ref{dec en irred}).

\bigskip

Comme corrollaire de la proposition \ref{res chern} on retrouve (avec une d\'emonstration similaire) un r\'esultat de Parthasarathy.
Soit $\hat{T}_s$ le tir\'e en arri\`ere du fibr\'e tautologique sur ${\Bbb G}_{p,p}$ en un fibr\'e sur ${\Bbb GS}_p$. Nous appellerons $\hat{T}_s$
le {\it fibr\'e tautologique} sur ${\Bbb GS}_p$. Soit $c_i (\hat{T}_s ) \in H^* ({\Bbb GS}_p )$ les classes de Chern de $\hat{T}_s$.

\begin{cor} \label{identification de chern}
La $i$-\`eme classe de Chern $c_i (\hat{T}_s )$ du fibr\'e tautologique $\hat{T}_s$ sur ${\Bbb GS}_p$ est un multiple non nul de $C_{(i , 1^{i-1})}$,
o\`u $(i,1^{i-1})$ d\'esigne la partition sym\'etrique de diagramme de Young associ\'e
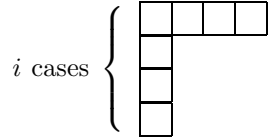
$$i \ {\rm cases}  \; \left\{
\begin{tabular}{l}
\begin{tabular}{|c|c|c|c|} \hline
 & & & \\ \hline
\end{tabular}
\\
\begin{tabular}{|c|} \hline
 \\ \hline
\end{tabular}
\\
 \begin{tabular}{|c|} \hline
  \\ \hline
\end{tabular}
\\
\begin{tabular}{|c|} \hline
 \\ \hline
\end{tabular}
\end{tabular} \right. $$
\end{cor}
{\it D\'emonstration.} La classe $c_i (\hat{T}_s )$ est obtenue par restriction de la $i$-\`eme classe de Chern du fibr\'e tautologique $\hat{T}$ sur ${\Bbb G}_{p,p}$.
Cette derni\`ere est un multiple non nul de $C_{(i)}^{{\rm Gr}}$. Puisque $(i,1^{i-1})^+ = (i)$, le corollaire \ref{identification de chern} d\'ecoule de la proposition
\ref{res chern}.

\bigskip

\subsubsection*{Action des classes de Chern sur la cohomologie}

Dans cette section on \'enonce et d\'emontre l'analogue de la proposition \ref{clef} dans le cas du groupe symplectique.

Nous noterons $\mathfrak{p}_L^+$ l'intersection $\mathfrak{p}^+ \cap \mathfrak{l} ( \lambda , \mu )$ et $E(G,L)$ le
sous-espace de $\bigwedge \mathfrak{p}^+$ engendr\'e par les translat\'es par $K$ du sous-espace $\bigwedge \mathfrak{p}_L^+$.

Enfin, nous dirons qu'une partition sym\'etrique $\nu$ {\it s'inscrit sym\'etriquement} dans un diagramme gauche sym\'etrique
$\mu / \lambda = (a_1 \times b_1 ) * \ldots * (a_m \times b_m ) * (p_0 \times p_0 ) * (b_m \times a_m ) * \ldots * (b_1 \times a_1)$ s'il existe une
image de la partition $\nu^+$ ou de sa transpos\'ee $(\nu^+)^*$ dans le diagramme gauche $(\mu / \lambda )^+ $ dont
le sous-diagramme contenu dans $p_0 \times p_0$ est \'egal \`a $\nu_0^+$ ou $(\nu_0^+)^*$ pour une certaine partition
sym\'etrique $\nu_0 \subset p_0 \times p_0$.

\begin{prop} \label{clef2}
Soient $\lambda$, $\mu$ et $\nu$ trois partitions sym\'etriques incluses dans $p\times p$ telles que $(\lambda , \mu)$ forme un couple compatible.
Notons $L = L( \lambda , \mu )$. Alors, les \'enonc\'es suivants sont \'equivalents~:
\begin{enumerate}
\item $C_{\nu} . V( \lambda , \mu ) \neq 0$ dans $\bigwedge \mathfrak{p}$;
\item la partition $\nu$ s'inscrit sym\'etriquement dans le diagramme gauche $\mu / \lambda $;
\item $V(\nu ) \subset E(G,L )$.
\end{enumerate}
De plus, les \'el\'ements $\{ C_{\nu} . v(\lambda ) \otimes w(\mu )^* \}$, o\`u $\nu$ d\'ecrit l'ensemble des partitions sym\'etriques $\subset p\times p$ qui s'inscrivent sym\'etriquement
dans $\mu / \lambda $, sont lin\'eairement ind\'ependants.
\end{prop}
{\it D\'emonstration.} La d\'emonstration de la proposition \ref{clef2} est similaire \`a celle de la proposition \ref{clef}, les diff\'erents lemmes \ref{lemvenky1}, \ref{lemvenky2} et \ref{lemvenky3}
emprunt\'es \`a Venkataramana sont d'ailleurs compl\`etement g\'en\'eraux. La seule diff\'erence est donc dans l'\'etude de l'application naturelle de restriction de
$\hat{X}_G$ vers $\hat{X}_L$. Pour comprendre celle-ci, commen\c{c}ons par remarquer que l'on a un diagramme commutatif~:
\begin{eqnarray} \label{diag comm}
\begin{array}{ccc}
H^* ({\Bbb G}_{p,p} ) & \longrightarrow & H^* ({\Bbb G}_{p_0 , p_0 }) \otimes \bigotimes_{i=1}^m H^* ({\Bbb G}_{a_i , b_i}) \\
\downarrow & & \downarrow \\
H^* ({\Bbb GS}_p ) & \stackrel{{\rm res}}{\longrightarrow} & H^* ({\Bbb GS}_{p_0} ) \otimes \bigotimes_{i=1}^m H^* ({\Bbb G}_{a_i , b_i} ),
\end{array}
\end{eqnarray}
o\`u l'on a conserv\'e les notations de (\ref{LL}) (le plongement de ${\Bbb G}_{a_i , b_i}$ dans ${\Bbb G}_{p,p}$ n'est donc pas standard).

\begin{lem} \label{res bizarre}
Soit $\nu\subset p\times p$ une partition sym\'etrique. Alors, l'image de $C_{\nu}$ par l'application de restriction en bas du diagramme commutatif (\ref{diag comm})
est une combinaison lin\'eaire \`a coefficients entiers strictement positifs de
$$ C_{\nu_0} \otimes C_{\alpha_1}^{{\rm Gr}} \otimes \ldots \otimes C_{\alpha_m}^{{\rm Gr}} ,$$
o\`u $\nu_0^+ * \alpha_1 * \ldots * \alpha_m$ ou $(\nu_0^+)^* * \alpha_1 * \ldots * \alpha_m$, avec $\nu_0 \subset p_0 \times p_0$ sym\'etrique et chaque
$\alpha_i \subset a_i \times b_i$, est une image de $\nu^+$ ou de sa transpos\'ee $(\nu^+)^*$.
\end{lem}
{\it D\'emonstration du lemme \ref{res bizarre}.}  \`A l'aide du diagramme (\ref{diag comm}), la d\'emonstration du lemme \ref{res bizarre} d\'ecoule imm\'ediatement
de la formule de Whitney sur les classes de Chern, du lemme \ref{res de C} et de la proposition \ref{res chern}.

\bigskip

\subsubsection*{D\'ecomposition ``\`a la Lefschetz''}

Fixons maintenant un sous-groupe de congruence $\Gamma$ dans $G$. Soit $(\lambda , \mu)$ un couple compatible de partitions sym\'etriques. Nous noterons
$H^{\lambda , \mu} (S(\Gamma )) = H^{|\lambda^+| +|\hat{\mu}^+ |} (A(\lambda , \mu ) : \Gamma )$ la $A(\lambda , \mu )$-composante fortement primitive de la cohomologie
de $S(\Gamma )$. La formule de Matsushima et la classification de Vogan-Zuckerman impliquent le th\'eor\`eme suivant.

\begin{thm}
Soit $\Gamma$ un sous-groupe de congruence dans $G$. Pour chaque couple d'entiers $(i,j)$ avec $i+j \leq pq$, on a~:
$$H^{i,j} (S(\Gamma )) = \bigoplus_{
\begin{footnotesize}
\begin{array}{c}
(\lambda , \mu ) \mbox{ couple compatible} \\
\mbox{de partitions sym\'etriques} \\
|\lambda^+| \leq i, \; |\hat{\mu}^+ |\leq j
\end{array}
\end{footnotesize}} \bigoplus_{
\begin{footnotesize}
\begin{array}{c}
\nu_0 \subset p_0 \times p_0 \\
\mbox{sym\'etrique} \\
\nu_i \subset b_i \times a_i \\
i= 1 , \ldots , m
\end{array}
\end{footnotesize}} E^{\lambda , \mu}_{\nu_0 , \nu_1 , \ldots , \nu_m } (S(\Gamma )),$$
o\`u le diagramme gauche $\mu / \lambda = (a_1 \times b_1) * \ldots * (a_{m} \times b_{m} )* (p_0 \times p_0 )* (b_{m} \times a_{m} )* \ldots *(b_1 \times a_1)$ et chaque
$ E^{\lambda , \mu}_{\nu_0 , \nu_1 , \ldots , \nu_m } (S(\Gamma ))$ est isomorphe \`a $H^{\lambda , \mu} (S(\Gamma ))$.
\end{thm}

L\`a encore, cette d\'ecomposition est en g\'en\'eral plus fine que celle induite par les classes de Chern. La d\'emonstration du th\'eor\`eme \ref{dec de Lefschetz} en rempla\c{c}ant
la proposition \ref{clef} par la proposition \ref{clef2} implique imm\'ediatement le th\'eor\`eme suivant.

\medskip

Commen\c{c}ons par remarquer que l'orbite de $x_0 \in {\Bbb GS}_p$ sous l'action du groupe $GS_p$ s'identifie \`a l'espace sym\'etrique hermitien $X_G$. Le groupe $GS_p$ agit sur $(\hat{T}_s)_{|X_G}$. En quotientant par l'action de $\Gamma$
sur $(\hat{T}_s)_{|X_G}$, on obtient un fibr\'e sur $S(\Gamma )$ que nous noterons $T_s$.

D'apr\`es un th\'eor\`eme classique de Cartan, l'espace $H^* ({\Bbb GS}_p )$ peut \^etre identifi\'e avec l'espace des formes diff\'erentielles $GSp(p)$-invariantes.
Soit $\omega$ une forme diff\'erentielle $GSp(p)$-invariante sur ${\Bbb GS}_{p}$ et soit
$\bar{\omega}$ une forme diff\'erentielle $GSp(p)$-invariante sur $X_G$ \'egale \`a $\omega$ au point $x_0$. Puisque
$\omega$ est, en particulier, $\Gamma$-invariante elle induit une forme (n\'ecessairement ferm\'ee) sur $S(\Gamma )$ qui d\'efinit
donc une classe de cohomologie dans $H^* (S(\Gamma ))$. On a ainsi construit une application
\begin{eqnarray} \label{eta2}
\eta : H^* ({\Bbb GS}_{p} ) \rightarrow H^* (S(\Gamma )).
\end{eqnarray}
Il est bien connu que $\eta$ est injective. Enfin, la $i$-\`eme classe de Chern $c_i (T_s )$ du fibr\'e $T_s$ est \'egale \`a $(-1)^i \eta (c_i (\hat{T}_s ))$.

\begin{thm} \label{dec de lef2}
Soit $\Gamma$ un sous-groupe de congruence dans $G$ et soit $\eta : H^* ({\Bbb GS}_p ) \rightarrow H^* (S(\Gamma ))$ l'application d\'efinie en (\ref{eta2}).
Fixons $\lambda$, $\mu$ et $\nu$ trois partitions sym\'etriques incluses dans $p\times p$ telles que le couple $(\lambda, \mu )$ soit compatible. Alors,
\begin{enumerate}
\item pour toute classe fortement primitive $s \in H^{\lambda , \mu} (S(\Gamma ))$, $\eta (C_{\nu} ) .s =0$ si et seulement si la partition $\nu$
ne s'inscrit pas sym\'etriquement dans le diagramme gauche $\mu / \lambda$, et
\item si $s \in H^{\lambda , \mu} (S(\Gamma ))$ est une classe non nulle, les \'el\'ements
$$\{ C_{\nu} .s \; : \; \nu \subset p\times p , \; \nu = \nu^* , \; \nu \mbox{ s'inscrit sym\'etriquement dans } \mu / \lambda \} $$
sont lin\'eairement ind\'ependants.
\end{enumerate}
\end{thm}

En sp\'ecialisant ce th\'eor\`eme au cas de la cohomologie holomorphe, on retrouve un th\'eor\`eme de Parthasarathy \cite[Theorem 4.1]{Parthasarathy}.

\subsection{Restriction stable \`a une sous-vari\'et\'e de Shimura}

Comme dans la premi\`ere partie, nous pourrions classifier les diff\'erents types possibles de sous-vari\'et\'es de Shimura. Nos m\'ethodes n'apporteraient pas
de r\'esultats int\'eressants en dehors du cas des sous-vari\'et\'es $Sh^0 H \subset Sh^0 G$ avec
$$H^{{\rm nc}} = GSp_{p_1} \times \ldots \times GSp_{p_m},$$
o\`u $p_1 + \ldots + p_m \leq p$ et dont le plongement dans $GSp_p$ est induit par
le plongement canonique (\ref{UidansU}) de $U(p_1 , p_1 ) \times \ldots \times U(p_m , p_m )$ dans $U(p,p)$.

Commen\c{c}ons par remarquer que la d\'emonstration de la proposition \ref{dualU} et le lemme \ref{res bizarre} impliquent la proposition suivante.

\begin{prop} \label{dualSS2}
Soient $p_j$ avec $j=1, \ldots , m$ des entiers strictement positifs tels que $p_1 + \ldots + p_m \leq p$. Supposons $H^{{\rm nc}} = GSp_{p_1} \times \ldots \times GSp_{p_m}$
et plong\'e dans $GSp_p$ comme au-dessus. Alors, \`a un multiple scalaire non nul pr\`es, la classe duale
$[\hat{X}_H ] \in  H^{\frac{p(p+1)}{2} - \sum_i \frac{p_i (p_i +1)}{2}} ({\Bbb GS}_p )$ est une combinaison lin\'eaire \`a coefficients entiers strictement positifs des classes
$C_{\hat{\nu}}$, o\`u $\nu$ d\'ecrit l'ensemble des partitions sym\'etriques $\subset p\times p$ telles que $\nu^+$ soit l'image d'un diagramme gauche
$\alpha_1 * \ldots * \alpha_m$ o\`u chaque $\alpha_i$ est \'egal au diagramme $(p_i \times q_i )^+$ ou \`a son transpos\'e.
\end{prop}

L'analogue du th\'eor\`eme \ref{res stable U} n'est non vide que dans le cas de la cohomologie holomorphe. Nos m\'ethodes ne permettent donc pas de
d\'emontrer de nouveaux r\'esultats, montrons n\'eanmoins comment retrouver les r\'esultats de Clozel et Venkataramana.

Le sous-espace $H^{\lambda , \mu} (Sh^0 G)$ apparait dans la cohomologie holomorphe si et seulement si $\mu = p\times p$. La partition $\lambda$
est alors naturellement param\'etr\'ee par un entier $r$ compris entre $0$ et  $p$ tel que
$$\lambda = ( \underbrace{p, \ldots , p}_{r \ {\rm fois}} , \underbrace{r , \ldots , r}_{p-r \ {\rm fois}} )$$
de diagramme de Young~:
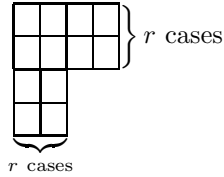
$$
\begin{array}{l}
\left. \hspace{0,035cm}
\begin{array}{|c|c|c|c|} \hline
 & & & \\ \hline
 & & & \\ \hline
\end{array}
\right\}  r \ {\rm cases} \\
\underbrace{
\begin{array}{|c|c|} \hline
 &  \\ \hline
 &  \\ \hline
\end{array}}_{r \ {\rm cases}}
\end{array}
$$
(Ici $p=4$.)

Dans ce cas $H^{\lambda , \mu} (Sh^0 G) = H^{rp- \frac{r(r-1)}{2}, 0} (Sh^0 G)$ et les espaces de cohomologie holomorphe sont triviaux dans tous les autres degr\'es.

La d\'emonstration du corollaire \ref{CV} se traduit facilement (\`a l'aide des propositions \ref{clef2} et \ref{dualSS2}) pour obtenir le corollaire suivant.

\begin{cor}[Clozel-Venkataramana]
Soit $Sh^0 H$ une sous-vari\'et\'e de Shimura de $Sh^0 G$ avec $H^{{\rm nc}} = GSp_{p_1} \times \ldots \times GSp_{p_m}$, $p_j \geq 1$ et $p_1 + \ldots + p_m \leq p$.
Soit $r$ un entiers naturel $\leq p$. Alors, l'application
$${\rm Res}_H^G : H^* (Sh^0 G) \rightarrow \prod_{G({\Bbb Q} )} H^* (Sh^0 H )$$
de restriction stable est {\bf injective} en restriction \`a $H^{rp - \frac{r(r-1)}{2}, 0} (Sh^0 G)$ {\bf si et seulement si} $p_1 + \ldots +p_m =p$
et $r=1$.
\end{cor}

\section{Cas des vari\'et\'es de Shimura associ\'e au groupe $O^* (2p)$}

Dans cette partie $G$ est un groupe alg\'ebrique r\'eductif connexe et anisotrope sur ${\Bbb Q}$ avec $G^{{\rm nc}} = O^* (2p)$, o\`u $p$ est un entier strictement positif.
Rappelons alors que
\begin{eqnarray} \label{Op}
G^{{\rm nc}} = \left\{ g = \left(
\begin{array}{cc}
A & B \\
C & D
\end{array} \right) \in U(p,p) \; : \; {}^t g \left(
\begin{array}{cc}
0 & 1_p \\
1_p & 0
\end{array} \right) g = \left(
\begin{array}{cc}
0 & 1_p \\
1_p & 0
\end{array} \right) \right\}.
\end{eqnarray}
Soit $K$ le groupe $U(p)$ plong\'e dans $G^{{\rm nc}}$ via l'application
$$k \mapsto \left(
\begin{array}{cc}
k & 0 \\
0 & {}^t k^{-1}
\end{array} \right) .$$
Le complexifi\'e $K_{{\Bbb C}}$ est le groupe $GL_p$. L'involution de Cartan $\theta$ est donn\'ee par $x \mapsto -{}^t \overline{x}$. Soit $T$ le sous-groupe
de $K_{{\Bbb C}}$ constitu\'e des matrices diagonales.

Comme dans la premi\`ere partie, nous noterons $\mathfrak{g}_0$, $\mathfrak{k}_0 ,\ldots$ les alg\`ebres de Lie de $G^{{\rm nc}}$, $K, \ldots$ et $\mathfrak{g}$, $\mathfrak{k}, \ldots$
leur complexifications. Nous noterons $(x_1 , \ldots , x_p ; -x_1 , \ldots \\, -x_p )$ les \'el\'ements de l'alg\`ebre de Lie de $T$ (vus comme \'el\'ements de $\mathfrak{g}$).

Dans l'alg\`ebre de Lie complexe $\mathfrak{g}$, on a~:
$$\mathfrak{p}^+ = \left\{
\left(
\begin{array}{cc}
0 & B \\
0 & 0
\end{array}
\right) \ {\rm avec } \ B \in M_{p\times p} ({\Bbb C} ) \mbox{ antisym\'etrique} \right\}$$
et
$$\mathfrak{p}^- = \left\{
\left(
\begin{array}{cc}
0 & 0 \\
C & 0
\end{array}
\right) \ {\rm avec } \ C \in M_{p\times p} ({\Bbb C} ) \mbox{ antisym\'etrique} \right\} .$$
Nous noterons $E={\Bbb C}^p$ la repr\'esentation standard de $U(p)$.
Alors, comme nous l'avons d\'ej\`a remarqu\'e au cours de la premi\`ere partie, comme repr\'esentation de $K_{{\Bbb C}}$, $\mathfrak{p}^+ = \bigwedge^2 (E )$.
De la m\^eme mani\`ere, $\mathfrak{p}^- = \bigwedge^2 (E^* )$.

Soit $(e_1 , \ldots , e_p)$ la base canonique de $E$. Choisissons comme sous-alg\`ebre de Borel $\mathfrak{b}_K$ dans $\mathfrak{k}$ l'alg\`ebre des matrices dans
$\mathfrak{k}$ qui sont triangulaires sup\'erieures sur $E$ par rapport \`a cette base. Alors l'ensemble des racines simples compactes positives
\begin{eqnarray} \label{rac cpctes2}
\Phi (\mathfrak{b}_K , \mathfrak{t} ) = \left\{ x_i - x_j \; : \; 1\leq i < j \leq p \right\} .
\end{eqnarray}
Les racines simples positives de $T$ apparaissant dans $\mathfrak{p}^+$ sont les formes lin\'eaires $x_i + x_j$ avec $1\leq i < j \leq p$.
Dans la suite nous noterons $e_{i,j} = e_i \wedge e_j$ la base canonique de $\bigwedge^2 ( E)$. (Vue comme matrice antisym\'etrique l'\'el\'ement $e_{i,j}=e_i \wedge e_j$ est donc \'egal \`a $E_{i,j} - E_{j,i}$.)

\subsection{D\'ecomposition ``\`a la Lefschetz'' de la cohomologie}

\subsubsection*{Modules cohomologiques et diagrammes de Young}

Nous avons vu dans la premi\`ere partie comment associer une sous-alg\`ebre parabolique $\theta$-stable $\mathfrak{q}$ \`a un \'el\'ement
$X= (x_1 , \ldots , x_p ; -x_1 , \ldots , -x_p )  \in i \mathfrak{t}_0$ (les $x_i$ sont donc tous r\'eels). Rappelons le choix fix\'e (\ref{rac cpctes2}) de racines simples
compactes positives. Apr\`es conjugaison par un \'el\'ement de $K$, on peut supposer, et nous le supposerons effectivement par la suite, que $X$ est dominant par
rapport \`a $\Phi (\mathfrak{b}_K , \mathfrak{t} )$, {\it i.e.} que $\alpha (X) \geq 0$ pour tout $\alpha \in \Phi (\mathfrak{b}_K , \mathfrak{t} )$; il satisfait alors aux in\'egalit\'es
$$x_1 \geq \ldots \geq x_p .$$

Nous associons maintenant \`a notre \'el\'ement $X\in i \mathfrak{t}_0$ un couple $(\lambda ,\mu )$ de partitions comme suit.
\begin{itemize}
\item La partition $\lambda \subset p\times p$ est associ\'ee au sous-diagramme de Young de $p\times p$ constitu\'e des cases de coordonn\'ees $(i,j)$ telles que
$x_i + x_j >0$.
\item  La partition $\mu \subset p\times p$ est associ\'ee au sous-diagramme de Young de $p\times p$ constitu\'e des cases de coordonn\'ees $(i,j)$ telles que
$x_i + x_j \geq 0$.
\end{itemize}

Il est imm\'ediat que le couple $(\lambda , \mu)$ est un couple compatible de partitions sym\'etriques et que, r\'eciproquement, tout
couple compatible de partitions sym\'etriques dans $p\times p$ est associ\'e \`a un \'el\'ement $X$ dans $i \mathfrak{t}_0$.

\'Etant donn\'ee une partition sym\'etrique $\lambda \subset p\times p$, nous noterons $\lambda^-$ la partition $( \lambda_1^- , \ldots , \lambda_p^-)$
o\`u pour $i=1, \ldots ,p$
\begin{eqnarray*}
\lambda_i^- & = & | \{ j > i \; : \; (i,j) \in \lambda \} | , \\
                       & = & \max (0 , \lambda_i -i ).
\end{eqnarray*}
Nous noterons $\check{\lambda}$ le diagramme de Young $\subset p\times p$ obtenu en soustrayant une case \`a chaque ligne intersectant la diagonale. Remarquons
alors que $|\check{\lambda}|$ est n\'ecessairement pair \'egal \`a $2|\lambda^- |$.

Comme dans la premi\`ere partie, on d\'eduit alors de la remarque suivant la d\'efinition des modules $A_{\mathfrak{q}}$ que chaque couple
de partitions $(\lambda^- , \mu^- )$ associ\'e \`a un couple compatible $(\lambda , \mu )$ de partitions sym\'etriques d\'efinit sans ambiguit\'e une classe
d'\'equivalence de $(\mathfrak{g} , K)$-modules que nous
noterons $A(\lambda^- , \mu^- )$. Nous nous autoriserons encore \`a parler de ``la'' sous-alg\`ebre parabolique $\mathfrak{q} (\lambda^- , \mu^- ) = \mathfrak{l} (\lambda^- , \mu^- )
\oplus  \mathfrak{u} (\lambda^- , \mu^- )$ de $(\mathfrak{g} ,K)$-module associ\'e $A(\lambda^- , \mu^- )$, l'important pour nous est qu'une telle sous-alg\`ebre existe. Nous supposerons
de plus, ce que l'on peut toujours faire, que le groupe $L(\lambda^- , \mu^- )$ associ\'e \`a la sous-alg\`ebre de Levi $\mathfrak{l} (\lambda^- , \mu^- )$ n'a pas de facteurs compacts non
ab\'elien. Il est alors facile de voir que
\begin{eqnarray} \label{LL2}
L(\lambda , \mu )/ (L(\lambda , \mu ) \cap K)  = O^* (2p_0 ) /U(p_0 ) \times \prod_{i=1}^m U(a_i, b_i ) / U(a_i ) \times U(b_i ) ,
\end{eqnarray}
o\`u le plongement du groupe $O^* (2p_0 ) \times \prod_{i=1}^m U(a_i , b_i )$ dans $O^* (2p)$ est (\`a conjugaison dans $O^* (2p)$ pr\`es) induit par le plongement
$$\begin{array}{lcl}
U(p_0 , p_0 ) \times U(a_1 , b_1) \times \ldots \times U(a_m , b_m ) & \longrightarrow & U(p,p) \\
(g_0 , g_1 , \ldots , g_m ) & \longmapsto & (g_0 , g_1 , \tilde{g}_1 , \ldots , g_m , \tilde{g}_m , id ).
\end{array} $$

Les r\'esultats de Parthasarathy, Kumaresan et Vogan-Zuckerman d\'ecrits dans la premi\`ere partie affirment alors que
$$\hat{G}^{{\rm nc}}_{{\rm VZ}} := \{ A(\lambda^- , \mu^- ) \; : \; (\lambda , \mu) \mbox{ compatible et } \lambda , \mu \mbox{ sym\'etriques} \}$$
est l'ensemble des $(\mathfrak{g} ,K)$-modules ayant des groupes de $(\mathfrak{g},K)$-cohomologie non nuls.

\medskip

\'Etant donn\'ee une partition sym\'etrique $\lambda \subset p\times p$, consid\'erons maintenant la repr\'esentation de $K_{{\Bbb C}}$
\begin{eqnarray} \label{U2}
V( \lambda^- ) := E^{\check{\lambda}} .
\end{eqnarray}
C'est une sous-repr\'esentation irr\'eductible de $\bigwedge^{|\lambda^-|} \left( \bigwedge^2 (E)\right) $; son vecteur de plus haut poids est
\begin{eqnarray} \label{ulambda2}
v(\lambda^- ) := \bigwedge_{i=1}^p \bigwedge_{j=i+1}^{\lambda_i} e_{i,j}
\end{eqnarray}
et son vecteur de plus bas poids est
\begin{eqnarray} \label{ylambda2}
w(\lambda^- ) := \bigwedge_{i=1}^p \bigwedge_{j=i+1}^{\lambda_i} e_{p-i+1 , p-j+1} .
\end{eqnarray}
Puisque d'apr\`es \cite[Lemma 3.5]{HottaWallach} il existe une correspondance bijective entre les sous-espaces irr\'eductibles de $\bigwedge \mathfrak{p}^+$ et les syst\`emes
positifs de racines contenant l'ensemble (\ref{rac cpctes2}), la repr\'esentation
\begin{eqnarray} \label{dec en irred2}
\bigwedge \mathfrak{p}^+ = \bigwedge \left( \bigwedge^2 (E) \right) = \bigoplus_{\begin{array}{c}
\lambda^- \\
\lambda \subset p\times p \\
\lambda = \lambda^*
\end{array} } V(\lambda^- ) ,
\end{eqnarray}
o\`u chaque sous-espace irr\'eductible $V(\lambda^- )$ apparait avec multiplicit\'e un.

Soit maintenant $(\lambda , \mu )$ un couple compatible de partitions sym\'etriques. Le vecteur
\begin{eqnarray} \label{uy2}
v(\lambda^- ) \otimes w( \hat{\mu}^- )^* \in \bigwedge^{|\lambda^-|} \left( \bigwedge^2 (E) \right) \otimes \bigwedge^{|\hat{\mu}^- |} \left( \bigwedge^2 (E) \right)^* =
\bigwedge^{|\lambda^-| , |\hat{\mu}^-|} \mathfrak{p} \subset \bigwedge^{|\lambda^-| + |\hat{\mu}^-|} \mathfrak{p}
\end{eqnarray}
est un vecteur de plus haut poids $2\rho (\mathfrak{u} (\lambda^- , \mu^- ) \cap \mathfrak{p} )$ et engendre donc sous l'action de $K_{{\Bbb C}}$ un sous-module irr\'eductible que l'on
note $V(\lambda^- , \mu^- )$. Ce module est isomorphe \`a $V(\mathfrak{q}(\lambda^- , \mu^- ))$.

\subsubsection*{Classes de Chern et diagrammes de Young}

Consid\'erons $G^{{\rm nc}}_{{\Bbb C}}$ le complexifi\'e du groupe $G^{{\rm nc}} = O^* (2p)$. Il s'identifie au sous-groupe de $GL(2p, {\Bbb C})$
constitu\'e des transformations de ${\Bbb C}^{2p}$ qui pr\'eservent la forme quadratique
$\omega = z_1  z_{p+1} + z_2 z_{p+2} + \ldots + z_p  z_{2p}$.
Nous noterons $O (2p)$ le groupe compact obtenu en intersectant $G^{{\rm nc}}_{{\Bbb C}}$ avec $U(2p)$ (c'est un sous-groupe compact maximal).

Soit toujours ${\Bbb G}_{p,p}$ la grassmannienne  des sous-espaces complexes de dimension $p$ dans ${\Bbb C}^{2p}$. Le dual compact ${\Bbb GO}_p := \hat{X}_G$ de
$X_G$ s'identifie au sous-ensemble de ${\Bbb G}_{p,p}$ constitu\'e de tous les espaces totalement isotropes (par rapport \`a la forme quadratique $\omega$)
de ${\Bbb C}^{2p}$. Le groupe $U(2p)$ agit transitivement sur ${\Bbb G}_{p,p}$, le sous-ensemble ${\Bbb GO}_p$ est une sous-vari\'et\'e projective lisse de ${\Bbb G}_{p,p}$
invariante sous l'action du groupe $O(2p)$ et cette action est transitive. Soit $x_0 \in {\Bbb GO}_p$ le point correspondant au sous-espace de dimension $p$ de ${\Bbb C}^{2p}$
constitu\'e des points dont les $p$ derni\`eres coordonn\'ees sont nulles.

Soit $\nu \subset p\times p$ une partition. On lui a associ\'e dans la premi\`ere partie une classe $C_{\nu } \in H^* ({\Bbb G}_{p,p} )$. Nous noterons dor\'enavant
$C_{\nu}^{{\rm Gr}}$ ces classes de cohomologie. Rappelons que pour tout entier $k\geq1$, la classe de cohomologie $C_{(1^k )}^{{\rm Gr}}$ (resp. $C_{(k)}^{{\rm Gr}}$) est
un multiple non nul de la $k$-i\`eme classe de Chern $\hat{C}_k$ (resp. $\hat{C}_k '$) du fibr\'e $\hat{T}$ (resp. $\hat{Q}$) sur la grassmannienne ${\Bbb G}_{p,p}$.

En utilisant (\ref{dec en irred2}) et son dualis\'e~:
$$\bigwedge \mathfrak{p} = \bigoplus_{\begin{array}{l}
\lambda^- \\
\lambda \subset p\times p \\
\lambda = \lambda^*
\end{array} } V(\lambda^- )^* ,$$
on peut d\'ecrire une base $\{ C_{\nu^-} \; : \; \nu \subset p\times p , \; \nu^* = \nu \}$ de l'espace $(\bigwedge \mathfrak{p} )^K$ des vecteurs $K$-invariants de $\bigwedge \mathfrak{p}$
param\'etr\'ee par l'ensemble des partitions $\nu^-$ o\`u $\nu$ est une partition sym\'etrique $\subset p\times p$. On prend
$C_{\nu^-} := \sum_l z_l \otimes z_l^* $ o\`u $\{ z_l \}$ est une base de
$V(\nu^- ) \subset \bigwedge \mathfrak{p}^+$ et $\{ z_l^* \}$ la base duale de $V(\nu^- )^* \subset \bigwedge \mathfrak{p}^-$.

Soit $\nu \subset p\times p$ une partition sym\'etrique.
Le th\'eor\`eme d\'ej\`a mentinonn\'e dans la premi\`ere partie identifie $(\bigwedge \mathfrak{p} )^K$ et $H^* ({\Bbb GO}_p )$. On peut donc voir $C_{\nu^-}$ comme
une classe de cohomologie dans $H^* ({\Bbb GO}_p )$.

Le plongement ${\Bbb GO}_p \rightarrow {\Bbb G}_{p,p}$ induit l'application (de restriction) $H^*({\Bbb G}_{p,p} ) \rightarrow H^* ({\Bbb GO}_p )$ en cohomologie.

\begin{prop} \label{res chern2}
Soit $\lambda \subset p\times p$ une partition. L'image de $C_{\lambda}^{{\rm Gr}}$ dans $H^* ({\Bbb GO}_p )$ par l'application de restriction
$H^*({\Bbb G}_{p,p} ) \rightarrow H^* ({\Bbb GO}_p )$ est non nulle si et seulement si $\lambda $ ou $\lambda^* = \nu^-$ pour une certaine partition sym\'etrique
$\nu \subset p \times p$. Et alors, l'image de $C_{\lambda }^{{\rm Gr}}$ dans $H^* ({\Bbb GS}_p )$ est \'egale \`a $C_{\nu^- }$.
\end{prop}
{\it D\'emonstration.} Dans la suite nous aurons \`a consid\'erer \`a la fois le cas du groupe $U(p,p)$ et du groupe $O^* ( 2p)$. Pour distinguer ces deux cas nous utiliserons
les notations $K_{{\rm Gr}}$, $\mathfrak{p}_{{\rm Gr}}$, $\mathfrak{p}_{{\rm Gr}}^+$, $\mathfrak{p}_{{\rm Gr}}^-$ et $V_{{\rm Gr}} (\lambda )$ comme dans la deuxi\`eme partie.

L'inclusion de l'espace des matrices antisym\'etriques dans l'espace de toutes les matrices induit une inclusion $\mathfrak{p}^- \rightarrow \mathfrak{p}_{{\rm Gr}}^-$ qui
commute \`a l'action de $K$. Par dualit\'e (pour la forme de Killing) cette inclusion induit \`a son tour l'application $\mathfrak{p}_{{\rm Gr}}^+ \rightarrow \mathfrak{p}^+$ que nous
dirons ``de restriction''. Cette derni\`ere application induit
\begin{eqnarray} \label{gr vers o}
\bigwedge \mathfrak{p}_{{\rm Gr}}^+ \rightarrow \bigwedge \mathfrak{p}^+ .
\end{eqnarray}

Comme dans la d\'emonstration de la proposition \ref{res chern}, il nous suffit de v\'erifier la proposition \ref{res chern2} au point base $x_0 \in {\Bbb GO}_p$. C'est l'objet du lemme suivant dont la d\'emonstration est similaire \`a celle du lemme \ref{rr}.

\begin{lem}
L'image de $V_{{\rm Gr}} (\lambda ) $ dans $\bigwedge \mathfrak{p}^+$ sous l'application (\ref{gr vers o}) est non triviale si et seulement
si $\lambda$ ou $\lambda^* = \nu^-$ pour une certaine partition sym\'etrique $\nu \subset p\times p$. Et alors, son image est $V(\nu^- )$.
\end{lem}

\bigskip

Soit $\hat{T}_o$ le tir\'e en arri\`ere du fibr\'e tautologique sur ${\Bbb G}_{p,p}$ en un fibr\'e sur ${\Bbb GO}_p$. Nous appellerons $\hat{T}_o$
le {\it fibr\'e tautologique} sur ${\Bbb GO}_p$. Soit $c_i (\hat{T}_o ) \in H^* ({\Bbb GO}_p )$ les classes de Chern de $\hat{T}_o$.

\begin{cor} \label{identification de chern2}
La $i$-\`eme classe de Chern $c_i (\hat{T}_o )$ du fibr\'e tautologique $\hat{T}_o$ sur ${\Bbb GO}_p$ est un multiple non nul de $C_{(i)}$ (ici $1\leq i \leq p-1$).
\end{cor}
{\it D\'emonstration.} La classe $c_i (\hat{T}_o )$ est obtenue par restriction de la $i$-\`eme classe de Chern du fibr\'e tautologique $\hat{T}$ sur ${\Bbb G}_{p,p}$.
Cette derni\`ere est un multiple non nul de $C_{(i)}^{{\rm Gr}}$. Puisque $(i+1,1^{i})^- = (i)$, le corollaire \ref{identification de chern2} d\'ecoule de la proposition
\ref{res chern2}.

\bigskip

\subsubsection*{Action des classes de Chern sur la cohomologie}

Comme dans les deux parties pr\'ec\'edentes, nous noterons $\mathfrak{p}_L^+$ l'intersection $\mathfrak{p}^+ \cap \mathfrak{l} ( \lambda^- , \mu^- )$ et $E(G,L)$ le
sous-espace de $\bigwedge \mathfrak{p}^+$ engendr\'e par les translat\'es par $K$ du sous-espace $\bigwedge \mathfrak{p}_L^+$.

Enfin, nous dirons qu'une partition sym\'etrique $\nu$ {\it s'inscrit antisym\'etriquement} dans un diagramme gauche sym\'etrique
$\mu / \lambda = (a_1 \times b_1 ) * \ldots * (a_m \times b_m ) * (p_0 \times p_0 ) * (b_m \times a_m ) * \ldots * (b_1 \times a_1)$ s'il existe une
image de la partition $\nu^-$ ou de sa transpos\'ee $(\nu^-)^*$ dans le diagramme gauche $(\mu / \lambda )^+ $ dont
le sous-diagramme contenu dans $p_0 \times p_0$ est \'egale \`a $\nu_0^-$ ou $(\nu_0^-)^*$ pour une certaine partition
sym\'etrique $\nu_0 \subset p_0 \times p_0$.

La d\'emonstration de la proposition \ref{clef2} se traduit facilement dans le cas du groupe $O^* (2p)$ pour d\'emontrer la proposition suivante.

\begin{prop} \label{clef3}
Soient $\lambda$, $\mu$ et $\nu$ trois partitions sym\'etriques incluses dans $p\times p$ telles que $(\lambda , \mu)$ forme un couple compatible.
Notons $L = L( \lambda^- , \mu^- )$. Alors, les \'enonc\'es suivants sont \'equivalents~:
\begin{enumerate}
\item $C_{\nu^-} . V( \lambda^- , \mu^- ) \neq 0$ dans $\bigwedge \mathfrak{p}$;
\item la partition $\nu$ s'inscrit antisym\'etriquement dans le diagramme gauche $\mu / \lambda $;
\item $V(\nu^- ) \subset E(G,L )$.
\end{enumerate}
De plus, les \'el\'ements $\{ C_{\nu^- } . v(\lambda^- ) \otimes w(\mu^- )^* \}$, o\`u $\nu^-$ d\'ecrit l'ensemble des parties n\'egatives des partitions sym\'etriques $\nu \subset p\times p$ qui s'inscrivent antisym\'etriquement dans $\mu / \lambda $, sont lin\'eairement ind\'ependants.
\end{prop}

\bigskip

\subsubsection*{D\'ecomposition ``\`a la Lefschetz''}

Fixons maintenant un sous-groupe de congruence $\Gamma$ dans $G$. Soit $(\lambda , \mu)$ un couple compatible de partitions sym\'etriques. Nous noterons
$H^{\lambda^- , \mu^-} (S(\Gamma )) = H^{|\lambda^-| +|\hat{\mu}^- |} (A(\lambda^- , \mu^- ) : \Gamma )$ la $A(\lambda^- , \mu^- )$-composante fortement primitive de la cohomologie
de $S(\Gamma )$. La formule de Matsushima et la classification de Vogan-Zuckerman impliquent le th\'eor\`eme suivant.

\begin{thm}
Soit $\Gamma$ un sous-groupe de congruence dans $G$. Pour chaque couple d'entiers $(i,j)$ avec $i+j \leq \frac{p(p-1)}{2}$, on a~:
$$H^{i,j} (S(\Gamma )) = \bigoplus_{
\begin{footnotesize}
\begin{array}{c}
(\lambda^- , \mu^- ) \\
(\lambda , \mu ) \mbox{ couple compatible} \\
\mbox{de partitions sym\'etriques} \\
|\lambda^-| \leq i, \; |\hat{\mu}^- |\leq j
\end{array}
\end{footnotesize}} \bigoplus_{
\begin{footnotesize}
\begin{array}{c}
\nu_0^- \\
\nu_0 \subset p_0 \times p_0 \\
\mbox{sym\'etrique} \\
\nu_i \subset b_i \times a_i \\
i= 1 , \ldots , m
\end{array}
\end{footnotesize}} E^{\lambda^- , \mu^-}_{\nu_0^- , \nu_1 , \ldots , \nu_m } (S(\Gamma )),$$
o\`u le diagramme gauche $\mu / \lambda = (a_1 \times b_1) * \ldots * (a_{m} \times b_{m} )* (p_0 \times p_0 )* (b_{m} \times a_{m} )* \ldots *(b_1 \times a_1)$ et chaque
$ E^{\lambda^- , \mu^- }_{\nu_0^- , \nu_1 , \ldots , \nu_m } (S(\Gamma ))$ est isomorphe \`a $H^{\lambda^- , \mu^- } (S(\Gamma ))$.
\end{thm}

L\`a encore, cette d\'ecomposition est en g\'en\'eral plus fine que celle induite par les classes de Chern. La d\'emonstration du th\'eor\`eme \ref{dec de Lefschetz} en rempla\c{c}ant
la proposition \ref{clef} par la proposition \ref{clef3} implique imm\'ediatement le th\'eor\`eme suivant.

\medskip

Commen\c{c}ons par remarquer que l'orbite de $x_0 \in {\Bbb GO}_p$ sous l'action du groupe $O^* (2p)$ s'identifie \`a l'espace sym\'etrique hermitien $X_G$. Soit
$\Gamma$ un sous-groupe de congruence de $G$. Le groupe $O^* (2p)$ agit sur $(\hat{T}_o)_{|X_G}$. En quotientant par l'action de $\Gamma$
sur $(\hat{T}_o)_{|X_G}$, on obtient un fibr\'e sur $S(\Gamma )$ que nous noterons $T_o$.

D'apr\`es un th\'eor\`eme classique de Cartan, l'espace $H^* ({\Bbb GO}_p )$ peut \^etre identifi\'e avec l'espace des formes diff\'erentielles $O (2p)$-invariantes.
Soit $\omega$ une forme diff\'erentielle $O (2p)$-invariante sur ${\Bbb GO}_{p}$ et soit
$\bar{\omega}$ une forme diff\'erentielle $O^* (2p)$-invariante sur $X_G$ \'egale \`a $\omega$ au point $x_0$. Puisque
$\omega$ est, en particulier, $\Gamma$-invariante elle induit une forme (n\'ecessairement ferm\'ee) sur $S(\Gamma )$ qui d\'efinit
donc une classe de cohomologie dans $H^* (S(\Gamma ))$. On a ainsi construit une application
\begin{eqnarray} \label{eta3}
\eta : H^* ({\Bbb GO}_{p} ) \rightarrow H^* (S(\Gamma )).
\end{eqnarray}
Il est bien connu que $\eta$ est injective. Enfin, la $i$-\`eme classe de Chern $c_i (T_o )$ du fibr\'e $T_o$ est \'egale \`a $(-1)^i \eta (c_i (\hat{T}_o ))$.

\begin{thm} \label{dec de lef2}
Soit $\Gamma$ un sous-groupe de congruence dans $G$ et soit $\eta : H^* ({\Bbb GO}_p ) \rightarrow H^* (S(\Gamma ))$ l'application d\'efinie en (\ref{eta3}).
Fixons $\lambda$, $\mu$ et $\nu$ trois partitions sym\'etriques incluses dans $p\times p$ telles que le couple $(\lambda, \mu )$ soit compatible. Alors,
\begin{enumerate}
\item pour toute classe fortement primitive $s \in H^{\lambda^- , \mu^-} (S(\Gamma ))$, $\eta (C_{\nu^-} ) .s =0$ si et seulement si la partition $\nu$
ne s'inscrit pas antisym\'etriquement dans le diagramme gauche $\mu / \lambda$, et
\item si $s \in H^{\lambda^- , \mu^-} (S(\Gamma ))$ est une classe non nulle, les \'el\'ements
$C_{\nu^-} .s$, o\`u $\nu^-$ d\'ecrit l'ensemble des parties n\'egatives de partitions sym\'etriques $\nu \subset p\times p$ qui s'inscrivent sym\'etriquement dans le diagramme
gauche $\mu / \lambda  $, sont lin\'eairement ind\'ependants.
\end{enumerate}
\end{thm}

En sp\'ecialisant ce th\'eor\`eme au cas de la cohomologie holomorphe, on corrige tr\`es l\'eg\`erement un th\'eor\`eme de Parthasarathy \cite[Theorem 5.1]{Parthasarathy}.

\begin{cor}[Parthasarathy]
Supposons $p\geq4$. Soient $C_1 , C_2 , \ldots , C_{p-1} \in H^* (S(\Gamma ))$ les classes de Chern du fibr\'e tautologique $T_o$ au-dessus de $S(\Gamma )$. Pour
$j= 3 , 4, \ldots , p$ notons
$$Q_j = \{ s \in H^{\frac{p(p-1)}{2} - \frac{j(j-1)}{2} , 0} (S(\Gamma )) \; : \; C_{j+1} .s = C_{j+2} .s = \ldots = 0 \}.$$
Pour $i=0 , 1 , \ldots , p-1$ notons
$$Q_i ' = \left\{ s \in H^{\frac{p(p-1)}{2} -i , 0} (S(\Gamma )) \; : \; \left( C_3 - C_1 . C_2 \right) . s = 0 \right\} .$$
Alors,
$$H^{l,0} (S(\Gamma )) = \bigoplus_{3 \leq j \leq p} Q_j \oplus \bigoplus_{0 \leq i \leq p-1} Q_i ' .$$
\end{cor}

\subsection{Restriction stable \`a une sous-vari\'et\'e de Shimura}

Comme dans la premi\`ere partie, nous pourrions classifier les diff\'erents types possibles de sous-vari\'et\'es de Shimura. Nos m\'ethodes n'apporteraient pas
de r\'esultats int\'eressants concernant l'injectivit\'e mais en se restreignant \`a la cohomologie holomorphe le lecteur d\'eduira facilement des
m\'ethodes de la premi\`ere partie le th\'eor\`eme suivant qui renforce des r\'esultats ant\'erieurs de Clozel et Venkataramana.

\begin{thm}
Soit $Sh^0 H$ une sous-vari\'et\'e de Shimura de $Sh^0 G$.
Alors, l'application
$${\rm Res}_H^G : H^{*} (Sh^0 G) \rightarrow \prod_{G({\Bbb Q} )} H^{*} (Sh^0 H )$$
de restriction stable est identiquement {\bf nulle} en restriction \`a la cohomologie {\bf holomorphe} $H^{*, 0} (Sh^0 G)$.
\end{thm}

Concluons cet article en remarquant que
le sous-espace $H^{\lambda^- , \mu^-} (Sh^0 G)$ apparait dans la cohomologie holomorphe si et seulement si $\mu^- = (p-1 , p-2 , \ldots , 1)$.
On a alors deux cas possibles.
\begin{enumerate}
\item La partition $\mu = p\times p$. La partition $\lambda$
est alors naturellement param\'etr\'ee par un entier $r$ compris entre $0$ et  $p$ tel que
$$\lambda = ( \underbrace{p, \ldots , p}_{r \ {\rm fois}} , \underbrace{r , \ldots , r}_{p-r \ {\rm fois}} )$$
de diagramme de Young~:
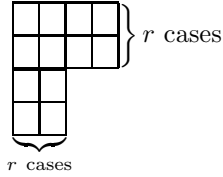
$$
\begin{array}{l}
\left. \hspace{0,035cm}
\begin{array}{|c|c|c|c|} \hline
 & & & \\ \hline
 & & & \\ \hline
\end{array}
\right\}  r \ {\rm cases} \\
\underbrace{
\begin{array}{|c|c|} \hline
 &  \\ \hline
 &  \\ \hline
\end{array}}_{r \ {\rm cases}}
\end{array}
$$
(Ici $p=4$.)

Dans ce cas $H^{\lambda^- , \mu^-} (Sh^0 G) \subset H^{rp- \frac{r(r+1)}{2}, 0} (Sh^0 G)$.

\item La partition $\mu = (\underbrace{p, \ldots , p}_{p-1 \ {\rm fois}} , p-1 )$. La partition $\lambda$ est alors naturellement param\'etr\'ee par un entier $s$ compris entre $0$ et
$p-1$ tel que
$$\lambda = (\underbrace{p, \ldots , p}_{s \ {\rm fois}} , \underbrace{p-1 , \ldots , p-1}_{p-s-1 \ {\rm fois}} , s )$$
de diagramme de Young~:
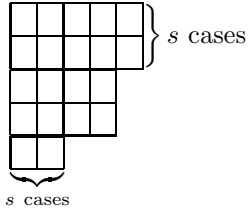
$$
\begin{array}{l}
\left. \hspace{0.035cm}
\begin{array}{|c|c|c|c|c|} \hline
 & & & &  \\ \hline
 & & & & \\ \hline
\end{array}
\right\} s \ {\rm cases} \\
\hspace{0.075 cm}
\begin{array}{|c|c|c|c|} \hline
 & & & \\ \hline
 & & & \\ \hline
\end{array} \\
\underbrace{
\begin{array}{|c|c|} \hline
 & \\ \hline
\end{array}}_{s \ {\rm cases}}
\end{array}
$$
(Ici $p=5$.)

Dans ce cas $H^{\lambda^- , \mu^- } (Sh^0 G) \subset H^{\frac{(p-1)(p-2)}{2} +s ,0} (Sh^0 G)$.
\end{enumerate}

Toute la cohomologie holomorphe de $Sh^0 G$ est obtenue ainsi et les seules redondances viennent du fait que
\begin{itemize}
\item si $r=p$ ou $p-1$ dans le premier cas on obtient le m\^eme sous-espace que si $s= p-1$ dans le deuxi\`eme cas, et
\item si $r=p-2$ dans le premier cas on obtient le m\^eme sous-espace que si $s=p-2$ dans le deuxi\`eme cas.
\end{itemize}

\begin{footnotesize}

\bigskip

\noindent
Unit\'e Mixte de Recherche 8628 du CNRS, \\
Laboratoire de Math\'ematiques, B\^at. 425, \\
Universit\'e Paris-Sud, 91405 Orsay Cedex, France \\
{\it adresse electronique :} \texttt{Nicolas.Bergeron@math.u-psud.fr}

\end{footnotesize}

\end{document}